 \documentclass[a4paper,10pt]{amsart}
 \usepackage{exscale}
 \usepackage{rotate}
 \usepackage{color}
 \usepackage{fancybox}
 \usepackage[dvips]{graphics}
 \newcommand{\N}{{\mathbb N}}
 \newcommand{\C}{{\mathbb C}}
 \newcommand{\Z}{{\mathbb Z}}
 \newcommand{\R}{{\mathbb {R}}}
 \newcommand{\co}{{c_0^2}}
 \newcommand{\un}{{\mathbf{u}}}
 \newcommand{\vn}{{\mathbf{v}}}
 \newcommand{\G}{{\mathbf{G}}}
 \newcommand{\Gn}{{\mathbf{G}}}
 \newcommand{\X}{{\mathbf{X}}}
 \newcommand{\Xn}{{\mathbf{X}}}
 \newcommand{\Tn}{{\mathbf{T}}}
 \newcommand{\nn}{{\mathbf{n}}}
 \newcommand{\bn}{{\mathbf{b}}}
 \newcommand{\Nn}{{\mathbf{N}}}

 \newcommand{\An}{{\mathbf{A}}}
 \newcommand{\Bn}{{\mathbf{B}}}
 
 \newcommand{\A}{{\mathcal{A}}}
 \newcommand{\B}{{\mathcal{B}}}
 \newcommand{\I}{{\mathcal{I}}}
 \newcommand{\rt}{{{\sqrt t}}}
 \newcommand{\tg}{{\tilde\gamma}}
 \newcommand{\twp}{{w_1}}
 \newcommand{\tws}{{w_2}}
 \newcommand{\fdm}{{|f|^2_{+\infty}}} 
 \newcommand{\fcm}{{|f|^4_{+\infty}}} 
 \newcommand{\pp}{{\phi_1(s)}}
 \newcommand{\ps}{{\phi_2(s)}}
 \newcommand{\pt}{{\phi_3(s)}}
 \newtheorem{proposition}{\bf Proposition}
 \newtheorem{corollary}{\bf Corollary}
 \newtheorem{lemma}{\bf Lemma}
 \newtheorem{theorem}{\bf Theorem}
 \theoremstyle{remark}
 \newtheorem{remark}{\bf Remark}
 \title[ Self-similar solutions of the Localized Induced Approximation]
 {Self-Similar solutions of the Localized Induction
 Approximation: Singularity Formation}
 \date{April 01, 2004}
 \email{S.Gutierrez@ed.ac.uk\and mtpvegol@lg.ehu.es }
 \keywords{LIA; Vortex filament, Schr\"odinger equations}
 \thanks{\textit{Mathematics Subject Classification.} 35Q35, 35J10.
 \newline
 L. Vega and S. Guti\'errez were supported by
 DGESIC BFM2001-0458. S. Guti\'errez was also partially funded by the
 Basque Country Government and TMR-Network HARP (HPRN-CT-2001-00273-HARP).
 }
 
\begin{document}
 \maketitle
 \begin{center}
  {\sl{SUSANA GUTIERREZ$^1$ \& LUIS VEGA$^2$}}\\
  {\sl{$^1$University of Edinburgh,
  School of Mathematics\\ 
  JCMB, King's Buildings\\
  Edinburgh, EH9 3JZ, U.K.}}\\[1ex]
  $^2$Departamento de Matem\'aticas, Facultad de Ciencias\\
  Universidad del Pa\'{\i}s Vasco\\
  Aptdo 644, 48080 Bilbao, Spain\\[2ex]
 \end{center}
 \tableofcontents
%
%
 \section{Introduction}
 \label{introduction}
In this paper we continue the investigation started in \cite {GRV}
about the creation of singula\-ri\-ties in a self-similar form for
the {\it {binormal flow}}
 \begin{equation}\label{in1}
 \Xn_t=\Xn_s\times \Xn_{ss},
 \end{equation}
with $\Xn(\cdot\,,t)$ a curve in $\R^3$ parametrized by $s$ the
arclength parameter for all time $t$. This flow was proposed by Da Rios
\cite{Da}, and rediscovered later on by Arms and Hama \cite{AH},
and Betchov \cite{Be}, as an approximation of the evolution of a
vortex tube of infinitesimal cross section $\Xn(\cdot\,,t)$ under
Euler equations. This approximation\footnote{We
refer the reader to \cite{Ba}, \cite{MaBe} and \cite{Saf} for a detailed
analysis of the model.} only takes into consideration the local
effects of the corresponding Biot-Savart integral. For this reason
(\ref{in1}) is also known as the {\it {Localized Induction
Approximation}} (LIA). Notice however that, if use is made of the
Serret-Frenet formulae
\begin{equation}
 \label{star}
 \left\{
 \begin{array}{l}
 {\mathbf T}_s=c {\mathbf n}  \\
 {\mathbf n}_s= -c{\mathbf T} + \tau {\mathbf b}  \\
 {\mathbf b}_s= -\tau {\mathbf n},
 \end{array}
 \right.
\end{equation}
then (\ref{in1}) can be also written as
$$
 \Xn_t=c \bn,
$$
where $c$ stands for the curvature, and $\bn$ for the binormal
vector.

Intimately related to the above flow is the one described by the
tangent vector $\Tn(s,t)$,
$$
 \Tn_t=\Tn\times \Tn_{ss}.
$$
Calling $\Gamma$ the circulation of the vortex tube, the singular
vectorial measure ${\overrightarrow {\omega}}=\Gamma \Tn\,ds$
gives the corresponding vorticity. Therefore, if we want to keep
in the model the fundamental property that the circulation along the
tube be constant, we have to take  $\Tn$ a unit vector-a property 
which is preserved under the flow. This has as a serious
drawback -see~\cite{Saf}-, that if the filament is closed, then 
the total length has to be preserved (at least if a
minimal regularity of the filament is assumed). We wish to emphasise that
$\Tn$ represents the direction of the vorticity vector
${{\overrightarrow \omega}\over{|\overrightarrow \omega|}}$. It
was proved in \cite{CFM} that, among other conditions, the
divergence of the integral
\begin{equation}\label{in2}
 \int_{I} \sup_x|\nabla_x(
 {{\overrightarrow  \omega}\over{ |{\overrightarrow \omega}| } }
 )|^2\,dt
 =\infty
\end{equation}
is necessary for the formation of a singularity in Euler equations for
some time in the interval $I$. In our setting the above integral reduces to
 \begin{equation}\label{in3}
 \int_{I}
 \sup_s|\Tn_s(s,t)|^2\,dt=\int_{I} \sup_s|c(s,t)|^2\,dt=\infty.\
 \end{equation}
In~\cite{GRV} we looked at self-similar solutions of (\ref{in1})
with respect to the unique scaling that preserves arclength. In this
case, it can be easily proved -see~\cite{Bu}, ~\cite{LD}, and 
~\cite{LRT}- that
 \begin{equation}\label{in4}
 \Xn(s,t)=\sqrt t\,\Gn_{c_0}(s/\sqrt t),
 \end{equation}
with $\Gn_{c_0}$ the curve determined by $c(s)=c_0$, $c_0$ denoting a
constant, and $\tau(s)=s/2$. Here $\tau$ stands for the torsion.
As a consequence, $c(s,t)=c_0/\sqrt{t}$ and  (\ref{in3}) 
holds true in this case.

The main purpose of this paper is to characterize all the possible
solutions of (\ref{in1}) such that (\ref{in3}) is verified in a
self-similar way, that is to say such that the curvature is
self-similar with respect to the scaling
$$
 \tilde c(s,t)={{1}\over {\sqrt{t}}} \,c(s/\sqrt t).
$$
Accordingly we take $\tilde\tau(s,t)=\tau(s/\sqrt t)/\sqrt{t}$.
Also recall that from the Hasimoto transformation
-see~\cite{Has}-, if
 \begin{equation}\label{in5}
 \tilde\psi(s,t)=\tilde c(s,t)e^{i\int_0^s\tilde
 \tau(s',t)\,ds'},
 \end{equation}
then
 \begin{equation}\label{in6}
 i\tilde\psi_t+\tilde\psi_{ss}+{\tilde\psi\over 2}
 (|\tilde\psi|^2+\alpha(t))=0.
 \end{equation}
Writing $\tilde\psi(s,t)=\psi(s/\sqrt t)/\sqrt{t}$ and assuming
$\alpha(t)=\alpha/t$ for some constant $\alpha$, we get that
 \begin{equation}\label{in7}
 \psi''
 -{i\over 2}\,
 (\psi+s\psi')+ {{\psi}\over {2}}\, (|\psi|^2+\alpha)=0.
 \end{equation}
Either looking at (\ref{in6}) or at (\ref{in7}) we conclude that 
for any $\alpha$ there is in principle a two parameter family 
of possible solutions.
As we have already seen -see ~\cite{GRV}-, the solutions (\ref{in4}),
except by rotations, are characterized just by one parameter
($c_0$). However, if we consider (\ref{in6}), the initial conditions
that can develop into self-similar solutions have to be homogenous of
degree $-1$, which gives two free parameters. One is the $\delta$
function which on integrating  the Frenet system of
equations leads to the kink solution of (\ref{in1}) obtained
in~\cite{GRV}. The other candidate is to take  
$\tilde\psi(s,0)=p.v.(1/s)$, as initial condition
of (\ref{in6}). In this case,
when we perform the first integration in $s$, a logarithmic term
appears which breaks the scaling symmetry. Therefore, in order to
see these ``self-similar solutions" in the framework of
(\ref{in1}) we have to make a modification of the usual ansantz
given in (\ref{in4}).

Take $\A$ a real antisymmetric $3\times 3$ matrix and define for
some $\Gn$
 \begin{equation}\label{in8}
 \Xn(s,t)=e^{{{\A}\over {2}}\log t}\sqrt t\Gn(s/\sqrt t).
 \end{equation}
If we ask that $\Xn(s,t)$  be a solution of (\ref{in1}), we get that
$\Gn$ has to solve the system of ordinary differential equations
(O.D.E. for short)
\begin{equation}\label{in9}
(\I+\A) \Gn-s\Gn'=2\Gn'\times \Gn'', \qquad |\Gn'(s)|^2=1.
\end{equation}
Multiplication by $\Gn'$ and from the vectorial identity
$$
 {\mathbf F}\times ({\mathbf G}\times {\mathbf H})=
 ({\mathbf F}\cdot {\mathbf H})\, {\mathbf G}-
 ({\mathbf F}\cdot {\mathbf G})\, {\mathbf H},
$$
we get
\begin{equation}\label{in10}
    \Gn''={1 \over 2}(\I + \A)\Gn\times \Gn'.
\end{equation}
Given any initial condition $(\Gn(0), \Gn'(0))$ with
$|\Gn'(0)|^2=1$, it is very easy to prove global existence of
a solution $\Gn$ which solves the above equation 
(see section~\ref{main}).

In section~\ref{main},  associated with the solution $\Gn$, we
shall prove the existence of regular $\Xn(s,t)$ of the form
(\ref{in8}) solving the following problem:
\begin{eqnarray} \label{in11}
 \left\{
 \begin{array}{l}
 \X_t= \X_s\times \X_{ss},\qquad t>0,\quad s\in\R, \\[2ex]
 \X(s,0)=  se^{\A\log|s|}{\mathbf A}^{+}\chi_{[0 ,+\infty)}(s)
        + se^{\A\log|s|}{\mathbf A}^{-}\chi_{(-\infty, 0]}(s) \qquad ,
 \end{array}
 \right.
\end{eqnarray}
for some $({\mathbf A}^{+}$, ${\mathbf A}^{-})\in \R^3\times
\R^3$. Here, $\chi_{E}$ denotes the characteristic function of the
set $E$ and, due to the invariance of LIA under rotations, we can
assume without loss of generality that
\begin{equation}
 \label{in12}
  \A=
 \left(
 \begin{array}{ccc}
  0 & -a & 0 \\
  a & 0 & 0  \\
  0 & 0 & 0
 \end{array}
 \right),
 \quad a\in\R.
\end{equation}

The initial data $\Xn(s,0)$ in (\ref{in11}) includes a
wide variety of 3d-spirals, whose rotation axis is the OZ-axis under
the assumption (\ref{in12}) on the matrix $\A$.
The initial data $\Xn(s,0)$ is singular at $s=0$ whenever, either $a=0$ and
$\An^{+}-\An^{-}\neq {\mathbf {0}}$, or $a\neq 0$ and the first two 
components of either $\An^{+}$ or $\An^{-}$ are different from zero. The
singularity of $\Xn(s,0)$ becames clear on studying the behaviour of
the tangent vector $\Tn(s,0)$ at $s=0$: in the case when $a\neq
0$, the singularity at $s=0$ comes from the non-existence 
the limit $\lim_{s\rightarrow 0} \Tn(s,0)$, whereas if
$a=0$ we find that $\Tn(s,0)$ has a jump singularity at $s=0$ (see~\cite{GRV}).

It is important to mention that, because
$\X(s,t)$ is of the form (\ref{in8}), the
shape of the curve $\Xn(s,0)$ is directly related to the
asymptotic behaviour of the solutions of (\ref{in10}). In this
setting we have the following theorem:
\medskip

\begin{theorem}\label{mainG1}
Given $(\Gn(0),\Gn'(0))$ with $| \Gn'(0)|^2=1$, let $\Gn(s)$ be
the solution of (\ref{in10}) asso\-cia\-ted to this initial data.
Then, there exist unique vectors ${\Bn^{\pm}}$ with
$|\Bn^{\pm}|=1$ and $\An^{\pm}=(\I+\A)^{-1}\Bn^{\pm}$ such that
the following asymptotics hold as $s\rightarrow \pm\infty$:
\begin{eqnarray*}
 & &
 \noindent {\it {i)}}
 \quad
 \G(s)=
 s\, e^{\A\log{|s|}}\An^{\pm}+
 {{e^{\A\log{|s|}}}\over s}
 \{ 2c^2_{\pm\infty}\Bn^{\pm}-\A \Bn^{\pm}\times\Bn^{\pm}\}
 - 4 \,{{c_{\pm \infty}\nn}\over{s^2}}
 +O\left({1\over {s^3}}\right),
    \\[2ex]
 & &
 \noindent {\it {ii)}}
 \quad
 \Tn(s)=
 e^{\A\log{|s|}}\Bn^{\pm}-2{{c_{\pm \infty}\bn}\over {s}}
 + O\left( {1\over s^2}\right).
      \\[2ex]
 & &
 \noindent {\it {iii)}}
 {\hbox { Moreover, if }} a\neq 0, B_3^{\pm}\neq \pm 1,
 {\hbox { and }} c_{\pm\infty}\neq 0, {\hbox { then}}
     \\[2ex]
 & &
 \hspace{0.8cm}
 c_{\pm \infty}(\nn-i\bn)(s)=
 \,{{b_\pm e^{ia_\pm}}\over {|\A\Bn^{\pm}|^2}}\,
 e^{i\phi(s)}
 e^{\A\log{|s|}}
 \{ \A \Bn^{\pm}\times \Bn^{\pm}
 -i\A\Bn^{\pm}  \}+O\left({1\over|s|}\right).   
     \\[1ex]
 & &
 \noindent {\it {iv)}}
 {\hbox { If }} a\neq 0  {\hbox { and  either }} B_3^{+} 
 {\hbox { or }} B_3^{-}\in\{+1,-1 \},
 {\hbox { then }} \Gn(s)=(0,0,\pm s).
\end{eqnarray*}

 \noindent Here\footnote{$\{ \Tn, \nn, \bn\}$ 
 is the Frenet frame associated
 to $\Gn$, $c(s)$ the curvature function, and 
 $c_{\pm\infty}=\lim_{s\rightarrow \pm\infty.}c(s)$. The existence
 of the limits $c_{\pm\infty}$ will be proved in
 Corollary~\ref{cofor}},
 $\phi(s)=(s^2/ 4)-\gamma_{\pm}\log |s|$, $a_{\pm}\in[0,2\pi)$,
 $$
 \begin{array}{l}
  \gamma_{\pm}= 3a B_3^{\pm}+\alpha,
  {\hskip0.5cm}
  c^2_{\pm\infty}= -a B_3^{\pm}-\alpha,
             \\[1ex]
  b^2_{\pm}= a^2 (-aB_3^{\pm}-\alpha)\,
          (1- (B_3^{\pm})^2)
  \quad {\hbox {with}}\quad
   b_{\pm}\geq 0
  \qquad {\hbox {and}}
             \\[1ex]
  \alpha= -aT_3(0)-{1\over 4} |(\I+\A)\Gn(0)|^2.
 \end{array}
 $$
 \end{theorem}
\begin{remark}
 \label{remainG1}
 In {\it (i)} and {\it (ii)} we understand that $c_{\pm\infty}\nn=
 c_{\pm\infty}\bn={\mathbf 0}$, whenever $c_{\pm\infty}=0$.
\end{remark}

Conversely, we can fix the data at infinity and find 
$(\Gn(0), \Gn'(0))$ such that the asymptotics at
infinity of the corresponding solution $\Gn$ is prescribed 
by the given data. More precisely:
\medskip

\begin{theorem}\label{mainG2}
Given $a\neq 0$, $\Bn^+=(B_1^+,B_2^+,B_3^+) \in {\mathbb {S}}^2$
with $B_3^+\neq\pm1$, $a_+\in[0,2\pi)$ and $b_+\geq0$, there
exists a solution of (\ref{in9}) satisfying
$(i)$, $(ii)$ and $(iii)$ of Theorem~~\ref{mainG1}, if
$s\to+\infty$ and $\Bn^+ = (\I+\A)\An^+$. Moreover, if $b_{+}>0$, 
the solution is unique.
A similar result can be obtained at $s=-\infty$.
\end{theorem}
\medskip

The proofs of Theorems~\ref{mainG1}  and~\ref{mainG2} are related to
the asymptotic behaviour of the solutions of the complex
O.D.E. (\ref{in7}). In order to integrate (\ref{in7}), it
turns out that, it is better to introduce the variable $f(s)$
defined by
\begin{equation}
 \label{in12'}
   \psi(s)= f(s) e^{is^2/4},
\end{equation}
and consider the equation
$$
  f'' + i {s\over 2} f' +{f\over 2} (|f|^2+\alpha)=0.
$$
The theorem below gathers the main properties that we prove for
the solutions of the latter equation:
\begin{theorem}\label{mainf1}\
Let $f$ be a solution of the equation
\begin{equation}
 \label{in13}
 f''+i{s\over 2}f'+{f\over 2}(|f|^2+\alpha)=0,
 \qquad \alpha\in \R.
\end{equation}
Then
\begin{itemize}

\item[{\it i)}] There exists $E(0)\geq 0$ such that the identity
 $$
 |f'|^2+{1\over 4}(|f|^2+\alpha)^2=E(0)
 $$
 holds true for all $s\in\R$.
\item[{\it ii)}] The limits
$\lim_{s\rightarrow \pm\infty}|f|^2(s)=|f|^2_{\pm\infty}$ and
$\lim_{s\rightarrow \pm\infty}|f'|^2(s)=|f'|^2_{\pm\infty}$ do
exist.
\item[{\it iii)}] Moreover, if $|f|_{+\infty}\neq 0$ or
 $|f|_{-\infty}\neq 0$, then
 \begin{eqnarray*}
  f(s)
  &=&
  |f|_{\pm\infty}\, e^{ic_\pm}\, e^{i\ps}
  + 2i\, |f'|_{\pm\infty} \,
  {{e^{id_\pm}}\over {s}}\, e^{i\pt}
  + O\left( {1\over {|s|^2}} \right),
       \\
  f'(s)
  &=&
  |f'|_{\pm\infty} e^{id_{\pm}} e^{i\phi_3(s)}
  + i |f|_{\pm\infty}(|f|^2_{\pm\infty}+\alpha)\,
  {{e^{ic_\pm}}\over {s}}\, e^{i\phi_2(s)}
  +O\left( {{1}\over {|s|^2}} \right).
 \end{eqnarray*}
\end{itemize}
 Here,  $|f|_{\pm\infty}$, $|f'|_{\pm\infty}\geq 0$,
 $c_{\pm}$ and $\ d_{\pm}$ are arbitrary
 constants in $[0,2\pi)$,
 $$
  \ps= (|f|^2_{\pm\infty}+\alpha)\log |s|,
  \qquad {\hbox {and}}\qquad
  \pt= -(s^2/ 4)- (2|f|^2_{\pm\infty}+\alpha)\log |s|.
 $$
\end{theorem}
Moreover, we also prove the following converse of the above
theorem:
\begin{theorem}\label{mainf2}
\noindent With the same notation as in Theorem~\ref{mainf1}, given
complex numbers  $|f|_{+\infty}e^{i\theta_1}$ and
$|f'|_{+\infty}e^{i\theta_2}$, with $|f|_{+\infty}$, 
$|f'|_{+\infty}\geq 0$ and $\theta_1$, $\theta_2\in [0,2\pi)$, 
there exists a unique solution $f \in {\mathcal {\C}}^2(\R)$ 
of (\ref{in13}) satisfying
$(iii)$ of Theorem~\ref{mainf1}, if $s\to+\infty$. A similar
result can be obtained at $s\to -\infty$.
\end{theorem}
The asymptotic behaviour for the solutions of (\ref{in13}) given
in Theorem~\ref{mainf1} is not only used in the proof of
Theorems~\ref{mainG1} and \ref{mainG2}, but also it is used to prove
an ill-posedness result for the following initial value problem 
(IVP for short) related to non-linear cubic
Schr\"odinger equations:
\begin{eqnarray}
  \label{in14}
  \left\{
  \begin{array}{l}
   i\psi_t+\psi_{ss} +
   \displaystyle{
   {\psi\over 2}\,(|\psi|^2+{\alpha\over t})}
   =0,\qquad t>0, \quad s\in \R,\qquad
   \alpha\geq 0
        \\[2ex]
   \psi(s,0)= c_1 \,p.v. (1/s),
   \qquad {\hbox {with}}\qquad c_1\in\C \setminus\{0\}.
  \end{array}
  \right.
\end{eqnarray}

 This paper is organized as follows. Section~\ref{main} contains
 the results related to solutions of the binormal flow. For the
 reader's convenience, we split Section~\ref{main} into different
 subsections. In Subsection~\ref{initialdata}, we prove
 the existence of solutions of LIA that converge uniformly to an
 initial data in the shape of a 3d-spiral. The exhaustive
 asymptotics for $\Gn$ and $\Tn$ in Theorem~\ref{mainG1} will be
 proved in Subsection~\ref{refined}. Finally, the proof of
 Theorem~\ref{mainG2} is included in Subsection~\ref{scatteringG}.

The key point for the asymptotics of $\Tn$ and $c(\nn -i \bn)$ 
is to look at the quantities
 \begin{equation}\label{in15}
   y= \displaystyle{
   {{d|\Tn'|^2}\over {ds}}}
   \qquad {\hbox {and}}\qquad
   h=-{1\over 2}\, \A\Tn\cdot \Tn'.
 \end{equation}
 In fact, $(y,h)$ solves the system -see (\ref{cl9}) and Remark~\ref{reyh}-
 \begin{eqnarray}\label{in16}
  \left\{
  \begin{array}{ll}
   y'=sh+g(|\Tn'|^2);\qquad 
     g(|\Tn'|^2)=2E(0)-(3|\Tn'|^2+\alpha)(|\Tn'|^2+\alpha)/2,
     \\[2ex]
   h'=-{\displaystyle{{s\over 4}}}\,y,
  \end{array}
  \right.
 \end{eqnarray}
 where $E(0)=a^2/4$.

 In the case $c(s)\neq 0$, then 
 \begin{equation}\label{in15'}
   h=c^2(\tau-s/2), \qquad y={{dc^2}\over {ds}},
   \qquad {\hbox {and}}\qquad
   {y\over 2}+ih=\bar f f',
 \end{equation}
 where $f$ is given in (\ref{in12'}) and is a solution of
 (\ref{in13}). In fact, the equation (\ref{in13}) will be reduced to
 solve the system (\ref{in16}). The study of the 
 properties of the solutions $(y,h)$ of
 the latter system of equations, 
 together with the proofs of Theorem~\ref{mainf1} and 
 its partial converse Theorem~\ref{mainf2}
 are saved for Section~\ref{selfNLS}.  

 We conclude this paper with Section~\ref{fr}. This section is
 devoted to proving more specific facts and consequences on the
 results obtained in the previous sections.
 Here, we consider two special symmetric cases of solutions of
 LIA. Also we discuss the question of ill-posedness for the IVP 
 (\ref{in14}) and the binormal flow. In particular, in
 Proposition~\ref{uT1} we prove the lack of uniqueness of weak
 solutions for the IVP associated to LIA when the
 initial data considered is a curve in the shape of a corner.

 We finish this section giving some of the notation that will be
 used through this paper. In the sequel  $f'$, $f_s$ or ${df/ds}$
 will denote the derivative with respect to the variable $s$ of
 $f(s)$. Unless it is explicitly stated otherwise,
 we will use bold and gothic letters to denote vectors and
 matrices, respectively, and the overbar will indicate the complex
 conjugate. Finally, ${\mathbb{S}}^{n-1}$ will be the unit sphere
 in the euclidean space $\R^n$.

 \noindent {\bf {Acknowledgements.}} The authors would like to
 thank J. Rivas for her tireless help with the figures of this work.
 We also want to thank R.~Jerrard for enlightening conversations.
%
\section{Proofs of the Theorems}
\label{main}
%
Given $a\in\R$, define for some $\Gn$
\begin{equation}\label{11}
 \X_{a}(s,t)= {e}^{{\A\over 2}\log t} \sqrt{t}\G({s/{\sqrt{t}}}),
 \qquad t>0,
\end{equation}
where, as we have already said in the introduction, $\A$ is
assumed to be, with out loss of generality, the matrix
\begin{equation}
 \label{10}
  \A=
 \left(
 \begin{array}{ccc}
  0 & -a & 0 \\
  a & 0 & 0  \\
  0 & 0 & 0
 \end{array}
 \right),
 \quad a\in\R.
\end{equation}
It is easy to see that $\Xn_a(s,t)$ is a solution of LIA if and
only if $\Gn(s)$ satisfies
\begin{equation}
 \label{12}
 (\I+\A)\G- s\G'= 2 \G'\times \G'', \qquad |\G'|=1.
\end{equation}
Observe that the equation (\ref{12}) can be written 
equivalently as
\begin{equation}
 \label{13}
 \G''={1\over 2} (\I+\A)\G\times \G',
\end{equation}
whenever
\begin{equation}
 \label{13'}
 |\G'(0)|=1 \qquad {\hbox {and}}\qquad
 (\I+\A)\G(0)\cdot \G'(0)=0.
\end{equation}
Indeed, the equation (\ref{13}) easily follows from (\ref{12}) by
taking the outer product of (\ref{12}) and $\G'$, and using the
vectorial identity
\begin{equation}
 \label{14}
 {\mathbf F}\times ({\mathbf G}\times {\mathbf H})=
 ({\mathbf F}\cdot {\mathbf H})\, {\mathbf G}-
 ({\mathbf F}\cdot {\mathbf G})\, {\mathbf H}.
\end{equation}
Now, assume that $\G$ is a solution of (\ref{13}) satisfying the
initial conditions (\ref{13'}). Then, we first notice that
\begin{equation}
 \label{14a}
 {{d }\over {ds}}|\G'(s)|^2= 2\G''\cdot \G' =
 ((\I+\A)\G\times\G')\cdot \G'=0,
\end{equation}
so that
\begin{equation}
 \label{14c}
 |\G'(s)|^2=|\G'(0)|^2=1,\ \forall\, s\in \R.
\end{equation}
Secondly, the outer product of $\G'$ and (\ref{13}) together with
the above identity  yields
$$
 \G'\times \G''=
 {1\over 2}((\I+\A) \G -(\G'\cdot (\I+\A)\G)\,\G').
$$
As a consequence, it is enough to prove that $\G'\cdot(I+\A)\G=s$.
To this end, notice that
\begin{equation}
 \label{14b}
 {{d }\over {ds}} (\G'\cdot (\I+\A)\G)=
 \G''\cdot (\I+\A)\G+ \G'\cdot(\I+\A)\G'=1,
\end{equation}
and by taking into account that the outer product of (\ref{13}) and
$\G''$ gives that $\G''\cdot (\I+\A)\G=0$, and that 
$\A\vn\cdot \vn=0$. Here we have used the antisymmetry of $\A$.
Integrating (\ref{14b}) and using the initial conditions in
(\ref{13'}) allow us to conclude that $\G'\cdot (\I+\A)\G=s$.

The previous observation reduces the problem of finding solutions of
LIA of the form (\ref{11}) to proving the existence of solutions of
(\ref{13}) with initial conditions satisfying (\ref{13'}).
In this setting, notice that the local existence of solution of
the initial value problem (\ref{13})-(\ref{13'})
follows from the classical theory for first order O.D.E. system.
The global existence of ${\mathcal {C}}^{\infty}(\R;
\R^3)$-solution follows from (\ref{14c}).

The proposition below summarizes the obtained results.
\begin{proposition}
 \label{solucion}
 Given $a\in\R$, define
 $$
  \Xn_a(s,t)=e^{{\A\over 2}\log t} {\sqrt t}\, \G(s/{\sqrt
  t}),
  \quad t>0,
  \qquad {\hbox {with}}\qquad
  \A=
  \left(
  \begin{array}{ccc}
   0 & -a & 0 \\
   a & 0 & 0  \\
   0 & 0 & 0
  \end{array}
  \right),
 $$
where $\G(s)$ is the solution of
$$
  \G''= {1\over 2} (\I+\A) \G\times \G'  \\[2ex]
$$
associated to a given initial data $(\Gn(0),\Gn'(0))$ such that
$$
 |\G'(0)|=1,
 \qquad {\hbox {and}}\qquad
 (\I+\A)\G(0)\cdot \G'(0)=0.
$$
Then, $\X_a(s,t)$ is an analytic solution of LIA for all $t>0$, and
$c^2(s)$ is also analytic.
\end{proposition}
\begin{remark}\label{sol}
 The condition $(\I+\A)\Gn(0)\cdot \Gn'(0)=0$
 can be removed:
 Let $\G$ be a solution of (\ref{13}) associated to a given
 initial data $(\Gn(0), \Gn'(0))\in \R^3\times {\mathbb {S}}$, and
 consider $\tilde \Gn(s)=\Gn(s-s_0)$ with
 $s_0= (\I+\A)\Gn(0)\cdot\Gn'(0)$. Then,
 $\tilde\Gn$ satisfies (\ref{12}) and
 the result in Proposition~\ref{solucion} holds true for
 $\Xn_{a}(s,t)= e^{{\A\over 2}\log t} {\sqrt{t}}\tilde
 \Gn(s/{\sqrt{t}})$.

 Hereafter, we will assume without loss of generality that
 $s_0=0$, i.e. $(\I+\A)\Gn(0)\cdot \Gn'(0)=0$.
 In the case when $s_0\neq 0$, the results that will be  proved
 in the sequel for $\Gn$ will be valid for $\tilde\Gn$.
\end{remark}
 \begin{figure}[h]
  \label{otros}
  \begin{center}
   \scalebox{0.5}{\includegraphics{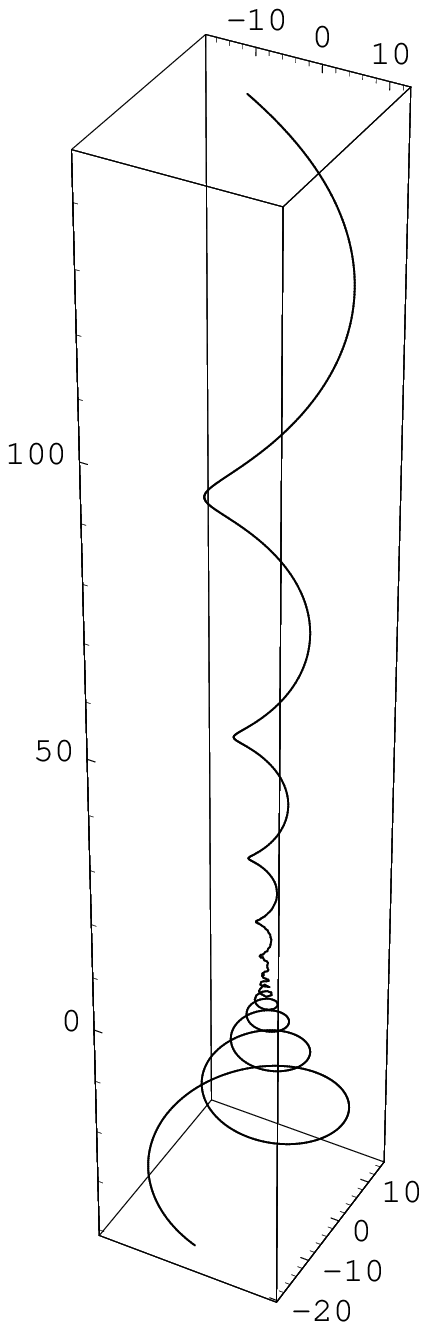}}
   \hspace{2cm} 
   \scalebox{0.5}{\includegraphics{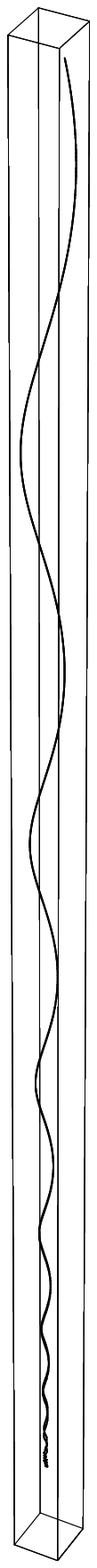}}
   \hspace{2cm}
   \scalebox{0.5}{\includegraphics{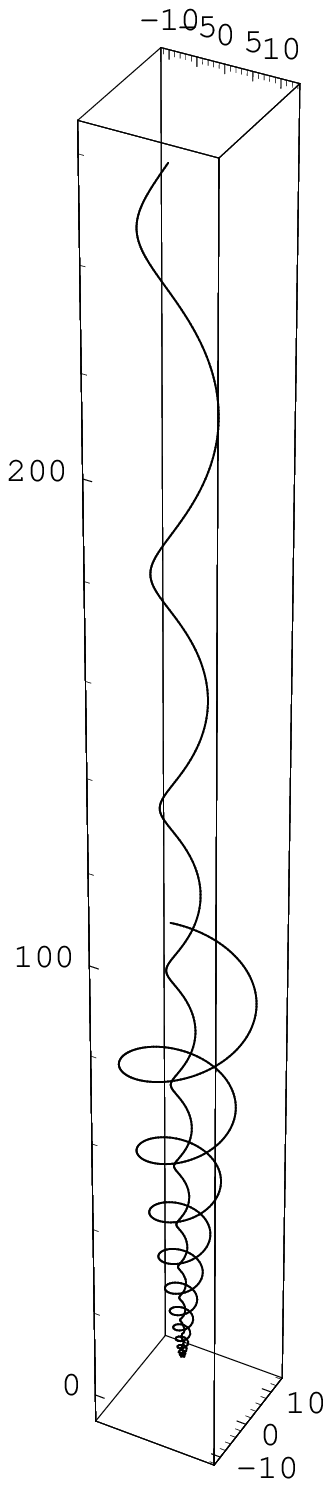}}
   \caption{$\Gn(0)=(0,0,2c_0)$, $\Gn'(0)=(1,0,0)$, with 
           $(a,c_0)=(10,1)$, $(15,5)$ and $(20,3)$, respectively.}
  \end{center}
 \end{figure}
Figure~1 depicts different examples of the curve $\Xn(s,t)$
at time $t=1$, that is of $\Gn(s)$. The curves $\Gn(s)$ have been
obtained solving the equation (\ref{13}) with the initial conditions
$\G(0)=(0,0,2c_0)$ and $\Gn'(0)=(1,0,0)$. The subsequent plots correspond
with different choices of the parameters $c_0$ and $a$. 
As we will continue to prove, the evolution of each of these curves
$\G(s)$ under the relation in (\ref{11}) leads to solutions
of LIA which converge as $t\rightarrow 0^{+}$ to different curves 
in the shape of a $3d$-spiral.

In what follows, we will often drop the subscript $a$ in the
definition of $\Xn_{a}(s,t)$ in Proposition~\ref{solucion} and
just write $\Xn(s,t)$.
Hereafter, $c(s)$ and $\tau(s)$ will be the curvature and torsion
functions related to $\Xn(s,1)=\G(s)$, respectively. Also, we will
write $c_0= c(0)$ and $c_{\pm\infty}=\lim_{s\rightarrow
\pm\infty}c(s)$. 
%
\subsection{Convergence to the initial data}\
\label{initialdata}
\smallskip

\noindent As in~\cite{GRV}, the study of $\lim_{t\rightarrow
0^+}\Xn(s,t)$ is based on finding an explicit expression of
$\Gn(s)$. Here, we also obtain a formula of the tangent vector
$\Tn(s)$ which will be used in the following subsection.

For simplicity of the exposition, through this subsection we will
write $\Gn'\times \Gn''= c\bn$ even when $c(s)=0$, in which case
the both sides of this identity will be understood to be zero.

We begin obtaining a closed formula for $\G(s)$. Consider the
quantity
$$
  e^{-\A\log{|s|}} {\G\over s},\qquad s\neq 0.
$$
Then, deriving with respect to $s$ and taking into account that
$\G$ is a solution of (\ref{12}), we obtain
\begin{eqnarray}\label{15}
 \left( e^{-\A\log{|s|}}{\G\over s}\right)_s
 &=&
 - e^{-\A\log{|s|}} {{\A\G}\over {s^2}}+ e^{-\A\log{|s|}}
 \left( {{s\G_s-\G}\over {s^2}}  \right)
    \nonumber \\
 &=&
 -2e^{-\A\log{|s|}}{{c\bn}\over {s^2}},
\end{eqnarray}
and therefore, integrating the above identity in the interval
$[s_1,s]$, we get
 \begin{equation}\label{15'}
 e^{-\A\log{|s|}} {{\G(s)}\over {s}} -
 e^{-\A\log{|s_1|}} {{\G(s_1)}\over {s_1}}=
 -2\int_{s_1}^{s} e^{-\A\log{|s'|}} {{c\bn}\over {(s')^2}}\, ds'.
 \end{equation}
We need the following lemma.
\begin{lemma} \label{for}
 If $\Gn$ is a  solution of (\ref{12}), then
 \begin{eqnarray}
 \label{formula}
 & &
 |\Tn'|^2(s)=
 -aT_3(s)-\alpha,
 \qquad {\hbox {with}} \qquad
 \alpha=-a T_3(0)- {1\over 4}\, |(\I+\A)\Gn(0)|^2
   \\[1ex]
 & &
 \qquad\qquad{\hbox {or, equivalently, }}
    \nonumber\\
 & &
 c^2(s)=
 -a T_3(s)-\alpha,
 \quad {\hbox {with}}\quad
 \alpha= -a T_3(0)-\co
 \quad {\hbox {and}}\quad
 \co= |(\I+\A)\Gn(0)|^2/4.
   \nonumber
 \end{eqnarray}
\end{lemma}
 \begin{remark}
 \label{reformula}
 If $\alpha>0$ and $a>0$ (resp. $a<0$), then from (\ref{formula}) it
 follows that  $G_3(s)$ is monotone decreasing (increasing), and 
 therefore $\Gn(s)$ has no self-intersections (see Figure~3).
 \end{remark}
\begin{proof}[Proof of Lemma \ref{for}]
 Firstly, notice that (\ref{12}) rewrittes as
 \begin{equation}
 \label{lemma1}
 (\I+\A)\Gn= s\Tn +2\, \Tn\times \Tn',
 \qquad {\hbox {and}}\qquad
 |\Tn|^2=1.
 \end{equation}
Deriving (\ref{lemma1}), we obtain 
$$
 \A\Tn= s\Tn' + 2\Tn\times \Tn'',\qquad \Tn\cdot \Tn'=0.
$$
Hence,
\begin{eqnarray*}
 \A\Tn\times \Tn
 &=&
 -s\Tn\times \Tn' - 2 (\Tn\cdot \Tn'') \Tn + 2 |\Tn|^2 \Tn''
    \\
 &=&
  - {s\over 2} (\I+\A)\Gn + {s^2\over 2}\, \Tn
  - 2 (\Tn\cdot \Tn'') \Tn + 2 \Tn'',
\end{eqnarray*}
by bearing in mind the vectorial identity (\ref{14}) and
(\ref{lemma1}).

Next, from (\ref{lemma1}), $(\I+\A)\Gn\cdot \Tn'=0$, so that if we
take the scalar product of the above identity and $\Tn'$, we
obtain 
$$
  (\A\Tn\times \Tn) \cdot \Tn'= {{d|\Tn'|^2}\over {ds}}.
$$
Finally, a simple calculation proves that, if $\A$ is defined by
(\ref{10}), then the latter identity rewrites
$$
   {{d|\Tn'|^2}\over {ds}}= -a T'_3(s).
$$
The result follows by integrating the previous identity and
noticing that, from (\ref{12}), $|\Tn'(0)|^2= |(\Tn'\times \Tn
)(0)|^2= |(\I+\A)\Gn(0)|^2/4$.
\end{proof}
Notice that, from the previous lemma, in particular  we get that
$|c(s)|\leq C,\ \, \forall s\in\R$, and $|a|\leq M$.

Then, if we choose $s_1=\pm 1$ in (\ref{15'}) and let
$s\rightarrow \pm \infty$, we conclude that the limits
$\lim_{s\rightarrow\pm \infty} e^{-\A\log{|s|}}\G(s)/s$ do exist
and are defined by
\begin{eqnarray} \label{15''}
 \An^{+}
 &=&
 \lim_{s\rightarrow +\infty}e^{-\A\log{|s|}} {{\G(s)}\over {s}}=
 \G(1)-2\int_{1}^{+\infty} e^{-\A\log{|s'|}} {{c\bn}\over {(s')^2}}\, ds',
       \nonumber   \\
 & &
          \\
 \An^{-}
 &=&
 \lim_{s\rightarrow -\infty}e^{-\A\log{|s|}} {{\G(s)}\over {s}}=
 -\G(-1)+2\int_{-\infty}^{-1} e^{-\A\log{|s'|}} {{c\bn}\over {(s')^2}}\, ds'.
        \nonumber
\end{eqnarray}
Now, choosing in (\ref{15'}) $s=+\infty, \ s_1=s>0$ and
$s=-\infty,\ s_1=s<0$, we respectively obtain 
\begin{eqnarray}\label{16}
 \G(s)
 &=&
 se^{\A\log{|s|}}\An^{+} +
 2se^{\A\log{|s|}}\int_{s}^{+\infty} e^{-\A\log{|s'|}} {{c\bn}\over
 {(s')^2}}\, ds', \qquad s>0,
     \nonumber \\
 & &
               \\
 \G(s)
 &=&
 se^{\A\log{|s|}}\An^{-} -
 2se^{\A\log{|s|}}\int_{-\infty}^{s} e^{-\A\log{|s'|}} {{c\bn}\over
 {(s')^2}}\, ds', \qquad s<0.
     \nonumber
\end{eqnarray}
Let us define
$$
 \Bn^{\pm}=(\I+\A)\An^{\pm}=
 \lim_{s\rightarrow \pm\infty} e^{-\A\log|s|}\Tn(s),
$$
(notice that, because $|\Tn(s)|=1$, from the above identity easily
follows that $|\Bn^{\pm}|=1$).

Then, deriving (\ref{16}) and using the above notation, we obtain
the following expressions for the tangent vector
 \begin{eqnarray}\label{17}
 \Tn(s)
 &=&
 e^{\A\log{|s|}}\Bn^{+}- 2{{c\bn}\over {s}}+
 2e^{\A\log{|s|}}\int_{s}^{+\infty}
 e^{-\A\log{|s'|}}{{(\I+\A)(c\bn)}\over {(s')^2}}\,ds',
 \quad s>0
  \nonumber \\
 & &
            \\
 \Tn(s)
 &=&
 e^{\A\log{|s|}}\Bn^{-}- 2{{c\bn}\over {s}}-
 2e^{\A\log{|s|}}\int_{-\infty}^{s}
 e^{-\A\log{|s'|}}{{(\I+\A)(c\bn)}\over {(s')^2}}\,ds',
 \quad s<0.
  \nonumber
 \end{eqnarray}
 \begin{corollary}
  \label{cofor}
  The limits $\lim_{s\rightarrow \pm\infty}T_3(s)$ do exist, and
  therefore (see (\ref{formula})) it follows the existence of
  $c_{\pm\infty}=lim_{s\rightarrow \pm\infty}c(s)$.
 \end{corollary}
We will now continue with the proof of the convergence of
$\Xn(s,t)$ to the initial data. Precisely, we will prove the
following result:
\begin{proposition}
 \label{convergencia}
 Let $\Xn(s,t)$ be the solution of LIA defined in
 Proposition~\ref{solucion}. Then,
 there exist ${\mathbf A}^{+}$ and ${\mathbf A}^{-}\in\R^3$ (not
 necessarily unitary vectors) such that
 $$
 \lim_{t\rightarrow 0^{+}} \Xn_{a}(s,t)=
 \lim_{t\rightarrow 0^{+}}
 e^{{\A\over 2}\log t}{\sqrt {t}}{\G}_{a}(s/{\sqrt t})=
 \left\{
  \begin{array}{ll}
   se^{\A \log |s|}{\mathbf A}^{+}, & s\geq 0 \\[2ex]
   se^{\A \log |s|}{\mathbf A}^{-}, & s\leq 0
  \end{array}
 \right.
 $$
 with
 $$
   |\Xn_{a}(s,t)-se^{\A\log|s|}\An^{\pm}|\leq 2\sqrt{t}\,
  \left(
  \underset{s\in\R}{sup} \ |c(s)|
  \right).
 $$
 Moreover, the maps
 $T^{\pm}:\ (\Gn(0),\Gn'(0),a)\in\R^3\times {\mathbb{S}}^2\times\R
 \longrightarrow \An^{\pm}
 $ are continuous.
 \end{proposition}
\begin{proof}
 From the definition of the matrix $\A$ in (\ref{10}), it can be
 easily shown that $e^{\A\log(|s|/{\sqrt t})}$ is an orthogonal
 matrix $\forall\, (s,t)$ such that $|s|/{\sqrt t}\neq 0$. Also, the norm
 $\| (\I+\A)^{-1} \|\leq 1$.

 Now, let $s>0$. By taking into account the just mentioned
 properties related to $\A$, from the expression of $\G(s)$ in
 (\ref{16}), we get that
 \begin{eqnarray*}
  |\Xn(s,t)- se^{\A\log|s|}{\mathbf A}^{+}|
  &=&
  | e^{{\A\over 2}\log t} {\sqrt t} \G(s/{\sqrt t})-
    se^{\A\log|s|}{\mathbf A}^{+}|
      \\
  &=&
  \left|
    2\, s e^{\A\log s}\int_{s/{\sqrt t}}^{+\infty}
    e^{-\A\log|s'|} {{(c\bn)(s')}\over {(s')^2}}\, ds'
  \right|
     \\
  &=&
  \left|
    2\, s \int_{s/{\sqrt t}}^{+\infty}
    e^{-\A\log|s'|} {{(c\bn)(s')}\over {(s')^2}}\, ds'
  \right|
      \\
  &\leq&
   2s{\sqrt t}\int_{s}^{+\infty}
   {{|c(s''/{\sqrt t})|}\over {(s'')^2}}\, ds''
      \\
  &\leq&
  2\, {\sqrt t}
  \left(
  \underset{s>0}{sup} \ |c(s)|
  \right).
 \end{eqnarray*}
The same argument as above, where we now use (\ref{16}) for $s<0$,
gives that
 $$
 |\Xn(s,t)- se^{\A\log|s|}{\mathbf A}^{-}|
 \leq
 2\, {\sqrt t}
 \left(
 \underset{s>0}{sup}\  |c(s)|
 \right).
 $$
 Next, consider $s=0$. Notice that, because
 $\Xn(s,t)=e^{{\A\over 2}\log t}{\sqrt t}\G(s/{\sqrt t})$, we know
 that $\G(s)$ is solution of (\ref{12}). In particular, the
 evaluation of (\ref{12}) at $s=0$ yields
 $$
   \G(0)= 2\, (\I+\A)^{-1}(\G'\times \G'')(0),
   \quad {\hbox {i.e.,}}\quad
   \G(0)=2\, (\I+\A)^{-1} (c_0\bn(0)).
 $$
Hence,
 \begin{eqnarray*}
  |\Xn(0,t)-{\mathbf 0}|
  &=&
  |e^{{\A\over 2}\log t} {\sqrt t}\G(0)-{\mathbf 0}|=
  {\sqrt t} |\G(0)|
       \\
  &=&
  {\sqrt t} \, |2\, (\I+\A)^{-1}(c_0\bn(0))|
  \leq
  2c_0\, {\sqrt t}
  \| (\I+\A)^{-1} \|
  =
  2c_0{\sqrt t},
 \end{eqnarray*}
 with $c_0=c(0)$. The convergence result of $\Xn(s,t)$ stated in
 the proposition is an immediate consequence of the above inequalities.

 Finally, the continuity property of the maps
 $T^{\pm}:\ (\Gn(0),\Gn'(0),a)\in\R^3\times {\mathbb{S}}^2\times\R
 \longrightarrow \An^{\pm}
 $
 easily follows from the expressions of $\An^{\pm}$ in
 (\ref{15''}), the fact that $c(s)$ is bounded (see
 Lemma~\ref{for}), and the continuous dependence of the
 solutions of the ODE (\ref{13}) with respect to
 the data $(\Gn(0),\Gn'(0),a)$.
\end{proof}
%
\subsection{Asymptotics}\
\label{refined}
\smallskip

\noindent In the sequel, we will be devoted to study the
asymptotic behaviour of $\G(s)$ and $\Tn(s)$, quantifying their
wavelike behaviour through the vector $c(\nn-i\bn)$. To this end,
we will need a more exhaustive study of the properties of the
curvature and torsion functions related to the curve $\Gn(s)$, 
in the case when $c_{\pm\infty}\neq 0$. 
In order to clarify the exposition we have included this analysis in
a separated section (see Section~\ref{stct}).

We will continue to prove the following theorem:
\begin{theorem}\label{asym} Let $\Gn$ be a solution of
(\ref{12}). Then, the following asymptotics hold  as
$s\rightarrow\pm\infty$:
\begin{eqnarray*}
 & &
 \noindent {\it {i)}}
 \quad
 \G(s)=
 s\, e^{\A\log{|s|}}\An^{\pm}+
 {{e^{\A\log{|s|}}}\over s}
 \{ 2c^2_{\pm\infty}\Bn^{\pm}-\A \Bn^{\pm}\times\Bn^{\pm}\}
 - 4 \,{{c_{\pm \infty}\nn}\over{s^2}}
 +O\left({1\over {s^3}}\right),
    \\[2ex]
 & &
 \noindent {\it {ii)}}
 \quad
 \Tn(s)=
 e^{\A\log{|s|}}\Bn^{\pm}-2{{c_{\pm \infty}\bn}\over {s}}
 + O\left( {1\over s^2}\right).
      \\[2ex]
 & &
 \noindent {\it {iii)}}
 {\hbox { Moreover, if }} a\neq 0, B_3^{\pm}\neq \pm 1
 {\hbox { and }} c_{\pm\infty}\neq 0, {\hbox { then}}
     \\[2ex]
 & &
 \hspace{0.8cm}
 c_{\pm \infty}(\nn-i\bn)(s)=
 \,{{b_\pm e^{ia_\pm}}\over {|\A\Bn^{\pm}|^2}}\,
 e^{i\phi(s)}
 e^{\A\log{|s|}}
 \{ \A \Bn^{\pm}\times \Bn^{\pm}
 -i\A\Bn^{\pm}  \}+O\left({1\over|s|}\right),
     \\[2ex]
 & &
 {\hbox {where}}\quad \phi(s)=(s^2/ 4)-\gamma_{\pm}\log |s|.
     \\[2ex]
 & &
 \noindent {\it {iv)}}
 {\hbox { If }} a\neq 0  {\hbox { and  either }} B_3^{+} 
 {\hbox { or }} B_3^{-}\in\{+1,-1 \},
 {\hbox { then }} \Gn(s)=(0,0,\pm s).
 \end{eqnarray*}
 Here, the vectors $\Bn^{\pm}$, $\An^{\pm}$, and the constants $a_{\pm}$,
 $b_{\pm}$, $c_{\pm\infty}$ and $\gamma_{\pm}$ are defined as in
 the main theorem (see Theorem~\ref{mainG1}).
\end{theorem}
 \begin{remark}
  \label{reasym}
  The  maps $(\Gn(0),\Gn'(0),a)\longrightarrow (b_{\pm}, a_{\pm})$ are
  continuous. 
 (See the expression of
  $z_{+}=b_{+}e^{ia_{+}}$ in (\ref{cl26'}), (\ref{cl25'}), and
  notice  the continuous dependence of the solutions of the
  O.D.E. which satisfies $\twp$ with respect to the initial data. A
  similar expression can be found for $z_{-}=b_{-}e^{ia_{-}}$).
 \end{remark}
\begin{proof}[Proof of Theorem~\ref{asym}]
In what follows, we will reduce ourselves to consider the case
$s\rightarrow +\infty$. The case $s\rightarrow -\infty$ will
follow using the same arguments.

Recall that in the previous subsection we have obtained the
following expressions for the vectors $\G(s)$ and $\Tn(s)$ (see
(\ref{16}) and (\ref{17}))
\begin{equation}\label{16*}
 \G(s)= se^{\A\log{|s|}}\An^{+}+
 2s\,e^{\A\log{|s|}} \int_{s}^{+\infty} e^{-\A\log{|s'|}}
 {{c\bn}\over {(s')^2}}\, ds',\qquad s>0,
 \qquad {\hbox {and}}
\end{equation}
\begin{equation}\label{17*}
 \Tn(s)=e^{\A\log{|s|}}\Bn^{+}-
 2{{c\bn}\over {s}} +
 2 e^{\A\log{|s|}}\int_{s}^{+\infty}
 e^{-\A\log{|s'|}}{{(\I+\A)(c\bn)}\over {(s')^2}}\, ds',\ s>0.
\end{equation}
Observe that, from the Serret-Frenet system (see (\ref{star}))
$\nn_s=-c\Tn +(\tau-s/2)\bn+(s/2)\bn$, with $|\Tn|=|\nn|=|\bn|=1$.
Then,
\begin{eqnarray}\label{asym1}
 \int_{s}^{+\infty}e^{-\A\log{|s'|}}{{c\bn}\over {(s')^2}}\, ds'
 &=&
 2\int_{s}^{+\infty} e^{-\A\log{|s'|}}
 {{c}\over {(s')^3}}\, (s'/2)\bn\, ds'
        \nonumber \\
 &=&
 2\int_{s}^{+\infty} e^{-\A\log{|s'|}}
 {{c}\over {(s')^3}}(\nn_s+c\Tn -(\tau-s'/2)\bn)\, ds'.
\end{eqnarray}
Now, notice that
$$
 e^{-\A\log{|s|}} {{c\nn_s}\over {s^3}}
 =
 \left( e^{-\A\log{|s|}} {{c\nn}\over {s^3}} \right)_s-
 e^{-\A\log{|s|}}
 \left\{-c{{\A\nn}\over {s^4}} + \left( {c_s\over s^3}-
 {{3c}\over {s^4}}\right)\,\nn \right\}.
$$
Therefore, after an integration by parts and having into account
that $c(s)$ is bounded (see Lemma~\ref{for}), the identity 
(\ref{asym1}) rewrites as
\begin{eqnarray}\label{asym2}
 \int_{s}^{+\infty}e^{-\A\log{|s'|}}{{c\bn}\over {(s')^2}}\, ds'
 &=&
 -2\, e^{-\A\log{|s|}} {{c\nn}\over {s^3}}
 +2\int_{s}^{+\infty} e^{-\A\log{|s'|}}
 {{(\A+3\I)(c\nn)}\over {(s')^4}}\, ds'
     \nonumber \\
 & &
 + 2\int_{s}^{+\infty}
 {{e^{-\A\log{|s'|}}}\over {(s')^3}}
 ( c^2\Tn-c_s\nn -c(\tau-s'/2)\bn )
  \,ds'.
\end{eqnarray}
Recall that $\Gn(s)$ satisfies that
$$
 (\I+\A)\G -s \G'= 2\, (\G'\times \G''),
 \qquad {\hbox {i.e.,}} \qquad
 (\I+\A)\G-s\Tn= 2c\bn.
$$
Then, by using the Serret-Frenet formulae, the derivation of the
above equation with respect to $s$ concludes that the unit vector
$\Tn(s)$ is a solution of
\begin{equation}\label{asym3}
  \A\Tn= c(s-2\tau)\nn + 2 c_s\bn,
\end{equation}
from which it follows that
\begin{equation}\label{asym4}
 \A\Tn\times \Tn = 2\,(c_s\nn + c(\tau-s/2)\bn).
\end{equation}
The substitution of (\ref{asym4}) into (\ref{asym2}) yields
\begin{eqnarray}
 \label{asym4'}
 \int_{s}^{+\infty}e^{-\A\log{|s'|}}
 {{c\bn}\over {(s')^2}}\,ds'
 &=&
 - \,2\,e^{-\A\log{|s|}} {{c\nn}\over {s^3}}
 + 2\int_{s}^{+\infty} e^{-\A\log{|s'|}}
 {{(\A+3\I)\Tn_s}\over {(s')^4}}\, ds'
     \nonumber \\
 &-&
 \int_{s}^{+\infty} e^{-\A\log{|s'|}}
 {{\A\Tn\times \Tn}\over {(s')^3}}\, ds'
 +2 \int_{s}^{+\infty} e^{-\A\log{|s'|}}
 {{c^2\Tn}\over {(s')^3}}\, ds',
\end{eqnarray}
 (notice that, since $c\nn=\Tn_{s}$, the above identity holds true
 even if $c=0$, what can be checked directly).
 
 As a consequence, we conclude the following expression 
 for $\Gn(s)$ (\ref{16*})
\begin{eqnarray}\label{asym5}
 \G(s)
 &=&
 s\,e^{\A\log{|s|}}{\mathbf {A}}^{+}
 - 4\, {{c\nn}\over {s^2}}
    \nonumber \\
 &+&
 2s\, e^{\A\log{|s|}}
 \left\{
  2\,\int_{s}^{+\infty} e^{-\A\log{|s'|}}
  {{(\A+3\I)\Tn_{s}}\over {(s')^4}}\, ds'
  -\int_{s}^{+\infty}e^{-\A\log{|s'|}}
  {{\A\Tn\times\Tn}\over {(s')^3}}\,ds'
  \right.
     \nonumber  \\
  &+&
  \left.
  2\,\int_{s}^{+\infty} e^{-\A\log{|s'|}}
  {{c^2\Tn}\over {(s')^3}}\, ds'
 \right\}.
\end{eqnarray}
 Also, from Lemma~\ref{for}, it is not difficult to see that the
 integral in (\ref{asym4'}) is an error term of order 
 $O(1/s^2)$, as $s\rightarrow +\infty$. Then, from (\ref{17*}), we get that
\begin{equation}\label{17asym}
 \Tn(s)= e^{\A\log{|s|}}\Bn^{+}- 2 {{c\bn}\over {s}}+
 O\left( {{1}\over {s^2}}  \right),
 \qquad s\rightarrow +\infty.
\end{equation}
Now, we come back to the proof of the asymptotics of $\G(s)$ in
{\it {i)}}. To this end, we will analyze each of the integrals in
(\ref{asym5}).

Firstly, since $|\Tn|=1$, notice that
$$
 e^{-\A\log|s|} {{(\A+3\I)\Tn_{s}}\over {s^4}}
 =
 \left(
 e^{-\A\log|s|} {{(\A+3\I)\Tn}\over {s^4}}
 \right)_{s}
 +
 O\left( {1\over {s^5}} \right),
 \qquad s\rightarrow +\infty.
$$
Then,
\begin{equation}\label{asym6}
 \int_{s}^{+\infty}
 e^{-\A\log|s'|}
 {{(\A+3\I)\Tn_{s}}\over {(s')^4}}\, ds'
 = O\left({1\over  {s^4}} \right).
\end{equation}
Secondly, by taking into account (\ref{17}) and the asymptotic
development of $c^2(s)$ in Theorem~\ref{ct} part {\it {iv)}},
precisely
$$
 c^2(s)= c^2_{+\infty} +2\, {b_{+}\over s} \sin \tilde\phi(s)
 + O\left( {1\over s^2}\right),
 \qquad  {\hbox {with}}\qquad
 \tilde \phi(s)= a_{+} + (s^2/4)-\gamma_{+}\log(s),
$$
we obtain
$$
 c^2\Tn(s)=
 e^{\A\log|s|} {\mathbf B}^{+}
 \left(
  c^2_{+\infty} + 2\, {b_{+}\over s}\sin\tilde\phi(s)
 \right)
 - 2 c^2_{+\infty} {{c\bn}\over {s}}
 +O\left( {1\over {s^2}}\right),
$$
so that
\begin{eqnarray}\label{asym7}
 \int_{s}^{+\infty}
 e^{-\A\log|s'|} {{c^2\Tn}\over {(s')^3}}\, ds'
 &=&
 c^2_{+\infty}\, {{{\mathbf B}^{+}}\over {2s^2}}+
 2 b_{+} {\mathbf B}^{+}
 \int_{s}^{+\infty}
  {{\sin \tilde\phi(s)}\over {(s')^4}}\, ds'
  \nonumber \\
 &-&
 2 \, c^2_{+\infty}\int_{s}^{+\infty}
 e^{-\A\log|s'|} {{c\bn}\over {(s')^4}}\, ds'
 + O\left(  {1\over s^4} \right).
\end{eqnarray}
 On the one hand, notice that
 $\tilde\phi'(s)= s/2-\gamma_{+}/s\neq 0$ for $s$ sufficiently
 large. Therefore, an integration by parts argument gives that
 \begin{eqnarray}\label{asym8}
 \int_{s}^{+\infty}
 {{\sin \tilde\phi}\over {(s')^4}}\, ds'
 =
 O\left( {1\over {s^4}} \right).
 \end{eqnarray}
 On the other hand, a similar argument to that given in obtaining
 (\ref{asym4'}), and using the fact that $c(s)$ is bounded (see Lemma
 \ref{for})  conclude that
 \begin{equation}\label{asym9}
 \int_{s}^{+\infty} e^{-\A\log|s'|}
 {{c\bn}\over {(s')^4}}\, ds'
 = O\left( {1\over s^4}  \right),
 \qquad{\hbox {as}}\qquad
 s\rightarrow +\infty,
 \end{equation}
 by using once again that $c(s)$ is bounded (see Lemma \ref{for}).
 Substituting (\ref{asym8}) and (\ref{asym9}) into (\ref{asym7}), one gets
 \begin{equation}\label{asym10}
  \int_{s}^{+\infty}
  e^{-\A\log|s'|} {{c^2\Tn}\over {(s')^3}}\, ds'=
  c^2_{+\infty}\, {{{\mathbf B}^{+}}\over {2s^2}}
  + O\left( {1\over {s^4}} \right).
 \end{equation}
 Finally, from (\ref{17asym})
 \begin{eqnarray*}
 \A\Tn\times \Tn
 &=&
 e^{\A\log|s|} (\A{\mathbf B}^{+}\times {\mathbf B}^{+})
 - {2\over s} \, \left( e^{\A\log|s|} \A{\mathbf B}^{+} \right)\times (c\bn)
     \\
 &-& {2\over s}\, \A(c\bn)\times \left(  e^{\A\log|s|} {\mathbf B}^{+}\right)
 + O\left( {1\over s^2} \right).
 \end{eqnarray*}
 Therefore, from (\ref{asym9})
 \begin{equation}\label{asym11}
  \int_{s}^{+\infty}  e^{-\A\log|s'|}
  {{\A\Tn\times \Tn}\over {(s')^3}}\, ds'
  =
  {1\over {2s^2}}\, \A {\mathbf B}^{+}\times {\mathbf B}^{+}
  + O\left( {1\over s^4} \right).
 \end{equation}
 The substitution of (\ref{asym6}), (\ref{asym10}) and (\ref{asym11}) into
 (\ref{asym5}) concludes the asymptotic behaviour of $\G(s)$ stated in {\it
 i)}, that is
 \begin{equation}\label{16asym}
 \G(s)=
 s\, e^{\A\log{|s|}}\An^{+}+
 {{e^{\A\log{|s|}}}\over s}
 \{ 2c^2_{+\infty}\Bn^{+}-\A \Bn^{+}\times\Bn^{+}\}
 - 4 \,{{c\nn}\over{s^2}}
 +O\left({1\over {s^3}}\right).
 \end{equation}
From (\ref{17asym}) and (\ref{16asym}), we see that giving more
accurate asymptotics of $\G$ and $\Tn$ implies the study of the
associated vectors  $c\bn$ and  $c\nn$.
We will continue to obtain a closed formula for $c(\nn-i\bn)$
which is valid for $s$ sufficiently large, $a\neq 0$, $B_3^{+}\neq
\pm 1$  and $c_{+\infty}\neq 0$. 

We firstly observe that, because $\A$ is an antisymmetric matrix,
$\A\Tn\cdot \Tn=0$. Then, under the above assumptions on $a$ 
and $B_{3}^{+}$, the expression for $\Tn(s)$ in (\ref{17asym})
asserts us that $|T_3(s)|<1$  when $s$ is sufficiently large. Then
$\A\Tn\neq {\mathbf 0}$ and we can consider
 $\{\A\Tn, \Tn, \A\Tn\times \Tn\}$ as a basis of
$\R^3$. Thus, the vector $c(\nn-i\bn)$ can be written as a linear
combination of the elements of this basis.

To this end, recall that (see (\ref{asym3}) and (\ref{asym4}))
$$
  \A\Tn= -2 c(\tau-s/2)\nn + 2 c_s\bn,
  \qquad {\hbox {and}} \qquad
  \A\Tn\times \Tn= 2c_s \nn + 2 c (\tau-s/2)\bn.
$$
Therefore, because $\nn\bot\bn$, $\bn\bot \Tn$, $\nn\bot \Tn$,
$|\nn|=|\bn|=|\Tn|=1$, we obtain 
\begin{eqnarray*}
 & &
 c(\nn-i\bn)\cdot \A\Tn= -2 c^2(\tau-s/2)
 - i {{dc^2}\over {ds}}=
 -i\left(
   {{dc^2}\over {ds}}- 2i c^2(\tau-s/2),
 \right)
    \\
 & &
 c(\nn-i\bn)\cdot \Tn=0,
 \qquad {\hbox {and}}\qquad
    \\
 & &
 c(\nn-i\bn)\cdot (\A\Tn\times \Tn)=
 {{dc^2}\over {ds}}- 2i c^2(\tau-s/2).
\end{eqnarray*}
As a consequence, we conclude that
\begin{equation}
 \label{asym12}
 c(\nn-i\bn)=
 \left(
   {{dc^2}\over {ds}}- 2i c^2(\tau-s/2)
 \right)\,
 \left(
  {{\A\Tn\times \Tn}\over {|\A\Tn\times \Tn|^2}}
   - i {{\A\Tn}\over {|\A\Tn|^2}}
 \right)
\end{equation}
(notice that (\ref{asym12})) is valid whenever $T_3(s)\neq \pm 1$
and $a\neq 0$. In particular, it is valid if $c=0$. In that case, both
sides are understood to be zero).

Now, by using the asymptotic behaviour of $\Tn(s)$ given in
(\ref{17asym}), we get that
\begin{equation}
 \label{asym13}
 \A\Tn\times \Tn - i \A \Tn =
 e^{\A\log|s|}
 \{ \A {\mathbf B}^{+}\times {\mathbf B}^{+}
 - i \A {\mathbf B}^{+}\} +
 O\left( {1\over {|s|}}  \right).
\end{equation}
Also, it is satisfied that
\begin{equation}\label{asym13'}
 |\A\Tn\times \Tn |^2(s)=|\A\Tn|^2(s)=
 |\A\Bn^{+}|^2+
 O\left( {{1}\over {|s|}} \right).
\end{equation}
On the other hand, from the asymptotics related to $c(s)$ and
$\tau (s)$ in Theorem~\ref{ct}, it follows that
\begin{equation}
 \label{asym14}
 {{dc^2 }\over {ds}}-2ic^2(\tau-s/2)=
 b_{+}\, e^{ia_+} e^{i\phi(s)}
 +O\left({1\over {|s|}}  \right),
 \qquad s\rightarrow +\infty,
\end{equation}
with $\phi(s)=(s^2/4)-\gamma_+\log|s|$.

 From (\ref{asym12})-(\ref{asym14}) we obtain
$$
 c(\nn-i\bn)= {{b_\pm e^{ia_+}}\over {|\A\Bn^{+}|^2}}\,e^{i\phi(s)}
 e^{\A\log{|s|}}
 \{  \A \Bn^{+}\times \Bn^{+}
 -i\A \Bn^{+}  \}+
 O\left({1\over {|s|}} \right),
$$
 for some constant $a_+\in[0,2\pi)$.

 The proof of {\it {(i)-(iii)}} is now an immediate 
 consequence of the above identity, (\ref{17asym}), (\ref{16asym}),
 and (\ref{cl5}).  
 Notice that (\ref{cl5}) follows from (\ref{cl3}),
 which is also true if $c_{\pm\infty}=0$ because it just involves
 $|f|^2=c^2$ and $h=-(\A \Tn\cdot \Tn')/2$. In fact, (\ref{cl3}) easily
 follows from Lemma \ref{for}, and using that, from the equation
 (\ref{12}), $h$ rewrites as $h=-a(G_3-sT_3)/4$ (see Remark~\ref{reyh}).
 Besides, observe that the constants in
 Theorem~\ref{asym} (or, equivalently in Theorem~\ref{mainG1} in the
 introduction) are directly deduced from the ones in
 Theorem~\ref{ct} and Lemma~\ref{for}. 
 Finally, we prove the part {\it {(iv)}}.
 Let $a\neq 0$ and $B^{+}_{3}=+1$. Then,
 $$
  a^2(1-T_3^2){|}_{s=\infty}=|\A\Tn\times \Tn|^2{|}_{s=\infty}=0,
 $$
 and therefore $h(+\infty)=y(+\infty)=0$ (recall $h=-(\A \Tn\cdot \Tn')/2$
 and $y=(\A \Tn\times \Tn)\cdot \Tn'$).

 Assume now that $c_{+\infty}=0$, so that $\alpha^2=a^2$, by using 
 Lemma~\ref{for}. Then, $y=h=0$ solves (\ref{cl9}) (recall that
 $E(0)=a^2/4$ -see (\ref{ct1})) so that $(dc^2/ds)=0$,
 and since $c_{+\infty}=0$, we obtain $c(s)=0$. Hence, $\Gn(0,0,s)$.

 Secondly, we assume that $c_{+\infty}\neq 0$. Then $c(s)\neq 0$ for
 $s$ large enough and (\ref{asym12}) is valid whenever $\A \Tn\neq 0$
 (i.e., $a\neq 0$ and $T_{3}\neq \pm 1$).

 Let us prove that $T_3(s)\neq \pm 1$ for $s$ large enough. Notice
 that, since $B_{3}^{+}=+1$, it is enough to see that $T_3(s)\neq
 1$. To this end, assume on the contrary that there exists $s_0$
 large enough such that $T_3(s_0)=1$. Then, there exists
 $\{s_n\}_{n}\rightarrow s_0$, as $n\rightarrow +\infty$, such that
 $T_3(s_n)\neq \pm 1$ and $T_3(s_n)\rightarrow 1$. Hence, from
 (\ref{asym12}), we get that $c(s_0)=0$, which contradicts the
 assumption $c_{+\infty}\neq 0$. Therefore, (\ref{asym12}) holds true
 for $s$ large enough and letting $s\rightarrow +\infty$ we get that
 $c_{+\infty}=0$, which also leads to contradiction.

 The same arguments as above prove that, in the cases  when $B_3^{-}=1$ or
 $B_3^{\pm}=-1$, $\Gn(s)=(0,0,+s)$ or $\Gn(s)=(0,0,-s)$,
 respectively. This concludes the proof.

\end{proof}
%
 \subsection{Scattering problem for the curve $\Gn$}\
 \label{scatteringG}
\smallskip

\noindent We will finish this section proving the
Theorem~\ref{mainG2} concerning with the scattering problem for
$\Gn(s)$. More precisely,
\begin{theorem}\label{st}
Given $a\neq 0$, $\Bn^{+}=(B_1^+, B_2^+, B_3^+)\in {\mathbb{S}}^2$
with $B_3^+\neq \pm 1$, $a_{+}\in [0,2\pi)$ and
$b_{+}\geq 0$, there exists a  solution $\Gn(s)$  of
\begin{equation}
 \label{st0}
 (\I +\A)\Gn -s\Gn '= 2\Gn'\times \Gn '',
 \qquad |\Gn'|=1
\end{equation}
satisfying the following identities:
\begin{eqnarray}
 & &
 \label{st1a}
 \noindent {\it {i)}}
 \quad
 \lim_{s\rightarrow +\infty}
 e^{-\A\log{|s|}} {{\Gn(s)}\over s}=\An^+,
       \\[2ex]
 & &
 \label{st1}
  \noindent {\it {ii)}}
 \quad
 \lim_{s\rightarrow +\infty} e^{-\A\log{|s|}} \Tn(s)=\Bn^+,
 \qquad \quad {\hbox {and}}
        \\[2ex]
 & &
 \label{st2}
 \noindent {\it {iii)}}
 \quad
 \lim_{s\rightarrow +\infty}e^{-i{\phi}(s)} e^{-\A \log{|s|}} (\Tn'-i\Tn\times
 \Tn')(s)=
 {{b_{+} e^{ia_+}}\over {|\A\Bn^{+}|^2}}\,(\A\Bn^{+}\times \Bn^{+}-i\A\Bn^{+}),
\end{eqnarray}
where  $\An^{+}=(\I+\A)^{-1}\Bn^{+}$,
 ${\phi}(s)= (s^2/4)-\gamma_{+}\log |s|$,
$$
 \gamma_{+}=3a B_3^{+}+\alpha
 \qquad {\hbox {and}}\qquad
 \alpha=-aB_3^+ -{{b_{+}^{2}}\over {a^2 [1-(B_3^+)^2]}}\ .
$$
Moreover, if $b_{+}>0$, the solution is unique. A similar result can 
be obtained at $s\rightarrow -\infty$.
\end{theorem}
%
\begin{proof}[Proof of Theorem~\ref{st}]

\noindent {\bf {Existence.}}
 Firstly, observe that, whenever $\Tn$ is a solution of
 \begin{equation}
 \label{s2}
   \A\Tn = s\Tn'+2\Tn\times \Tn ''
 \end{equation}
 with initial data
 \begin{equation}
  \label{s3}
 |\Tn(0)|=1 \qquad {\hbox {and}}\qquad
 \Tn(0)\cdot \Tn'(0)=0,
 \end{equation}
then $\Gn(s)$ such that
\begin{equation}
 \label{s4}
 \Gn'(s)=\Tn(s)
 \qquad {\hbox {and}}\qquad
 (\I+\A)\Gn(0)=2\, \Tn(0)\times \Tn'(0)
\end{equation}
satisfies (\ref{st0}), and $(\I +\A)\Gn(0)\cdot \Gn'(0)=0$.

Indeed, from (\ref{s2}),
\begin{eqnarray*}
 {{d\  }\over {ds}} |\Gn'|^2
 &=&
 {{d|\Tn|^2 }\over {ds}}= 2\, \Tn\cdot \Tn'=0,\quad s\neq 0
 \qquad {\hbox {and}}\qquad
    \\[1ex]
 {{d\ }\over {ds}} (\I+\A)\Gn
 &=&
 \A\Tn+\Tn =s\Tn' +2 \Tn\times \Tn''+\Tn =
 {{d\ }\over {ds}} (s\Tn +2\, \Tn\times \Tn').
\end{eqnarray*}
Then, (\ref{st0}) follows from the above identities and the
initial conditions in (\ref{s3}) and (\ref{s4}), respectively.

Moreover, from the results in Section~\ref{main}, the condition
(\ref{st1a}) is satisfied with
$$
  \An^{+}=(\I+\A)^{-1}\left(\lim_{s\rightarrow +\infty}
  e^{-\A\log|s|}\Tn(s)\right).
$$
Previous remarks allow us to reduce ourselves to prove the existence
of a solution $\Tn$ of (\ref{s2}) and (\ref{s3}) satisfying the
limiting conditions (\ref{st1}) and (\ref{st2}).

It is easy to see the existence of a global solution of (\ref{s2})
and (\ref{s3}): Given a fixed initial data $(\Tn(0), \Tn'(0))$ as
in (\ref{s3}), consider $\Gn(s)$ such that
\begin{equation}
 \label{s5}
 \Gn''={1\over 2}\, (\I+\A)\Gn\times \Gn',
\end{equation}
\begin{equation}
 \label{s6}
 (\I+\A)\Gn(0)=2\, \Tn(0)\times \Tn'(0)
 \qquad {\hbox {and}}\qquad
 \Gn'(0)=\Tn(0).
\end{equation}
Notice that (\ref{s3}) and (\ref{s6}) imply that $ |\Gn'(0)|=1$
and $(\I+\A)\Gn(0)\cdot \Gn'(0)=0$ and thus, there exists $\Gn(s)$
global solution of (\ref{s5})-(\ref{s6}) (see Section~\ref{main}).
Moreover, $\Gn\in \C^{\infty}(\R;\R^3)$ and $|\Gn'(s)|=1,\
\forall\, s\in\R$.

Now take $\Tn(s)=\Gn'(s)$. Then, deriving (\ref{st0}) we obtain
that $\Tn(s)$ is a solution of (\ref{s2}). On the other hand, from
(\ref{s5}) at $s=0$, $2\,\Gn''(0)=(\I+\A)\Gn(0)\times \Gn'(0)$.
Therefore,
$$
 \Tn(0)=\Gn'(0),\qquad |\Gn'(0)|=1
 \qquad {\hbox {and}}\qquad
$$
$$
 \Tn(0)\cdot \Tn'(0)= \Gn'(0)\cdot \Gn''(0)=
 \Gn'(0)\cdot {1\over 2}[(\I+\A)\Gn(0)\times \Gn'(0)]=0,
$$
so that the initial conditions in (\ref{s3}) are satisfied.

Finally, recall in the sequel that, if $\Tn$ is a solution of
(\ref{s2}) and (\ref{s3}), then $\Gn$ defined in (\ref{s4})
verifies (\ref{st0}). Hence, we have already proved (see
(\ref{14a}), (\ref{14c}) and Lema~\ref{for}) that $\Tn$ satisfies
the following properties:
\begin{equation}
 \label{s8}
 |\Tn(s)|=1,\qquad
 \Tn(s)\cdot\Tn'(s)=0,
 \qquad {\hbox {and}}\qquad
 |\Tn'(s)|^2
 =-aT_3(s)-\beta,\quad  \forall\, s\in\R,
\end{equation}
with $\beta= -|\Tn'(0)|^2-aT_3(0)$.

We now come back to the proof of Theorem~\ref{st}. For a fixed
$a\neq 0$, $\Bn^{+}\in {\mathbb{S}}^2$ with $B^{+}_3\neq \pm 1$,
$a_+\in[0,2\pi)$ and $b_+\geq 0$, there exists a unique $\alpha\in
\R$ such that
\begin{equation}
 \label{s11}
 b_{+}=\sqrt{a^2 (-aB_{3}^{+}-\alpha) (1-(B_{3}^{+})^2)},
 \qquad {\hbox {i.e.,}} \qquad
 \alpha= -aB_{3}^{+} -
 {{b_{+}^{2}}\over {a^2(1-(B_{3}^{+})^2)}}.
\end{equation}
In particular, as we have already observed at the beginning of the
proof, the result in Theorem~\ref{st} will follow from the
existence, for this value of $\alpha$, of initial data
$(\Tn(0),\Tn'(0))$ such that
\begin{equation}
 \label{s12}
 |\Tn(0)|=1,\qquad \Tn(0)\cdot \Tn'(0)=0
 \qquad {\hbox {and}}\qquad
 |\Tn'(0)|^2=-a T_3(0)-\alpha,
\end{equation}
and $\Tn(s)$ solution of
\begin{equation}
 \label{s13}
  \A \Tn= s\Tn' + 2\, \Tn\times \Tn'',
\end{equation}
associated to this initial data, satisfying the limiting
conditions (\ref{st1}) and (\ref{st2}). The existence of such this
solution is based on a compactness argument.

To this end, consider the compact set
$$
 K^{\alpha}=\{ (\un,\vn)\in\R^3\times\R^3\ / \
 |\un|=1,\
 \un\cdot \vn=0,\ |\vn|^2=-au_3-\alpha \},
$$
with $\un=(u_1, u_2, u_3)$, and for $s\geq 1$ define the operator
$W(s)$ given by
$$
 W(s)
 \left(
 \begin{array}{c}
  \Tn(0)\\
  \Tn'(0)
 \end{array}
 \right)=
 \left(
 \begin{array}{c}
  e^{-\A\log|s|}\Tn(s)\\
  e^{-i{\tilde\phi}(s)}e^{-\A\log|s|}(\Tn'-i \Tn\times \Tn')(s)
 \end{array}
 \right),
 \qquad
 \forall\, (\Tn(0),\Tn'(0))\in K^{\alpha},
$$
where $\Tn$ is a solution of (\ref{s12})-(\ref{s13}) and
$\tilde\phi(s)=(s^2/4)-(3T_3(s)+\alpha)\log|s|$.

Recall that we have already proved the global existence of
a solution of the problem (\ref{s12})-(\ref{s13}), so that $W(s)$ is
well-defined. Also, it is easy to see that $W(s)$ is a continuous
operator, $\forall s\geq 1$.

Assume momentarily that the following claim holds:

\noindent {\it {Claim. }} Let be  $\Bn^{+}$, $a_+$, $b_+$  and
$\alpha$ as above, and $s_0\geq 1$. Then, there exists initial
data $(\Tn(0), \Tn'(0))\in K^{\alpha}$, depending on $s_0$, and
$\Tn(s)$ solution of (\ref{s13}) associated to this initial data
such that
\begin{equation}
 \label{s14}
 e^{-\A\log|s_0|}\Tn(s_0)=\Bn^{+},
 \qquad \quad {\hbox {and}}
\end{equation}
\begin{equation}
 \label{s15}
 e^{-i{\tilde\phi}(s_0)} e^{-\A\log|s_0|} (\Tn'-i\Tn\times \Tn')(s_0)=
 {{b_+e^{ia_+}}\over {|\A\Bn^{+}|^2}}\,(\A\Bn^{+}\times \Bn^{+}-i \A\Bn^{+}),
\end{equation}
with $T_3(s_0)=B_3^{+}$.

Now, choose a sequence $\{s_n\ :\ s_n\geq 1  \}_{n\in\N}$ such
that $s_n\rightarrow +\infty$, as $n\rightarrow +\infty$. Then,
for any fixed $n\in\N$ the above claim ensures the existence of
initial data $(\Tn_{n}(0), \Tn'_{n}(0)) \in K^{\alpha}$, and
$\Tn_{n}(s)$ solution of (\ref{s13}) with this initial data such
that the following identity holds
\begin{eqnarray}
 \label{s19}
 W(s_n)
 \left(
 \begin{array}{c}
  \Tn_{n}(0) \\
  \Tn'_{n}(0)
 \end{array}
 \right)
 &=&
  \left(
 \begin{array}{c}
 e^{-\A\log|s_n|}\Tn_{n}(s_n) \\
 e^{-i{\tilde\phi}(s_n)}e^{-\A\log|s_n|}(\Tn'_{n}-i \Tn_{n}\times\Tn'_{n})(s_n)
 \end{array}
 \right)
        \nonumber \\
 &=&
 \left(
 \begin{array}{c}
  \Bn^+\\
  {{b_+e^{ia_+}}\over {|\A\Bn^{+}|^2}}\, (\A\Bn^+\times \Bn^{+}-i\A\Bn^{+})
 \end{array}
 \right),
 \qquad \forall n\in\N.
\end{eqnarray}
The compactness of the set $K^{\alpha}$ yields the existence of
$(\Tn(0),\Tn'(0))\in K^{\alpha}$ such that
\begin{equation}
 \label{s20}
  \Tn(0)=\lim_{k\rightarrow +\infty}\Tn_{n_k}(0)
  \qquad {\hbox {and}}\qquad
  \Tn'(0)=\lim_{k\rightarrow +\infty}\Tn'_{n_k}(0),
\end{equation}
for some subsequence of $\{(\Tn_{n}(0),\Tn'_{n}(0))\}_{n\in\N}$.

Next, consider the operator
 $W(\infty)=\lim_{s\rightarrow +\infty}W(s)$
defined as follows
\begin{eqnarray}
 \label{s21}
 W(\infty)
 \left(
 \begin{array}{c}
  \Tn(0) \\
  \Tn'(0)
 \end{array}
 \right)
 &=&
 \lim_{s\rightarrow +\infty}
 W(s)
 \left(
 \begin{array}{c}
  \Tn(0) \\
  \Tn'(0)
 \end{array}
 \right)
   \nonumber \\
 &=&
 \lim_{s\rightarrow +\infty}
 \left(
 \begin{array}{c}
  e^{-\A\log|s|} \Tn(s) \\
 e^{-i{\tilde\phi}(s)} \, e^{-\A\log|s|}
 (\Tn' -i \Tn\times \Tn')(s)
 \end{array}
 \right),
\end{eqnarray}
where $\Tn(s)$ satisfies (\ref{s13}) with initial data
$(\Tn(0),\Tn'(0))\in K^{\alpha}$.

Notice that, from (\ref{s8}) and the fact that
$(\Tn(0),\Tn'(0))\in K^{\alpha}$, it follows that
$$
 |\Tn (s)|=1
 \qquad {\hbox {and}}\qquad
 |\Tn'(s)|^2= -aT_3(s)-\alpha
 \qquad {\hbox {with}}\qquad
 \alpha= -aT_3(0)- |\Tn'(0)|^2.
$$
Also, from the Serret-Frenet formulae (see~(\ref{star})),
$c(\nn-i\bn)= \Tn'-i\Tn\times \Tn'$, where we understand that both
sides are zero if $\Tn'=0$.

The above remarks and the asymptotic behaviour 
of $\Tn$ and $c(\nn-i\bn)$ (see parts {\it {i)}}-{\it {iii)}} 
of Theorem~\ref{mainG1} and Remark~\ref{remainG1})
assert the existence of the limits in the definition of
$W(\infty)$, so that $W(\infty)$ is well-defined. Moreover, from
Remark~\ref{reasym} it follows that $W(\infty)$ is continuous.

Now take the initial data $(\Tn(0),\Tn'(0))$ in (\ref{s20}). Then,
the continuity property of $W(\infty)$ and the identity
(\ref{s19}) yield
\begin{eqnarray}
 \label{s23}
 W(\infty)
  \left(
 \begin{array}{c}
  \Tn(0) \\
  \Tn'(0)
 \end{array}
 \right)
 &=&
 \lim_{k\rightarrow \infty}
 W(\infty)
 \left(
 \begin{array}{c}
  \Tn_{n_k}(0) \\
  \Tn'_{n_k}(0)
 \end{array}
 \right)
    \nonumber   \\
 &=&
 \lim_{k\rightarrow \infty}
 \lim_{s\rightarrow +\infty}
 \left(
 \begin{array}{c}
  e^{-\A\log|s_{n_k}|} \Tn_{n_k}(s_{n_k}) \\
  e^{-i{\tilde\phi}(s_{n_k})} \,e^{-\A\log|s_{n_k}|}
 (\Tn'_{n_k} -i \Tn_{n_k}\times \Tn'_{n_k})(s_{n_k})
 \end{array}
 \right)
    \nonumber \\
 &=&
 \left(
 \begin{array}{c}
  \Bn^{+} \\
  {{b_{+}e^{ia_+}}\over {|\A\Bn^{+}|^2}}\, 
  (\A\Bn^{+}\times \Bn_{+}- i \A \Bn^{+})
 \end{array}
 \right),
\end{eqnarray}
because $s_{n_k}\rightarrow +\infty$, as $k\rightarrow +\infty$.

From (\ref{s23}) and the definiton of the action of $W(\infty)$ on
the initial data in (\ref{s20}), we conclude that the solution of
(\ref{s12}) associated to the initial data in (\ref{s20})
satisfies the limiting conditions in (\ref{st1}) and (\ref{st2}),
that is
$$
 \lim_{s\rightarrow +\infty}
 \left(
 \begin{array}{c}
  e^{-\A\log|s|} \Tn(s) \\
 e^{-i{\tilde\phi}(s)} \,e^{-\A\log|s|}
 (\Tn' -i \Tn\times \Tn')(s)
 \end{array}
 \right)
 =
 \left(
 \begin{array}{c}
  \Bn^{+} \\
  {{b_{+}e^{ia_+}}\over {|\A\Bn^{+}|^2}}(\A\Bn^{+}\times \Bn^{+}- i \A \Bn^{+})
 \end{array}
 \right).
$$
The result in Proposition~\ref{st} is an immediate consequence of
this identiy, by noticing that the first limiting condition gives
that $\lim_{s\rightarrow +\infty}T_3(s)=B_3^{+}$ .

We now come back to the proof of the claim:
Let $s_0\geq 1$. Firstly, we notice that the global existence of
solution of (\ref{s12})-(\ref{s13}) and the identities in
(\ref{s8}) reduce the problem to finding  $(\Tn(s_0),\Tn'(s_0))\in
K^{\alpha}$ satisfying the conditions in the statement of the
claim, that is (\ref{s14}) and (\ref{s15}).

Then, from (\ref{s14}), it follows that $\Tn(s_0)$ is given by
\begin{equation}
 \label{s16}
 \Tn(s_0)= e^{\A\log|s_0|}\Bn^{+},
\end{equation}
and in particular $T_3(s_0)=B_3^{+}$.

On the other hand, by taking into account (\ref{s14}) in
(\ref{s15}), we obtain that $\Tn'(s_0)$ should satisfy that
$$
 e^{-i{\tilde\phi}(s_0)} e^{-\A\log|s_0|}(\Tn'-i \Tn \times \Tn')(s_0)=
 {{b_+ e^{ia_+}}\over {|\A\Bn^{+}|}}\, e^{-\A \log|s_0|} 
 (\A\Tn\times \Tn- i\A\Tn)(s_0),
$$
or, equivalently,
\begin{equation}
 \label{s18}
 (\Tn'-i \Tn\times \Tn')(s_0)= \omega \,(\A\Tn\times
 \Tn-i\A\Tn)(s_0),
\end{equation}
for some known $\omega\in \C$ which depends on $B_{3}^{+}$, $a_+$,
$b_+$ and $\alpha$, and such that $|\omega|^2=
b_{+}^{2}/|\A\Bn^{+}|^4$.

In order to prove the existence of $\Tn'(s_0)$ solution of the
above equation, firstly observe that, because $B_{3}^{+}\neq \pm
1$ and $\Tn(s_0)=e^{\A\log|s_0|}\Bn^{+}$ (see~(\ref{s16})), it
follows that $T_3(s_0)\neq \pm 1$, $s_0\geq 1$, and therefore
$\A\Tn(s_0)\neq {\mathbf {0}}$. As a consequence, the set of
vectors
$$
 \B=\{ \Tn(s_0), (\A\Tn\times\Tn)(s_0), \A\Tn(s_0)
 \}
$$
is a basis in $\R^3$ (recall that, because $\A$ is an
antisymmetric matrix, $(\A\Tn\cdot \Tn)(s_0)=0 $).

In the coordinate system associated to $\B$, $\Tn=(1,0,0)$ and, by
taking into account that $(\Tn\cdot\Tn')(s)=0$ (see~(\ref{s8})),
$\Tn'(s_0)=\beta_1\,(\A\Tn\times\Tn)(s_0)+\beta_2 \A\Tn(s_0)$,
where the (real) scalars $\beta_1$ and $\beta_2$ are uniquely
determined by the equation $\beta_1 + i\beta_2=\omega$, just by
rewriting the equation (\ref{s18}) in the coordinate system
associated to $\B$.

Now it is easy to check that $(\Tn(s_0), \Tn'(s_0))\in
K^{\alpha}$. Indeed, on the one hand, from (\ref{s16}) and
(\ref{s18}),
$$
  |\Tn(s_0)|=1
  \qquad {\hbox {and}}\quad
  \Tn(s_0)\cdot \Tn'(s_0)=0.
$$
On the other hand, from (\ref{s18}), $T_3(s_0)=B_3^{+}$ and
$|\omega|^2= b^{2}_{+}/ |\A\Bn|^4$, we get that
$$
  2|\Tn'(s_0)|^2= 2|\omega|^2 |\A\Tn(s_0)|^2= 2\,
  {{b_{+}^2}\over {|\A\Bn^{+}|^4}}\, |\A\Bn^{+}|^2.
$$
Then, using the definition of $b_{+}$ in (\ref{s11}) and the fact
that $|\A\Bn^{+}|^2= a^2(1-(B_3^{+})^2)$, from previous identity
we conclude that
$$
  |\Tn'(s_0)|^2= -aB_3^{+}-\alpha= -aT_3(s_0)-\alpha
$$
 This finishes the proof of the claim.

 The proof of the uniqueness is based on a uniqueness result for
 the solutions of the self-similar Schr\"odinger equation that
 will be proved in Section~\ref{selfNLS} (see Lemma~{\ref{unicidad}}
 in Section~\ref{selfNLS}).

\noindent {\bf {Uniqueness.}}
Assume $\Gn^{j}$, $j=1,2$, are two solutions of (\ref{st0})
satisfying (\ref{st1a}), (\ref{st1}) and (\ref{st2}) for some given
$a\neq 0$, ${\mathbf {B}}^{+}\in{\mathbb {S}}^{2}$ with $B^{+}_{3}\neq
\pm 1$, $a_{+}\in [0, 2\pi)$ and $b_{+}>0$. Let $c_j$ and $\tau_j$
denote respectively the curvature and torsion related to

Consider also the functions $f_j(s)$, $j=1,2$ defined through the
relations (\ref{24}), i.e,
$$
 f_j(s)e^{i{s^2\over 4}}= c_j(s) e^{\int_{0}^{s}\tau_j(s')\, ds'},
 \qquad j=1,2,
$$
and the associated functions $h_j=\Im m (\bar f_j f'_j)$ and $y_j=
d|f_j|^2/ds$.

We begin to prove that $f_j,$ $j=1,2$ satisfy the hypothesis in
Lemma~\ref{unicidad}.

To this end, recall that it has been already shown that
\begin{equation}
 \label{p0}
 \A\Tn\times \Tn= 2\, c(\tau-s/2)\bn +2\, c_s \nn
 \qquad {\hbox {and}} \qquad
 c^2(s)= -aT_3(s)+\alpha,
\end{equation}
for any $\Gn(s)$ solution of equation (\ref{st0}) (from
(\ref{asym4}) and by using the Serret-Frenet formulae in
(\ref{s8})).

Firstly, notice that from (\ref{st2}), it is easy to check that
$$
 2\, {{|b_{+}|^2}\over {|\A \Bn^{+}|^2}}
 =
 \lim_{s\rightarrow +\infty} |\Tn'_j-i \Tn_j \times \Tn'_j|^2(s)
 =
 \lim_{s\rightarrow +\infty} 2\,|\Tn'_j|^2(s)= 
 2\,c^2_{+\infty,j},
 \qquad j=1,2.
$$
Then,
\begin{equation}
 \label{p2}
 c^2_{+\infty,j}= c^2_{+\infty}
 \qquad {\hbox {with}}\qquad
 c^2_{+\infty}=
 {{|b_{+}|^2}\over {|\A \Bn^{+}|^2}}\neq 0.
\end{equation}
because $B_3^{+}\neq \pm 1$, $a\neq 0$  and $b_{+}>0$ or, equivalently,
\begin{equation}
 \label{p3}
 |f_j|^2_{+\infty}= |f|^2_{+\infty}
 \qquad {\hbox {with}}\qquad
 |f|^2_{+\infty}= 
 {{|b_{+}|^2}\over {|\A \Bn^{+}|^2}}\neq 0.
\end{equation}
Also, from (\ref{st1}),
 \begin{equation}
  \label{p5}
  \lim_{s\rightarrow +\infty} |\A\Tn_{j}\times \Tn_{j}|(s)=
  |\A\Bn^{+}\times \Bn^{+}|,
  \qquad j=1,2.
 \end{equation}
Besides, from (\ref{p0}),
$$
 |\A\Tn_j\times \Tn_j|^2(s)=
 4\, ( (c'_j)^2 + c_j^2(\tau_j- s/2)^2),
 \qquad \forall\, s\in\R ,\qquad j=1,2,
$$
and, by taking into account the above identities together with
(\ref{25}), we obtain 
\begin{equation}
 \label{p6}
 |f'_j|_{+\infty}= |f'|_{+\infty},
 \qquad {\hbox {with}}\qquad
 |f'|_{+\infty}= |\A\Bn^{+}\times \Bn^{+}|/4.
\end{equation}
Finally, we will prove that
\begin{equation}
 \label{p7}
 \lim_{s\rightarrow +\infty}
 e^{i({s^2\over 4} -\gamma_{+,j}\log|s|)}
 (\bar f_jf'_j)(s)= {{b_+}\over {2}}\, e^{-i{a}_{+}},
 \qquad j=1,2.
\end{equation}
To this end, notice that the part {\it {iv)}} in
Theorem~\ref{yh}  asserts that
\begin{equation}
 \label{p8}
 \left(
 {1\over 2} \, {{d|f_j|^2}\over {ds}}-i\Im m(\bar f_j f'_j)
 \right)(s)
 = {{b_{+,j}}\over 2}\, e^{ia_{+,j}} e^{i\phi_{1,j}(s)}+
 O\left( { {1}\over {|s|}} \right),
 \qquad s\rightarrow +\infty,
\end{equation}
where $\phi_{1,j}(s)=(s^2/4)-\gamma_{+,j}\log|s|$, and the
constants $a_{+,j}$, $\gamma_{+,j}$ and $b_{+,j}$ are defined in
Theorem~\ref{yh}.

Let us continue by proving that $a_{+,j}$ and $b_{+,j}$ do not depend
on $j=1,2$ under the limiting conditions (\ref{st1}) and
(\ref{st2}).

Indeed, from the Serret-Frenet equations, (\ref{25}),
(\ref{asym12}), (\ref{p8}) and (\ref{st1}), it is easy to see that
\begin{eqnarray*}
 & &
 \lim_{s\rightarrow +\infty}e^{-i{\tilde\phi}(s)} e^{-\A \log{|s|}}
(\Tn'-i\Tn\times \Tn')(s)=
 \lim_{s\rightarrow +\infty}
 e^{-i{\tilde\phi}(s)}
 e^{-\A\log|s|} c_j(s)\, (\nn_j -i \bn_j)(s)=
    \\
 & &
 \lim_{s\rightarrow +\infty}
 e^{-i({s^2\over 4} + (3|c_j(s)|^2+2\alpha)\log|s|)}
 e^{-\A\log|s|}
 \left(
 {{dc_j^2}\over {ds}} -2i c^2_j (\tau_j-s/2)
 \right)
 (\A\Tn_j\times \Tn_j -i \A\Tn_j)(s)
    \\
 & &
 \qquad \qquad
 = b_{+,j} e^{ia_{+,j}}(\A\Bn^{+}\times \Bn^{+}-i\A\Bn^{+}).
\end{eqnarray*}
Here, it has been used that $ -(3|c_j(s)|^2+2\alpha)-\gamma_{+,j}=
-3(|c_j(s)|^2-c^2_{+\infty,j})=o(1)$, as $s\rightarrow +\infty$,
in obtaining the last identity.

Then, the above identity and (\ref{st2}) yield
$$
  b_{+,j}=b_+
  \qquad {\hbox {and}}\qquad
  a_{+,j}=a_{+},
  \qquad j=1,2.
$$
The identity (\ref{p7}) easily follows by substituing the above
identities into (\ref{p8}).

Now, notice that under the conditions (\ref{p3}),
(\ref{p6}) and (\ref{p7}), Lemma~\ref{unicidad} concludes that
$$
  {{d|f_1|^2}\over {ds}}(s)=  {{d|f_2|^2}\over {ds}}(s)
 \qquad {\hbox {and}}\qquad
 \Im m(\bar f_1 f'_1)= \Im m (\bar f_2 f'_2),
 \qquad s\geq s_0\gg 1.
$$
Hence,
$$
  {{dc_1^2}\over {ds}}(s)=  {{dc_2^2}\over {ds}}(s)
  \qquad {\hbox {and}}\qquad
  c_1^2(\tau_1-s/2)= c^2_2(\tau_2 - s/2),
$$
by taking into account (\ref{25}).

Next recall that $c^2_{+\infty, j}=c^2_{+\infty}\neq 0$
(see~(\ref{p2})). Then, from the above identities we obtain 
$$
   c_1(s)=c_2(s)\neq 0
   \qquad {\hbox {and}}\qquad
   \tau_1(s)=\tau_2(s),
   \qquad \forall\, s\geq s_0\gg 1.
$$
Therefore, there exists $\rho$ rotation in $\R^3$ such that
\begin{equation}
 \label{p8'}
  \Tn^{1}(s)= \rho(\Tn^2(s)),
  \qquad \forall\, s\geq s_0.
\end{equation}
Now observe that, from (\ref{p8'}) and (\ref{st1}), we get that
$$
 \lim_{s\rightarrow +\infty} e^{-\A\log|s|}\rho(\Tn^2(s))=
 \lim_{s\rightarrow +\infty} e^{-\A\log|s|}\Tn^1(s)= \Bn^{+}.
$$
Then, writting $\rho=(a_{i,j})^{3}_{i,j=1}$ from the above
identity follows that
$$
 \lim_{s\rightarrow +\infty}(a_{31}T^1_{1}(s) +
 a_{32}T^1_{2}(s)+ a_{33}T^1_{3}(s))=B^{+}_{3},
$$
where $\Tn^1$ satisfies (\ref{st1}), so that
$a_{31}=a_{32}=0$ and $a_{33}=1$.

As a consequence,
 $
  \rho=\left(
  \begin{array}{cc}
  e^{i\delta}& 0\\
  0 & 1
  \end{array}
  \right),
 $
for some $\delta\in [0,2\pi)$, and
$$
   e^{-\A\log|s|}\rho= \rho\, e^{-\A\log|s|},
   \qquad s\neq 0.
$$
Next, using the above identity and the fact that
$\Tn^{1}=\rho(\Tn^{2})$, from (\ref{st2}) we get that
\begin{equation}
 \label{p9}
 \A\Bn^{+}\times \Bn^{+}-i \A\Bn^{+}=
 \rho (\A\Bn^{+}\times \Bn^{+}-i \A\Bn^{+} ).
\end{equation}
 Denoting as before $(B_1^{+},B_2^{+}, B_3^{+})$ the components of 
 $\Bn^{+}$ and taking into account that
  $
  \rho=
  \left(
  \begin{array}{cc}
  e^{i\delta}& 0\\
  0 & 1
  \end{array}
  \right),
  $
 the first two components of the vectorial identity (\ref{p9}) can
 be rewritten as
 $$
 a\, (B^{+}_{1}+i\,B^{+}_{2})\, (B^{+}_{3}+1)\, (e^{i\delta}-1)=0,
 $$
and, because $B^{+}_{3}\neq \pm 1$ and $a\neq 0$, we conclude that
$e^{i\delta}=1$, $\delta\in[0,2\pi)$, that is $\delta=0$. Then,
$\rho\equiv {\huge {1}}_{3\times 3}$, so that
$$
   \Tn^{1}(s)=\Tn^{2}(s),\qquad \forall\, s\geq s_0\gg 1,
$$
and therefore there exist ${\mathbf {a}}\in\R^3$ such that
$$
 \Gn^{1}(s)=\Gn^2(s)+{\mathbf {a}},
 \qquad \forall\, s\geq s_0\gg 1,
$$
Finally, observe that from (\ref{st1a})
$$
 \Gn^{j}(s)= se^{\A\log|s|}\An^{+}+ o(|s|),
 \qquad {\hbox {as}}\quad
 s\rightarrow +\infty\quad
 {\hbox {and}}\quad j=1,2.
$$
Then, ${\mathbf {a}}=(0,0,0)$ and we conclude that
$$
 \Gn^1(s)=\Gn^2(s),
 \qquad s\geq s_0\gg 1.
$$
\end{proof}
%
 \section{Properties of the curvature and torsion associated
 to the curve $\Gn$}
 \label{stct}
%
In this section, we analyze the properties of the curvature and
torsion functions related to the curve $\Gn$. 
We will assume through this section that $c_{\pm\infty}\neq 0$.

We will now continue to see how the Hasimoto
transform allows to reduce the study of these properties to
analyzing those of the (self-similar) solutions of the cubic
non-linear Schr\"odinger equation.
Though the argument can be found in some classical references as
~\cite{Has} and ~\cite{Lam}, we have included it here for the sake
of completeness.

Firstly, notice that if
\begin{equation}\label{18}
\X(s,t)=e^{{\A\over 2}\log{t}}\sqrt{t}\G(s/\sqrt{t}),
\end{equation}
then the associated curvature and torsion have the self-similar
form
\begin{equation}\label{19}
 c(s,t)={1\over \sqrt{t}}\,c(s/\sqrt{t})
 \qquad {\hbox {and}}\qquad
 \tau(s,t)={1\over \sqrt{t}}\,\tau(s/\sqrt{t}),
\end{equation}
with $c(s)=c(s,1)$ and $\tau(s)=\tau(s,1)$.

Indeed, from (\ref{18}) we get that
$$
 \Tn(s,t)= e^{{\A\over 2}\log{t}}\Tn(s/\sqrt{t}),\qquad
 \nn(s,t)= e^{{\A\over 2}\log{t}}\nn(s/\sqrt{t})
 \qquad {\hbox{and}}
$$
\begin{equation}\label{20}
 \bn(s,t)= e^{{\A\over 2}\log{t}}\bn(s/\sqrt{t}),
\end{equation}
where $\{ \Tn(s),\nn(s),\bn(s) \}$ is the Serret-Frenet frame of
the curve $\G(s)=\X(s,1)$ (i.e. $\Tn(s)=\Tn(s,1)$,
$\nn(s)=\nn(s,1)$ and $\bn(s)=\bn(s,1)$). Then,
 $$
 \Tn_s(s,t)={{1}\over {\sqrt{t}}}
 e^{{\A\over 2}\log{t}}
 \Tn_s(s/\sqrt{t}),
 {\hspace{1cm}}
 \bn_s(s,t)={{1}\over {\sqrt{t}}}
 e^{{\A\over 2}\log{t}}
 \bn_s(s/\sqrt{t})
 $$
and, from the Serret-Frenet system (see (\ref{star})) we conclude
that
\begin{eqnarray*}
 c(s,t)
 &=&\nn(s,t)\cdot\Tn_s(s,t)=
 {1\over \sqrt{t}}(\nn\cdot\Tn_s)(s/\sqrt{t})=
 {1\over \sqrt{t}}c(s/\sqrt{t})
 \qquad {\hbox {and}} \qquad
               \\[2ex]
 \tau(s,t)
 &=&-\nn(s,t)\cdot\bn_s(s,t)=
 -{1\over \sqrt{t}}(\nn\cdot\bn_s)(s/\sqrt{t})=
 {1\over \sqrt{t}}\tau(s/\sqrt{t}).
\end{eqnarray*}
Secondly, defining as in~\cite{Has} (see also~\cite[p.p~195]{Lam})
\begin{eqnarray}\label{21}
\left\{
\begin{array}{ll}
 \psi(s,t)=c(s,t)\, e^{i\int_{0}^{s}\tau (s',t)\, ds'}
       \\[2ex]
 \Nn(s,t)= (\nn+i\bn)(s,t)\,e^{i\int_{0}^{s}\tau(s',t)\, ds'},
\end{array}
\right.
\end{eqnarray}
then $\psi$ solves the nonlinear cubic Schr\"odinger equation
\begin{equation}\label{22}
 i\psi_t+\psi_{ss}+{{\psi}\over {2}} (|\psi|^2+\alpha(t))=0,
\end{equation}
with
$$
 \alpha(t)=
 -i(\Nn_t\cdot {\overline {\Nn}})(0,t)-{{\co}\over {t}},
$$
or, equivalently,
 \begin{equation}\label{23}
     \alpha (t)= ((A\bn\cdot \nn)(0)-\co)/t.
 \end{equation}
Indeed, from the identities in (\ref{19}) and (\ref{20}), and the
definition of the complex vector $\Nn$ in (\ref{21}), it follows
that
$$
 \Nn(s,t)= (\nn+i\bn)(s,t)e^{i\int_{0}^{s}\tau(s',t)\, ds'}
 = e^{{\A\over 2}\log t} (\nn+i\bn)(s/\sqrt{t})
 e^{\int_{0}^{s/\sqrt{t}}\tau(s')\, ds'}.
$$
Now, by deriving this identity with respect to the time variable,
we get that
$$
 {\displaystyle {
 \left.
 (\Nn)_{t}(s,t)= {1\over {2t}}\, e^{{\A\over 2}\log t}
 (\A(\nn+i\bn)(\eta) -\eta (\nn+i\bn)_{t}(\eta)-i \eta (\nn
 +i\bn)(\eta)\,\tau(\eta))\,
 e^{i\int_{0}^{\eta}\tau(s')\, ds'} \right|_{\eta=s/\sqrt{t}}.
 }}
$$
Then, evaluating these identities at $s=0$,
$$
  \Nn(0,t)= e^{{\A\over 2}\log t} (\nn+i\bn)(0),
  \qquad
  {\Nn}_{t}(0,t)= {1\over {2t}}
  e^{{\A\over {2}}\log t} \A(\nn+i\bn)(0),
$$
and therefore,
\begin{eqnarray*}
  (\Nn_{t}\cdot{\overline{\Nn}})(0,t)
  &=&
  {1\over {2t}}\,
  e^{{\A\over {2}}\log t} \A(\nn+i\bn)(0)
  \cdot
  e^{{\A\over {2}}\log t} (\nn-i\bn)(0)
     \\
  &=&
  {1\over {2t}}\,
  \A(\nn+i\bn)(0) \cdot (\nn-i\bn)(0)
  = {i\over t} \,\A\bn\cdot \nn(0).
\end{eqnarray*}
Here, we have used that $\A$ is an antisymmetric matrix, so that
$\A \vn\cdot \vn=0$, and $e^{{\A\over 2}\log t}$ is a rotation in
$\R^3$.

Next, notice that, because  $c(s,t)$ and $\tau(s,t)$ satisfy
(\ref{19}), the function $\psi$ can be written as
$$
 \psi(s,t)={{1}\over {\sqrt{t}}}\psi(s/\sqrt{t}),
$$
where
\begin{equation}\label{psi}
 \psi(s)=c(s) e^{i\int_{0}^{s}\tau(s')\, ds'},
\end{equation}
and, from the equation (\ref{22}) and (\ref{23}), it follows that
$\psi(s)$ solves
\begin{equation}\label{eqpsi}
 \psi''-{i\over 2}(\psi+s\psi')+{\psi\over 2}
 (|\psi|^2+\alpha)=0,
\end{equation}
with $\alpha=\alpha(1)$, that is
\begin{equation}
 \label{23'}
 \alpha= \A \bn\cdot \nn(0) -\co.
\end{equation}
Finally, if we introduce the function $f$ through the definition
\begin{equation}\label{f}
 \psi(s)=f(s) e^{is^2/4},
\end{equation}
then $f$ satisfies
$$
 f''+i{s\over 2}f'+{f\over 2}(|f|^2+\alpha)=0.
$$
On the one hand, recall that from the definitions of $\psi$ and $f$ in
(\ref{psi}) and (\ref{f}), it follows that
\begin{equation}\label{24}
 f(s)e^{is^2/4}= c(s)e^{\int_{0}^{s}\tau(s')\, ds'},
\end{equation}
so that
\begin{equation}\label{25}
 |f|^2=c^2, \qquad
 |f'|^2=c_{s}^2 + c^2(\tau- s/2)^2 \qquad
 {\hbox {and}}\qquad
 \Im m(\bar f f')=c^2(\tau-s/2).
\end{equation}
On the other hand, we have already seen that if $\Xn(s,t)$ of the form
(\ref{18}) is a solution of LIA, then $\Tn=\Gn'$ satisfies-see Lemma
\ref{for} and (\ref{asym3})-
$$
 \A \Tn= c(s-2\tau)\nn+2c_{s}\bn
 \qquad {\hbox {and}}\qquad
 c^2=-aT_3(s)-\alpha,
 \footnote{  Define ${\mathbf {v}}=2c_0\bn(0)$, 
  so that $|\vn|^2=4c_{0}^{2}$. From 
 the definiton of $\alpha$ in (\ref{23'}) and the evaluation of
 (\ref{12}) at $s=0$, it is easy to see that
 $
 \alpha=
 \Tn(0)\cdot 
 \left[
   \A\left({{\mathbf{v}}\over {|{\mathbf{v}}|}} \right)\times
   {{\mathbf{v}}\over {|{\mathbf{v}}|}}
 \right]
 -{{|{\mathbf{v}}|^2}\over{4}}=
 -aT_{3}(0)-c_{0}^{2}=
 -aT_3(0)-{1\over 4}\, |(\I+\A)\Gn(0)|^2.
 $
 }
$$
from which the above identities rewrite as
\begin{equation}\label{25'}
 |f|^2=-aT_3(s)-\alpha,
 \qquad
 |f'|^2={1\over 4}\, |\A\Tn \times \Tn|^2 \qquad
 {\hbox {and}} \qquad
 \Im m(\bar f f')= -{1\over 2}\, \A\Tn\cdot \Tn'.
\end{equation}
Finally, recall the definitions of the pair $(y,h)$ in (\ref{in15})
and (\ref{in15'}):
\begin{equation}\label{25''}
  y={{d|\Tn'|^2}\over {ds}},
  \qquad 
  h=-{1\over 2}\,\A \Tn\cdot \Tn'
  \qquad {\hbox{and}}\qquad
  {y\over 2}+ih=\bar f f'.
\end{equation}
The identities in (\ref{25}) reduce the problem of studying 
the properties of $c(s)$ and $\tau(s)$  
to analize the properties related to the solutions $f$ of the latter equation. 
In fact, the following theorem is an immediate
consequence of (\ref{25}) and the results that will be proved in
the next section for the solutions of the latter O.D.E. (see
Theorem ~\ref{yh}).

\begin{theorem}\label{ct}\
Let $\Xn(s,t)=e^{{\A\over 2}\log t}\sqrt{t}\Gn(s/{\sqrt{t}})$ and
let $c(s)$ and $\tau(s)$ denote the curvature and torsion related
to the curve $\Gn(s)=\Xn(s,1)$. Then,
\begin{itemize}
 \item[{\it i)}] There exists 
  $  \footnote{ Notice that, from (\ref{25}) and (\ref{25'}), it follows
  that $E(0)=a^2/4$.}E(0)\geq 0$ 
  such that the identity
  $$
  c_s^2+c^2(\tau-s/2)^2+{1\over 4}(c^2+\alpha)^2=E(0)
  $$
  holds true for all $s\in\R$.
 \item[{\it ii)}] $c(s)$, $c_s(s)$ and $\tau -s/2$ are bounded
 globally defined functions.
 \item[{\it iii)}] The limits
 $\lim_{s\rightarrow\pm\infty}c(s)=c_{\pm\infty}$do  exist and
  $$
  c^2(s)-c^2_{\pm\infty}=O({1/|s|}),
  \qquad s\rightarrow\pm\infty.
  $$
 \item[{\it iv)}] Moreover, the following asymptotics hold as
 $s\rightarrow\pm\infty$:
  \begin{eqnarray*}
  c^2(s)
  &=&
  c^2_{\pm\infty}-{4\over s}\,c^2(\tau-s/2)-2
  {{{{\tg}}_{\pm}}\over {s^2}}+O\left({1\over |s|^3}
  \right),
     \\[2ex]
  \left( {1\over 2}\, {{dc^2}\over {ds}}+i c^2(\tau-{s\over 2})\right)(s)
  &=&
   {b_\pm\over 2}\, e^{-ia_\pm} e^{-i\phi(s)} -i\, {{\tg_\pm}\over s}+
   i\,{b_\pm\over 2}\,
   \left( 3\tg_{\pm} - {{\gamma^2_{\pm}}\over {2}} \right)\,
   {{e^{-ia_\pm}}\over {s^2}}\, {e^{-i\phi(s)}}
      \\[2ex]
  & &
    \qquad
   - \,{{b_\pm\gamma_\pm}\over 2}\,
   {{e^{ia_\pm}}\over {s^2}}\, e^{i\phi(s)}
   + \,{{b_\pm\gamma_\pm}\over 2}\,
   {{e^{-ia_\pm}}\over {s^2}}\, e^{-i\phi(s)}
   +O\left( {{1}\over {|s|^3}}  \right).
 \end{eqnarray*}
 \end{itemize}
 Here, $\phi(s)= (s^2/4) -\gamma_{\pm}\, \log |s|$,
 $a_{\pm}$ is a constant in the interval $[0,2\pi)$,
  \begin{eqnarray}\label{ct1}
  \gamma_{\pm}
  &=&
  -3 \,c^2_{\pm\infty} - 2\alpha,
  \qquad\qquad
  {{\tg}}_{\pm}
  = 2E(0)-(3c^{2}_{\pm\infty}+\alpha)\,(c^{2}_{\pm\infty}+\alpha)/2,
      \nonumber \\[1ex]
  b_{\pm}^2
  &=&
  4c^2_{\pm\infty}(E(0)-(c^2_{\pm\infty}+\alpha)^2/4)
  \qquad {\hbox {with}}\qquad
   b_{\pm}\geq 0,
      \nonumber \\[1ex]
  \alpha
  &=&
  -aT_3(0)-{1\over 4}\, |(\I+\A)\Gn(0)|^2
  \qquad {\hbox {and}}\qquad E(0)={a^2\over 4}.
 \end{eqnarray}
\end{theorem}

\section{Self-similar solutions of the Cubic Non-linear
Schr\"odinger equation} \label{selfNLS}
%
\begin{theorem}\label{yh}\
Let $f$ be a solution of the equation
\begin{equation}\label{eqf}
 f''+i{s\over 2}f'+{f\over 2}(|f|^2+\alpha)=0,
 \qquad \alpha\in \R.
\end{equation}
Then
\begin{itemize}

\item[{\it i)}] There exists $E(0)\geq 0$ such that the identity
 $$
 |f'|^2+{1\over 4}(|f|^2+\alpha)^2=E(0)
 $$
 holds true for all $s\in\R$.
\item[{\it ii)}] $f(s)$, $f'(s)$ and $\I m(\bar f f')(s)$ are bounded
globally defined functions.
\item[{\it iii)}] The limits
$\lim_{s\rightarrow \pm\infty}|f|^2(s)=|f|^2_{\pm\infty}$ and
$\lim_{s\rightarrow \pm\infty}|f'|^2(s)=|f'|^2_{\pm\infty}$ do
exist and
$$
 |f|^2 -|f|^{2}_{\pm\infty}= O\left({{1}\over {|s|}}  \right),
 \qquad s\rightarrow \pm \infty.
$$
\item[{\it iv)}] The following asymptotics hold as $s\rightarrow \pm\infty$:
 \begin{eqnarray*}
 |f|^2(s)
 &=&
 |f|^2_{\pm\infty}-{4\over s}\,\Im m (\bar f f')(s)-
 2 {{{{\tg}}_{\pm}}\over {s^2}}
 +O\left({1\over |s|^3}\right),
    \\[2ex]
 \bar f f'(s)
 &=&
  \left( {1\over 2}\, {{d|f|^2}\over {ds}}+i \Im m(\bar f f')\right)(s)
    \\[1ex]
 &=&
  {b_\pm\over 2}\, e^{-ia_\pm} e^{-i\pp} -i\, {{\tg_\pm}\over s}+
  i\,{b_\pm\over 2}\,
  \left( 3\tg_{\pm} - {{\gamma^2_{\pm}}\over {2}} \right)\,
  {{e^{-ia_\pm}}\over {s^2}}\, {e^{-i\pp}}
     \\[2ex]
 & &
   \qquad
  - \,{{b_\pm\gamma_\pm}\over 2}\,
  {{e^{ia_\pm}}\over {s^2}}\, e^{i\pp}
  + \,{{b_\pm\gamma_\pm}\over 2}\,
  {{e^{-ia_\pm}}\over {s^2}}\, e^{-i\pp}
  +O\left( {{1}\over {|s|^3}}  \right).
 \end{eqnarray*}
 \item[{\it {v)}}] Moreover, if $|f|_{+\infty}\neq 0$ or
 $|f|_{-\infty}\neq 0$, then
 \begin{eqnarray*}
  f(s)
  &=&
  |f|_{\pm\infty}\, e^{ic_\pm}\, e^{i\ps}
  + 2i\, |f'|_{\pm\infty} \,
  {{e^{id_\pm}}\over {s}}\, e^{i\pt}
  - |f|_{\pm\infty} (|f|^2_{\pm\infty} +\alpha)
  {{e^{ic_\pm}}\over {s^2}}\, e^{i\ps} + O\left( {1\over {|s|^3}} \right),
       \\
  f'(s)
  &=&
  |f'|_{\pm\infty} e^{id_{\pm}} e^{i\phi_3(s)}
  + i |f|_{\pm\infty}(|f|^2_{\pm\infty}+\alpha)\,
  {{e^{ic_\pm}}\over {s}}\, e^{i\phi_2(s)}
       \\
  & &
  \qquad
  + |f'|_{\pm\infty}
  \left(
  (|f|^2_{\pm\infty}+\alpha)
  +i(3\tg_\pm -{\gamma_\pm^2\over 2}+ \gamma_\pm)
  \right)\,
  {{e^{id_\pm}}\over {s^2}}\, e^{i\phi_3(s)}
       \\
  & &
  \qquad \qquad
  -2|f'|_{\pm\infty}(|f|^2_{\pm\infty}+\alpha)\,
  {{e^{i(c_\pm+ a_\pm)}}\over {s^2}}\,
  e^{i(\phi_2+\phi_1)(s)}
  +O\left( {{1}\over {|s|^3}} \right).
 \end{eqnarray*}
\end{itemize}
 Here,
 $$
  \left\{
  \begin{array}{l}
  \pp= (s^2/4) -\gamma_{\pm}\, \log |s|,
      \\[1ex]
  \ps= (|f|^2_{\pm\infty}+\alpha)\log |s|,
      \\[1ex]
  \pt=
  -(s^2/ 4)- (2|f|^2_{\pm\infty}+\alpha)\log |s|,
  \qquad (i.e.\  \phi_3=\phi_2-\phi_1)
      \\[1ex]
 \end{array}
 \right.
 $$
 $|f|_{\pm\infty}$, $|f'|_{\pm\infty}\geq 0$,
 $a_{\pm}$ and $\ c_{\pm}$ are arbitrary
 constants in $[0,2\pi)$, $d_{\pm}=c_{\pm}-a_{\pm}$,
 \begin{eqnarray}\label{yh1}
  \gamma_{\pm}
  &=& -3 \,|f|^2_{\pm\infty} - 2\alpha,
  \qquad 
  \tg_{\pm}= 2E(0)- (3|f|^2_{\pm\infty}+\alpha)\, 
  (|f|^2_{\pm\infty} +\alpha)/2,
   \quad {\hbox {and}}
    \nonumber \\
  b^2_{\pm}
  &=&
  4 |f|^2_{\pm\infty}\,
  (E(0)-(|f|^2_{\pm\infty}+\alpha)^2/4)
  \qquad {\hbox {with}}\qquad
  b_{\pm}\geq 0.
 \end{eqnarray}
\end{theorem}
\begin{remark}
 \label{reyh0}
 The maps $(f(0), f'(0))\longrightarrow (|f|_{\pm\infty}, c_{\pm},
 |f'|_{\pm\infty}, a_{\pm})$ are continuous. This follows from the
 construction of the solution. In particular, from (\ref{cl3}),
 (\ref{cl12}), (\ref{cl23}), and
 from the convergence of the integrals (\ref{cl26})-(\ref{cl26'}).
\end{remark}
\begin{proof}[Proof of Theorem~\ref{yh}]
Firstly, by multiplying the equation (\ref{eqf}) by $\bar f'$ and
taking the real part, it is easy to see that
$$
{{dE}\over {ds}}= {{d }\over {ds}} [|f'|^2+{1\over
4}(|f|^2+\alpha)^2]=0,
$$
so that the following quantity is preserved for all $s$
\begin{equation}\label{cl}
 |f'|^2+{1\over 4}(|f|^2+\alpha)^2=E(0).
\end{equation}
As a consequence, we get that $f$ and $f'$ are globally
well-defined,
$$
 |f(s)|,\ |f'(s)|\leq C
 \qquad {\hbox {and}} \qquad
 |\Im m(\bar f f')(s)|\leq C,\qquad \forall s\in\R.
$$
This concludes {\it {i)}} and {\it {ii)}}.
We will now continue with the proof of {\it {iii)}}-{\it {v)}} in
the case $s\rightarrow +\infty$. The case $s\rightarrow -\infty$
follows using the same arguments. 
To this end, we consider the functions $h(s)$ and $y(s)$ defined
through the following identities:
\begin{equation}\label{cl1}
 h=\Im m(\bar f f')
 \qquad {\hbox {and}}\qquad
 y={{d|f|^2}\over {ds}}.
\end{equation}
On the other hand, since $f$ solves  (\ref{eqf}), we obtain 
$$
 h'(s)=\Im m(\bar f' f'+\bar f f'')=-{s\over 4}{{d|f|^2}\over {ds}}
    =-{s\over 4} y(s),
$$
and
$$
 y'=2\, |f'|^2 + 2\Re e (f''\bar f)= sh+ 2|f'|^2 -|f|^2 \,
 (|f|^2+\alpha)=
 sh+g(|f|^2),
$$
where $g(|f|^2)= 2E(0)-(3|f|^3+\alpha)(|f|^2+\alpha)/2$, by using the
conservation law (\ref{cl}).

Therefore, we conclude that the pair $(y,h)$  satisfies the coupled
system of equations:
 \begin{eqnarray}\label{cl9}
 \left\{
 \begin{array}{ll}
  \displaystyle{{{d|f|^2}\over {ds}}}=y &\qquad
        \\[2ex]
  y'=sh+g(|f|^2);  &\qquad
  g(|f|^2)=2E(0)-(3|f|^2+\alpha)(|f|^2+\alpha)/2
        \\[2ex]
  h'=-\displaystyle{{s\over 4} }y. &
 \end{array}
 \right.
 \end{eqnarray}
Also, since $f$ and $f'$ are bounded,  $h(s)$ and $y(s)$
are bounded functions. Notice also that multiplication by $\bar f$
in (\ref{eqf}) yields
$$
 f''\bar f +i {s\over 2} f'\bar f +
 {{|f|^2}\over {2}} (|f|^2+\alpha)=0,
$$
from which we obtain 
$$
 \Im m(f''\bar f)(s) +{s\over 4} {{d |f|^2}\over {ds}}=0.
$$
Then,
$$
 \int_{0}^{s} \Im m(f''\bar f)(s')\, ds' +
 {1\over 4}\int_{0}^{s} s'{{d|f|^2 }\over {ds'}}(s')\, ds'=0,
$$
and integration by parts gives
\begin{equation}\label{cl2}
 |f|^2-{1\over s}\int_{0}^{s}|f|^2=
 -{4\over s} (h(s)-h(0)),
\end{equation}
so that
\begin{equation}\label{cl3}
  \left({1\over s}\int_{0}^{s}|f|^2  \right)_s=
 + {1\over s}\left( |f|^2 -{1\over s}\int_{0}^{s}|f|^2  \right)=
 -{4\over {s^2}}(h(s)-h(0)),\qquad s\neq 0.
\end{equation}
Since  $h$ is bounded, from the
above identity we get that the limit $\lim_{s\rightarrow +
\infty}{(1/s)}\int_{0}^{s}|f|^2$ exists and, from (\ref{cl2}), it
follows that
$$
 \lim_{s\rightarrow + \infty}
 {1\over s}\int_{0}^{s}|f|^2 = |f|^2_{+\infty}.
$$
Notice  that, from the conservation law (\ref{cl}) and the above
identity, we also obtain the existence of the limit
$\lim_{s\rightarrow +\infty}|f'(s)|=|f'|_{+\infty}$.

Now the integration of (\ref{cl3}) from $s>0$ to $+\infty$ yields
\begin{equation}\label{cl4}
 |f|^2_{+\infty}- {1\over s}\int_{0}^{s} |f|^2(s')\, ds'=
 4{h(0)\over s} - 4 \int_{s}^{+\infty} {h(s')\over {(s')^2}}\,
 ds',
 \qquad s>0.
\end{equation}
From (\ref{cl2}) and (\ref{cl4}), we get
\begin{equation}\label{cl5'}
 |f|^2(s)-|f|^2_{+\infty}=
 -4 {h\over s} + 4 \int_{s}^{+\infty} {{h(s')}\over {(s')^2}}\,
 ds',
\end{equation}
and, by taking into account that $h$ is bounded, we conclude {\it
  {(iii)}}, that is
\begin{equation}\label{cl5}
 |f|^2(s)-|f|^2_{+\infty}=O\left( {{1}\over {|s|}} \right),
 \qquad s\rightarrow +\infty.
\end{equation}
We now continue to prove the asymptotics in {\it (iv)}.

Defining  ${{\tg}}_{+}$ to be the limiting value 
of $g(|f|^2(s))$ as $s\rightarrow +\infty$, that is
\begin{equation}\label{cl10}
 {{\tg}}_{+} =
 2E(0)-{1\over 2}\,(3 |f|^{2}_{+\infty}+\alpha)\,( |f|^{2}_{+\infty}+\alpha)
 =2|f'|^2_{+\infty} -|f|^2_{+\infty}(|f|^2_{+\infty}+\alpha),
\end{equation}
(recall the conservation law (\ref{cl})) and taking into account
that $h=y'/s-g(|f|^2(s))$ (see~(\ref{cl9})), it can be shown that
\begin{eqnarray}
 \label{cl11}
 & &
 |f|^2(s)-|f|^{2}_{+\infty}
 = -4 {h\over s}+
  4\int_{s}^{+\infty} {{h(s')}\over {(s')^2}}\,ds'=
  -4 {h\over s} + 4 \left.{{y(s')}\over {(s')^3}}\right|_{s}^{+\infty}
     \\[1ex]
 & &
 \qquad\qquad
 +12\int_{s}^{+\infty} {{y(s')}\over {(s')^4}}\,ds'
 - 2 {{{{\tg}}_{+}}\over {s^2}}
 -4\int_{s}^{+\infty} {{g(|f|^2(s'))-{{\tg}}_{+}}\over {(s')^3}}\,
 ds',
 \qquad s>0. 
    \nonumber 
\end{eqnarray}
Also, defining $\gamma_+$ to be the limiting value of
$g'(|f|^2(s))$ as $s\rightarrow +\infty$,
\begin{equation}\label{cl12}
 \gamma_{+}=g'(|f|^2_{+\infty})=-3|f|^{2}_{+\infty}-2\alpha,
\end{equation}
from the definition of $g(|f|^2)$ in (\ref{cl9}) and (\ref{cl5}),
it follows that
$$
 g(|f|^2)-{{\tg}}_{+}={\gamma_+}(|f|^2(s)-|f|^2_{+\infty})
 -{3\over 2}(|f|^2(s)-|f|^2_{+\infty})^2
 =O\left({{1}\over {|s|}}  \right),
 \qquad {\hbox {as}}\qquad s\rightarrow +\infty.
$$
Therefore, taking into account the above observations in
(\ref{cl11}), we conclude that
\begin{equation}\label{cl13}
 |f|^2(s)=|f|^{2}_{+\infty}-4 {h\over s}
 -2 {{{{\tg}}_+}\over {s^2}}
 +O\left({{1}\over {|s|^3}}  \right),
 \qquad {\hbox {as}}\qquad s\rightarrow +\infty.
\end{equation}
This finishes the proof of the asymptotic development related to
$|f|^2(s)$ in {\it {iv)}}.
In order to derive the behaviour of $y(s)$ and $h(s)$ for $s$
sufficiently large, we will consider in (\ref{cl9}) two new
variables $u$ and $v$ defined through $y(s)$ and $h(s)$ as follows
\begin{equation}\label{cl14}
 y(s)=u(s^2/4)
 \qquad {\hbox {and}} \qquad
 h(s)=v(s^2/4).
\end{equation}
Then, defining $t=s^2/4$, i.e., $|s|=2\sqrt{t}$, from (\ref{cl9})
and (\ref{cl14}), we obtain that $u$ and $v$ satisfy:
\begin{equation}\label{cl15}
 u'(t)= 2v(t)+{{g(|f|^2(2\rt))}\over {\rt}},
 \qquad {\hbox {and}}\qquad
 v'(t)=-{1\over 2}u(t).
\end{equation}
Therefore,
\begin{equation}\label{cl17}
 u''=-u+{\gamma_+\over t} u
     -{{\tg_+}\over {2{t}^{3/2}}}
     + {{(g')(|f|^2(2\rt))-\gamma_+}\over {t}}\, u
     -{{g(|f|^2(2\rt))-\tg_+}\over {2{t}^{3/2}}},
\end{equation}
Observe that, because $f$ and $f'$ are bounded, $|f|^2$, $y$ and
$h$ are bounded. Then, from (\ref{cl14}) and (\ref{cl15}), it
follows that
\begin{equation}\label{cl22'}
 |u(t)|\leq C
 \qquad {\hbox {and}}\qquad
 |u'(t)|\leq C,
 \qquad \forall\, t>1,
 \qquad {\hbox {and}}
\end{equation}
\begin{equation}\label{cl16'}
 h(2\rt)=v(t)
 ={{u'}\over 2}-{{{\tg}_+}\over {2\rt}}-
 {{g(|f|^2(2\rt))-\tg_+}\over {2\rt}}
 ={{u'}\over 2}+O\left( {1\over \rt} \right),
 \qquad t\rightarrow +\infty.
\end{equation}
Next, notice that (\ref{cl17}) rewrites as
\begin{eqnarray} \label{cl18'}
 {\left(
 \begin{array}{l}
  u\\
  u'
 \end{array}
 \right)}'
 &=&
 \left(
 \begin{array}{ll}
  0& 1  \\
  -1+{\gamma_+\over t}& 0
 \end{array}
 \right)
 \left(
 \begin{array}{l}
 u\\
 u'
 \end{array}
 \right)
 +
 \left(
 \begin{array}{l}
 0\\
 F_2
 \end{array}
 \right),
 \qquad {\hbox {with}}
\end{eqnarray}
$$
  F_2(t)= -{{\tg_+}\over {2{t}^{3/2}}}
     + {{(g')(|f|^2(2\rt))-\gamma_+}\over {t}}\, u
     -{{g(|f|^2(2\rt))-\tg_+}\over {2{t}^{3/2}}}.
$$
From the definition of $g$, $\gamma_{+}$ and $\tg_{+}$ (see
(\ref{cl9}), (\ref{cl10}) and (\ref{cl12}), respectively), we get
that
\begin{eqnarray}
 & &
 g(|f|^2(2\rt))-\tg_+=\gamma_+ (|f|^2(2\rt)-|f|^2_{+\infty})
 -{3\over 2}\, (|f|^2(2\rt)-|f|^2_{+\infty})^2,
 \qquad {\hbox {and}}
     \nonumber \\[2ex]
 & &
 (g')(|f|^2(2\rt))- \gamma_+= -3\, (|f|^2(2\rt)-|f|^2_{+\infty}).
 \label{cl18}
\end{eqnarray}
where, using (\ref{cl13}) and (\ref{cl16'}),
$$
 |f|^2(2\rt)-|f|^2_{+\infty}
 = -{{u'}\over {\rt}}+ {{\tg_+}\over {2t}}
   +{{g(|f|^2(2\rt))-\tg_+}\over {t}}
   +O\left( {{1}\over {{t}^{3/2}}} \right).
$$
Substituting the above identity into  expression (\ref{cl18}),
after some straightforward calculations we find that
\begin{eqnarray}
 \label{cl20}
 g(|f|^2(2\rt))-\tg_+
 &=&
 -\gamma_+ {{u'(t)}\over {\rt}}+{{\gamma_+\tg_+}\over {2t}}
 -{3\over 2}\, {{(u')^2(t)}\over {t}}
 +O\left( {{1}\over {{t}^{3/2}}} \right),
 \quad t\rightarrow +\infty,
     \\
 (g')(|f|^2(2\rt))-\gamma_+
 &=&
 3{{u'(t)}\over {\rt}}- {3\over 2}\,{{\tg_+}\over t}
 +O\left( {{1}\over {{t}^{3/2}}} \right),
     \nonumber
\end{eqnarray}
so that
\begin{equation}
 \label{cl22}
 F_2(t)=-{{\tg_+}\over {2{t}^{3/2}}}
        + 3\, {{u\, u'}\over {t^{3/2}}}
        + {\gamma_+\over 2} {{u'}\over {t^2}}
        -{{3\tg_+}\over {2}}\, {u\over {t^2}}
       +O\left( {{1}\over {{t}^{5/2}}} \right),
 \qquad {\hbox {as}}\qquad
 t\rightarrow +\infty.
\end{equation}
The diagonalisation of the matrix defining (\ref{cl18'}), and the
change of variables:
\begin{equation}\label{cl23}
  \left(
  \begin{array}{l}
  u\\
  u'
  \end{array}
  \right)= P
  \left(
  \begin{array}{l}
  \twp\\
  \tws
  \end{array}
 \right),
 \qquad {\hbox {with}}\qquad
 P= \left(
 \begin{array}{ll}
  1& 1\\
  i\lambda_{+} & -i\lambda_{+}
 \end{array}
 \right);
 \qquad \lambda_{+}=\sqrt{1-\gamma_{+}/t}
\end{equation}
give that the new variables $\twp$ and $\tws$ satisfy
\begin{eqnarray*}
  & &
  {\left(
  \left(
  \begin{array}{ll}
  {e}^{-i\int_{1}^{t}\lambda_{+}}& 0\\
  0& {e}^{-i\int_{1}^{t}\lambda_{+}}
  \end{array}
  \right)
  \left(
  \begin{array}{l}
  \twp\\
  \tws
  \end{array}
  \right)
 \right)}^{\prime}(t)=
     \\[3ex]
 & & \qquad \qquad
 \left(
  \begin{array}{ll}
  {e}^{-i\int_{1}^{t}\lambda_{+}}& 0\\
  0& {e}^{i\int_{1}^{t}\lambda_{+}}
  \end{array}
 \right)
 \left\{
 -{i\over {2\lambda_{+}}}
 \left(
  \begin{array}{l}
  F_2\\
  -F_2
  \end{array}
 \right)
 + {{\gamma_+}\over {4t^2(1-{\gamma_+/ t})}}
 \left(
  \begin{array}{l}
  \tws-\twp\\
  -(\tws-\twp)
  \end{array}
 \right)
 \right\}.
\end{eqnarray*}
Notice that,  because $u(t)$ and $u'(t)$ are real valued
functions, from the definition of $\twp(t)$ and $\tws$ in
(\ref{cl23}), it is easy to see that
\begin{equation}\label{cl24}
 \tws=\overline{\twp}.
\end{equation}
As a byproduct, in order to analyze the solution of the latter
system of equations, we can reduce ourselves to study its first
component, that is
$$
 {\left(
 {e}^{-i\int_{1}^{t}\lambda_{+}} \twp
 \right)}^{'}(t)=
 {\displaystyle {
 {e}^{-i\int_{1}^{t}\lambda_{+}}
 }}
 \left(
 - {i\over {2\lambda_{+}}}\, F_2(t)
 + {{\gamma_+}\over {4t^2(1-\gamma_+/t)}}
 (\overline{\twp}-\twp)(t)
 \right)=
 {\displaystyle {
 {e}^{-i\int_{1}^{t}\lambda_{+}}
 }}\, I(t).
$$
To this end, we integrate the previous identity from $1$ to $t\gg
1$ to obtain
\begin{equation}\label{cl25}
 {e}^{-i\int_{1}^{t}\lambda_{+}} \twp (t)=
 \twp(1)+\int_{1}^{t} {e}^{-i\int_{1}^{t'}\lambda_{+}}
 \, I(t')\, dt'
\end{equation}
Notice that (\ref{cl23}) and (\ref{cl24}) yield
$$
 \left(
  \begin{array}{l}
   \twp\\
   \overline{\twp}
  \end{array}
 \right)=
 P^{-1}
  \left(
  \begin{array}{l}
   u\\
   u'
  \end{array}
 \right)
  = {1\over 2}
  \left(
  \begin{array}{l}
   u-{i\over \lambda_{+}}u'\\
   u+{i\over \lambda_{+}}u'
   \end{array}
 \right),
 \quad {\hbox{so that}}\quad
 \overline{\twp} -\twp={i\over \lambda_{+}}u'.
$$
Then, from (\ref{cl22}) and the above observation, we get that
 \begin{equation}\label{cl25'}
 I(t)=
 {i\over 2\lambda_{+}}
 \left(
 {{\tg_+}\over {2{t}^{3/2}}}
 - 3\, {{u\, u'}\over {t^{3/2}}}
 +{{3\tg_+}\over {2}}\, {u\over {t^2}}
 \right)+O\left( {{1}\over {{t}^{5/2}}} \right).
 \end{equation}
 Here we have used the asymptotic behaviour of
 ${(1-\gamma_+/t)}^{-1}$
 when $t\gg 1$, and the fact that both $u$ and $u'$ are bounded
 for $t$ sufficiently large (see Remark in (\ref{cl22'})).

 Therefore, from (\ref{cl25}) and (\ref{cl25'}), we obtain 
 \begin{equation}\label{cl26}
 {e}^{-i\int_{1}^{t}\lambda_{+}} \twp(t)= z_+ -{i\over 2}
 \int_{t}^{+\infty}
 {{{e}^{-i\int_{1}^{t'}\lambda_{+}}}\over {\lambda_{+}}}
 \left\{
 { {\tg_+}\over  {2{(t')}^{3/2}} }
 - 3 { {u\,u'}\over {{(t')}^{3/2}}  }
 +{3\over 2} \tg_+ {u\over {{(t')}^2}}
 \right\}\, dt'
 +O\left( {{1}\over {{t}^{3/2}}} \right),
 \end{equation}
 with 
 $z_{+}=\lim_{t\rightarrow +\infty} {e}^{-i\int_{1}^{t}\lambda_{+}}
 \twp(t)$,
 that is,
 \begin{equation}\label{cl26'}
 z_{+}=
 \twp(1)+\int_{1}^{+\infty} {e}^{-i\int_{1}^{t'}\lambda_{+}}
  \left(
  -{i\over {2\lambda_{+}}} F_2(t')+
  {{\gamma_{+}}\over {4 (t')^2(1-\gamma_+/t')}}
  (\overline{\twp}-\twp)(t')\, dt'
  \right).
 \end{equation}
 In particular, 
 \begin{equation}\label{cl27}
  \twp(t)= z_+ {e}^{i\int_{1}^{t}\lambda_{+}}
  (1 + O( {1/ \rt} )),
  \qquad {\hbox {as}}\qquad t\rightarrow +\infty.
 \end{equation}
 Now, recall that $\lambda_{+}={(1-\gamma_+/t)}^{-1/2}= 1+ O(1/t)$, as 
 $t\rightarrow +\infty$. Also, from (\ref{cl23}) and the fact that
 $\tws=\overline{\twp}$, we obtain 
 \begin{equation}\label{cl27'}
 \left(
  \begin{array}{l}
  u\\
  u'
  \end{array}
 \right)= P
 \left(
  \begin{array}{l}
  \twp\\
  \tws
  \end{array}
 \right)=
 \left(
  \begin{array}{l}
  \twp+\overline{\twp}\\
  i\lambda_{+} (\twp-\overline{\twp})
  \end{array}
 \right).
 \end{equation}
 Then, from (\ref{cl27}), it follows that
\begin{eqnarray*}
 u
 &=&
 \twp+\overline{\twp}= z_+ e^{i\int_{1}^{t}\lambda_{+}} +
 \overline{z_+} e^{-i\int_{1}^{t}\lambda_{+}}+
 O( {{1}/ {\rt}} ),
 \quad \quad {\hbox {and}}
     \\
 u\,u'
 &=&
 i\lambda\, ( z_+^2 {e}^{2i\int_{1}^{t}\lambda_{+}}
 -(\overline{z_+})^2 {e}^{-2i\int_{1}^{t}\lambda_{+}}
 +
 O( {{1}/ {\rt}})).
\end{eqnarray*}
 On the one hand, if we write
 $e^{-i\int_{1}^{t'}\lambda_{+}}=
 (e^{-i\int_{1}^{t'}\lambda_{+}})_s/(-i\lambda_{+})$,
 an integration by parts argument shows that
\begin{equation} \label{cl27''}
 \int_{t}^{+\infty} {{{e}^{-i\int_{1}^{t'}\lambda_{+}}}\over {\lambda_{+}}}
 \, {{dt'}\over {(t')^{3/2}}}
 = \int_{t}^{+\infty}
 {{{e}^{-i\int_{1}^{t'}\lambda_{+}}}\over {(t')^{3/2}}}
 \, dt'
 + O\left( {1\over t^{3/2}}  \right)
 =O\left( {{1}\over {t^{3/2}}}  \right).
\end{equation}
On the other hand, a similar integration by parts argument in the
expression of $uu'(t)$ yields
\begin{eqnarray*}
 & &
 \int_{t}^{+\infty} {{{e}^{-i\int_{1}^{t'}\lambda_{+}}}\over
 {\lambda_{+}}}\,
 {{uu'}\over {(t')^{3/2}}}\, dt'=
 O\left({1\over t}\right),
 \qquad {\hbox {as}}\qquad
 t\rightarrow +\infty.
\end{eqnarray*}
Hence, from (\ref{cl26}), we get that
$$
 \twp(t)= z_+\, {e}^{i\int_{1}^{t}\lambda_{+}}
 (1+ O(1/ t) ),
 \qquad {\hbox {as}}\qquad
 t\rightarrow +\infty.
$$
Next, arguing with this new expression for $\twp(t)$ instead of
with (\ref{cl27}), the same ideas as before show that
$$
 \int_{t}^{+\infty} {{{e}^{-i\int_{1}^{t'}\lambda_{+}}}\over
 {\lambda_{+}}}\,
 {{uu'}\over {(t')^{3/2}}}\, dt'=
 O({t^{-{3/2}}})
 \quad {\hbox {and}}\quad
 \int_{t}^{+\infty}
 {{{e}^{-i\int_{1}^{t'}\lambda_{+}}}\over {\lambda_{+}}}\,
 {{u}\over {(t')^{2}}}\, dt'=
 {{z_+}\over t}+
 O(t^{-2}).
$$
By replacing these identities into (\ref{cl26}) and taking into
account (\ref{cl27''}), it follows that
$$
 {e}^{-i \int_{1}^{t}\lambda_{+}}\twp(t)= z_+\,
 \left( 1- {3\over 4}\, {{\tg_+}\over t}\, i  \right)
 +O\left( {1\over {t^{3/2}} } \right),
 \qquad {\hbox {as}} \qquad
 t\rightarrow +\infty,
$$
so that, writing  $z_+= b_+{e}^{i a_+}$, for certain $a_+\in
[0,2\pi)$ and $b_+\geq 0$, we conclude
\begin{equation}\label{cl28}
 \twp(t)= b_+\,
 {e}^{i(\int_{1}^{t'}\lambda_{+}+a_+)}\,
 \left( 1- {3\over 4}\, {{\tg_+}\over t}\, i  \right)
 +O\left( {1\over {t^{3/2}} } \right),
 \qquad {\hbox {as}} \qquad
 t\rightarrow +\infty.
\end{equation}
 In the sequel, the constant $a_+\in [0,2\pi)$ may change its value
 at each occurrence.

 Coming back to (\ref{cl27'}), from (\ref{cl28}) we obtain the
 following asymptotics for $u(t)$ and $u'(t)$, as
 $t\rightarrow +\infty$:
\begin{eqnarray}\label{cl29}
 u(t)&=&
 b_+\cos \phi(t)+ {b_+\over 4}\left( 3\tg_+ -{{\gamma^2_+}\over 2} \right)
 {{\sin \phi(t)}\over {t}} + O\left( {1\over {t^{3/2}}} \right),
            \\
 u'(t)&=&
 -\, b_+\sin \phi(t) + {b_+\over 4}
 \left( 3\tg_+ - {{\gamma_+^2}\over 2} \right)
 {{\cos \phi(t)}\over {t}} + {b_+\over 2}\, \gamma_+ {\sin \phi(t)\over t}
 + O\left( {1\over t^{3/2}}  \right),
     \nonumber
\end{eqnarray}
with $a_{+}\in[0, 2\pi)$ and $\phi(t)=a_+ + t -{\gamma_+\over
2}\log t$, $t>0$.

Finally, recall that $y(s)=u(s^2/4)$ and $h(s)=v(s^2/4)$ (see
(\ref{cl14})).
Also, from (\ref{cl16'}) and (\ref{cl20}), it follows that
\begin{equation}\label{cl29''}
 v= {{u'}\over 2}- {{\tg_+}\over {2\rt}} + \gamma_+\,
 {{u'}\over {2t}}+ O\left( {1\over {t^{3/2}}}\right).
\end{equation}
Therefore, from the asymptotics related to $u$ and $u'$ and these
two observations, we get that
\begin{eqnarray}\label{cl30}
 y(s)
 &=&
 b_+\cos \phi(s^2/4)
 + b_+\,\left( 3\tg_+-{{\gamma^2_+}\over 2} \right)\,
 {{\sin \phi(s^2/4)}\over {s^2}}
 + O\left( {1\over s^3}  \right),
      \nonumber \\
 h(s)
 &=&
 -\,{b_+\over 2}\sin \phi(s^2/4) -\, {\tg_+\over s}
 +\, {b_+\over 2} \,
 \left( 3\tg_+ - {{\gamma^2_+}\over {2}} \right)
 {{\cos \phi(s^2/4)}\over {s^2}}
    \nonumber  \\
 &-&
 \gamma_+ \, b_+ {{\sin \phi (s^2/4)}\over {s^2}}
 +O \left( {1\over s^3}  \right)
 \qquad {\hbox {as}}\qquad
 s\rightarrow +\infty.
 \end{eqnarray}
 The asymptotics of $\bar f f'$ in {\it {iv)}} is now an immediate
 consequence of (\ref{cl30}), and the fact that $\bar ff'=y/2 + i
 h$ (see~(\ref{cl1})). Besides, in particular from (\ref{cl30}),
 $b_{+}=\lim_{s\rightarrow +\infty} |y+2ih|(s)=
 2|f|_{+\infty}|f'|_{+\infty}$,
 that is
 \begin{equation}\label{cl30''}
 b_{+}= 2\, |f|_{+\infty}|f'|_{+\infty}\geq 0
 \qquad {\hbox {and}}\qquad
 b^2_+=4\, |f|^2_{+\infty}\, (E(0)- (|f|^2_{+\infty}+\alpha)^2/4),
 \end{equation}
 using the conservation law (\ref{cl}).

Finally, we give the proof of {\it {v)}}. We start the study of
$f(s)$ in the case when $|f|_{+\infty}\neq 0$. To this end, by
using polar coordinates in the plane, we write $f(s)$ as
$$
 f(s)= \rho(s)\, e^{i\varphi(s)}.
$$
The asymptotics for $f(s)$ will follow from the ones concerning
with $\rho (s)$ and $\varphi (s)$.
Notice that $\rho^2=|f|^2$. Then, from (\ref{cl13}) and
(\ref{cl30}), we obtain 
\begin{eqnarray}\label{cl31}
 \rho(s)
 &=&
 |f|_{+\infty}
 \left\{
 1+ {b_+\over \fdm}\, {\sin\phi(s^2/4)\over s}+
 {1\over \fdm}\,
 \left(
  \tg_{+}- {b^2_+\over {4\fdm}}
 \right)
 \, {1\over s^2}
 \right.
   \nonumber \\
 & &
 \left.
 \qquad
 +{{b_+^2}\over {4\fcm}}\,
 {{\cos 2\phi(s^2/4)}\over {s^2}}
 + O\left( {1\over {|s|^3}} \right)
 \right\}.
\end{eqnarray}
Secondly, since $f(s)=\rho(s)\, {e}^{i\varphi(s)}$, $\Im m(\bar
ff')=\varphi'\rho^2$ and, from the asymptotics of $h(s)=\Im m(\bar
f f')$ in (\ref{cl30}), we obtain 
$$
 \varphi'\rho^2=
 -{b_+\over 2} \sin\phi(s^2/4) -{\tg_+\over s}
 +{b_+\over 2}\,
 \left(
   3\tg_+ - {\gamma_+^2\over 2}
 \right)\, {{\cos \phi(s^2/4)}\over {s^2}}
 -\gamma_+ b_+\, {{\sin \phi(s^4/2)}\over {s^2}}
 +O\left( {1\over {|s|^3}} \right).
$$
If we now use (\ref{cl31}) in the previous identity and integrate
the result from $1$ to $s\gg 1$, we obtain 
\begin{eqnarray}
 \label{cl33}
 \varphi(s)
 &=&
 c_+ -{1\over \fdm}
 \left(
 \tg_+ -{{b_+^2}\over {2\fdm}}
 \right)\, \log |s|
 + {{b_+}\over \fdm}\, {{\cos \phi(s^2/4)}\over {s}}
    \nonumber \\
 & &
 \qquad \qquad
 - {{b_+^2}\over {2\fcm}}\, {{\sin 2\phi(s^2/4)}\over {s^2}}
 +O\left( {1\over {|s|^3}}  \right),
\end{eqnarray}
with (see (\ref{cl30''}) and (\ref{cl10}))
\begin{equation}\label{cl33'}
 {1\over \fdm}\,
 \left(
 {{b_+^2}\over {2\fdm}}- \tg_+
 \right)= |f|^2_{+\infty}+\alpha.
\end{equation}
Finally, recall that $f(s)=\rho(s)\, {e}^{i\varphi(s)}$. Then,
after a few simplifications where we use the identities
$\phi(s^2/4)=a_++\pp$, $b_+=2|f|_{+\infty}|f'|_{+\infty}$
 and (\ref{cl33'}),
from (\ref{cl31}) and (\ref{cl33}) we easily get that
\begin{eqnarray}\label{cl34}
 f(s)
 &=&
  |f|_{+\infty}\, e^{ic_+}\, e^{i\ps} + 2i\, |f'|_{+\infty} \,
 {{e^{id_+}}\over {s}}\, e^{i\pt}
   \nonumber \\
 & &
 \qquad \qquad
 - |f|_{+\infty} (|f|^2_{+\infty} +\alpha)
 {{e^{ic_+}}\over {s^2}}\, e^{i\ps} + O\left( {1\over {|s|^3}} \right),
\end{eqnarray}
with $\ps= (\fdm +\alpha)\, \log |s|$, $ \pt=(\phi_2 -\phi_1)(s)
     =-(s^2/ 4) -(2|f|^2_{+\infty}+\alpha)\,\log|s|
$ and $d_+=c_+-a_+$. This concludes the proof of the asymptotic
behaviour of $f(s)$, when $|f|_{+\infty}\neq 0$.

We now consider the case $|f|_{+\infty}=0$. Then, $b_{+}=0$ (see
(\ref{cl33'})), and from the results in the part {\it {iv)}} it
follows that
\begin{equation}\label{cl34'}
 |f(s)|^2=2\,{{\tilde\gamma_{+}}\over{s^2}}+
 O\left( {{1}\over {s^3}}
 \right)
 \qquad s\rightarrow +\infty,
\end{equation}
 with ${\tilde\gamma_{+}}=2|f'|_{+\infty}$.

 Next, define $g=e^{-i\beta \log s}f$, where
 $\beta=|f|^{2}_{+\infty}+\alpha=\alpha$. Since $f$ solves
 (\ref{eqf}), we get that $g$ solves
 \begin{equation}\label{cl34''}
   \left[
   e^{i({{s^2}\over {4}}+2\beta\log s)} g'
   \right]'
   = e^{i({{s^2}\over {4}}+2\beta\log s)}\,G(s),
 \end{equation}
 with
 $$
   G(s)=-{g\over 2}\, (|g|^2-|g|^2_{+\infty})+
   {{\beta(\beta+i)}\over {s^2}}\, g=
   {g\over s^2}\, (\beta(\beta+i)-{\tilde\gamma}_{+})+
   O\left( {{1}\over {|s|^3}}  \right),
 $$
 bearing in mind(\ref{cl34'}).

 Simple computations give the asymptotic behaviour of $f(s)$ in
 the case $|f|_{+\infty}=0$.

In order to obtain the asymptotics of $f'(s)$, we observe that
$$
 |f|^2 f' = f(\bar f f')= f
 \left(
  {1\over 2}\,{{d|f|^2}\over {ds}}+ i\Im m (\bar f f')
 \right),
$$
so that
\begin{equation}\label{cl35}
 f'={f\over {|f|^2}}\,
 \left({1\over 2}\, {{d|f|^2}\over {ds}} + i\Im m (\bar ff')\right),
\end{equation}
whenever $|f|^2(s)\neq 0$.

Recall that $\rho(s)=|f(s)|$, $\phi(s^2/4)=a_+ +\pp$,
$b_+=2|f|_{+\infty}|f'|_{+\infty}$ and
$\tg_+=2|f'|^2_{+\infty}-|f|^2_{+\infty}(|f|^2_{+\infty}+\alpha)$.
Then, the behaviour related to $f'(s)$ in the case
$|f|_{+\infty}\neq 0$
easily follows from (\ref{cl35}), the asymptotis of $\bar f f'$
given in  the part {\it {(iv)}},
(\ref{cl31}) and (\ref{cl34}). The same argument is valid in the
case $|f|_{+\infty}=0$, by using (\ref{cl34'}) instead of
(\ref{cl31}). To this end, it is enough to observe that
$|f'|_{+\infty}\neq 0$, if $|f|_{+\infty}=0$, unless $f=0$ (see
Lemma~\ref{unicidad} in the following pages), so that
${\tilde\gamma}_{+}=2|f'|_{+\infty}\neq 0$.
This finishes the proof of the Theorem~\ref{yh}.

\end{proof}
\begin{remark}\label{reyh}
 The parts {\it (iii)} and {\it (iv)} also hold true for the
 quantities
 $$
  y=\displaystyle{{d|\Tn'|^2}\over {ds}}
  \qquad {\hbox {and}}\qquad
  h=-{1\over 2}\, \A \Tn\cdot \Tn',
 $$
 with $\Tn=\Gn'$ and $\Gn$ solving the equation (\ref{12}), that is
 \begin{equation}\label{reyh1}
 (\I+\A)\Gn -s\Gn'= 2 \Gn' \times \Gn'', \qquad |\Gn'|^2=1.
 \end{equation}
 In order to prove this remark, it is enough to observe that the
 identity (\ref{cl3}) is also satified because it just involves  
 $|f|^2=c^2$ and $h$. In fact,
 (\ref{cl3}) is a consequence of Lemma~\ref{for} and the identity 
 $h=-a(G_3-sT_3)/4$, which easily follows from the equation
 (\ref{reyh1}).

 Also, notice that the pair $(y,h)$ defined as above is a solution of
 the system (\ref{cl9}) with $E(0)=a^2/4$.
 To this end, we firstly observe that the derivation of (\ref{reyh1})
 gives that
 \begin{equation} \label{reyh2}
  \A \Tn = s\Tn' =2 \Tn\times \Tn'', \quad |\Tn|=1,
 \end{equation}
 from which it follows that 
 \begin{equation} \label{reyh3}
  y= \displaystyle{{d|\Tn'|^2}\over {ds}}
   =(\A \Tn \times\Tn)\cdot \Tn'.
 \end{equation}
 Then,
 $$
 h'= -{1\over 2}\, \A \Tn\cdot \Tn''
   = -{1\over 2}\, (s\Tn'+2\Tn\times \Tn'')\cdot \Tn''
   = - {s\over 4}\, y.
 $$
 Also,
 $$
 y'-sh= (\A \Tn\times \Tn)'\cdot \Tn'
      + (\A \Tn\times \Tn)\cdot \Tn''
      + {s\over 2}\, \A \Tn \cdot \Tn',
 $$
 where, from Lemma~\ref{for}, we get that
 $$
 (\A \Tn\times \Tn)'=((0,0,-a)+aT_3\Tn)'\cdot \Tn'
                    =-(|\Tn'|^2+\alpha) |\Tn'|^2
 \qquad {\hbox {and}}
 $$
 \begin{eqnarray*}
 (\A \Tn\times \Tn)\cdot \Tn''+ {s\over 2}\, \A \Tn\cdot \Tn'
 &=&
 -\A\Tn\cdot (\Tn''\times \Tn)+ {s\over 2}\, \A\Tn \cdot \Tn'
    \\
 &=&
 -{1\over 2}\, |\A\Tn|^2={a^2\over 2}\, (1-T_{3}^{2})=
 {1\over 2}\, (a^2-(|\Tn'|^2+\alpha)^2),
 \end{eqnarray*}
 so that
 $$
  y'-sh = {a^2\over 2}-{1\over 2}\, (3|\Tn'|^2 +\alpha)(|\Tn'|^2+\alpha).
 $$
 Finally, from Lemma~\ref{for}, observe that $y$ and $h$ are bounded functions.
\end{remark}
\begin{remark}\label{wunicidad}
 Using a fixed point argument in (\ref{cl26}), we obtain:
 Given $z_+\in\C$, $\gamma_{+}$ and $\tg_{+}$, there exists $t_0$
 sufficiently large and a unique $\twp\in{\mathcal{C}}(t\geq
 t_0)$, with $\sup_{t\geq t_0}|\twp(t)|\leq 2|z_+|$, solving the
 integral equation (\ref{cl26}). Moreover, $\twp(t)$ satisfies the
 following limiting condition
 $$
 \lim_{t\rightarrow +\infty}
 e^{-i(t-{\gamma_+\over 2} \log|t|)}
 \, \twp(t)=z_+.
 $$
To this end, consider $X$ to be the following set
$$
 X=\{
 w_1\in{\mathcal {C}}(t\geq t_0)\ / \
 ||w_1||_{X}=\sup_{t\geq t_0} |w_1(t)|\leq 2 |z_{+}|
 \},
$$
and define the operator $T w_1$ as the righthand side of
(\ref{cl26})-see also (\ref{cl18'})-(\ref{cl20})-. Then it is 
easy to prove that $T:X\longrightarrow X$ is a contraction on $X$.
\end{remark}
\begin{theorem}\label{yhreciproco}
\noindent With the same notation as in Theorem~\ref{yh}, given
complex numbers $|f|_{+\infty}e^{i\theta_1}$ and
$|f'|_{+\infty}e^{i\theta_2}$, with $|f|_{+\infty}$, 
$|f'|_{+\infty}\geq 0$ and $\theta_1$, $\theta_2\in [0,2\pi)$, 
there exists a unique  
solution  $f\in{\mathcal {\C}}^2(\R)$ of (\ref{eqf}) such that
 \begin{eqnarray}
 \label{ast0}
 & &
 \lim_{s\rightarrow +\infty} e^{-i\phi_2(s)}
 f(s)= |f|_{+\infty} e^{i\theta_1}
 \qquad \qquad {\hbox {and}}
    \nonumber\\[-1ex]
 & &
     \\[-1ex]
 & &
 \lim_{s\rightarrow +\infty} e^{-i\phi_3(s)}
 f'(s)= |f'|_{+\infty} e^{i\theta_2}.
    \nonumber
 \end{eqnarray}
 A similar result can be obtained as $s\rightarrow -\infty$.
\end{theorem}
\begin{proof}[Proof of Theorem~\ref{yhreciproco}]

\noindent {\bf {{Existence.}}} Given $\theta _1$, $\theta_2\in [0,2\pi)$,
$|f'|_{+\infty}$ and $|f|_{+\infty}\geq 0$, define the compact
set
$$
  D^{\alpha}=
  \{  (z_1,z_2)\in\C\times\C\, / \,
  |z_2|^2+{1\over 4} (|z_1|^2+\alpha)^2=
  |f'|^2_{+\infty}+{1\over 4} (|f|^2_{+\infty}+\alpha)^2
  \}.
$$
We continue to prove that for fixed $\alpha$,
$|f|_{+\infty}e^{i\theta_1}$, $|f'|^2_{+\infty}e^{i\theta_2}$ and
$s_0\geq 1$, there exists $(f_{s_0}(0), f'_{s_0}(0))\in
D^{\alpha}$ and $f$ solution of (\ref{eqf}) associated to this
initial data such that
\begin{eqnarray}
 \label{ast8}
 e^{-i(|f(s_0)|^{2} +\alpha)\log|s_0|} f(s_0)
 &=&
  |f|_{+\infty} e^{i\theta_1}
 \qquad \qquad {\hbox {and}}
    \nonumber\\[-1ex]
 & &
     \\[-1ex]
 e^{i( {{s_0^2}\over {4}}+(2|f(s_0)|^{2}+\alpha)\log|s_0|)}
 f'(s_0)
 &=&
 |f'|_{+\infty} e^{i\theta_2}.
    \nonumber
\end{eqnarray}
To this end, firstly define
$$
 f(s_0)= |f|_{+\infty}\, e^{i\theta_1} e^{i(|f|^2_{+\infty}+\alpha)\log|s_0|}
 \qquad {\hbox {and}}\qquad
 f'(s_0)=|f'|_{+\infty} \,e^{i\theta_2}
 e^{-i ({{s^2_0}\over 4}+(2|f|^2_{+\infty}+\alpha)\log|s_0|)}.
$$
Notice that, from the above identities, it follows that
$|f(s_0)|=|f|_{+\infty}$ and $|f'(s_0)|=|f'|_{+\infty}$. Hence
$(f(s_0), f'(s_0))\in D^{\alpha}$.

Secondly, the ODE's theory asserts us the existence of a local
solution of (\ref{eqf}), in particular for any $D^{\alpha}$.
Moreover, by taking into account the conservation law (\ref{cl}),
we observe that such this solution is globally well-defined. As a
by-product, we can consider $(f_{s_0}(0), f'_{s_0}(0))$, and from
(\ref{cl}) we conclude that $(f_{s_0}(0), f'_{s_0}(0))\in
D^{\alpha}$ (because $(f(s_0), f'(s_0))\in D^{\alpha}$).

Next, for $s\geq 1$, define the map $F(s)$ by
$$
 F(s)
 \left(
 \begin{array}{c}
 f(0) \\[2ex]
 f'(0)
 \end{array}
 \right)
 =
 \left(
 \begin{array}{c}
 e^{-i(|f(s)|^2+\alpha)\log|s|} f(s) \\[2ex]
 e^{i({s^2\over 4}+(2|f(s)|^2\alpha)\log|s|)} f'(s)
 \end{array}
 \right),
 \qquad
 \forall\, (f(0), f'(0))\in\C\times \C,
$$
where $f$ is a solution of (\ref{eqf}) associated to the initial
data $(f(0), f'(0))$.

Let $\{s_n\, :\, s_n\geq 1  \}_{n\geq 1}$ be a sequence such that
$s_n\rightarrow+\infty$, as $n\rightarrow +\infty$.
For any fixed $n\in\N$, we have proved the existence of $(f_n(0),
f'_n(0))\in D^{\alpha}$ and $f_n$ solution of (\ref{eqf})
associated to this initial data such that (see (\ref{ast8}))
\begin{equation} \label{ast9}
 F(s_n)
 \left(
 \begin{array}{c}
 f_n(0)\\[2ex]
 f'_n(0)
 \end{array}
 \right)
 =
 \left(
 \begin{array}{c}
 e^{-i(|f_n(s_n)|^2+\alpha)\log|s_n|} f_n(s_n)
   \\[2ex]
 e^{i( {{s_n^2}\over {4}}+(2|f_n(s_n)|^2+\alpha)\log|s_n|)} f'_n(s_n)
 \end{array}
 \right)
 =
 \left(
 \begin{array}{c}
 |f|_{+\infty}\ e^{i\theta_1}\\[2ex]
 |f'|_{+\infty} \, e^{i\theta_{2}}
 \end{array}
 \right).
\end{equation}
Notice that, because $D^{\alpha}$ is a compact set, there exists a
subsequence of $\{(f_n(0), f'_n(0)) \}_{n\in\N}$ such that
\begin{equation}\label{ast10}
 \lim_{k\rightarrow +\infty} f_{n_k}(0)= f(0)
 \qquad  {\hbox {and}} \qquad
 \lim_{k\rightarrow +\infty} f'_{n_k}(0)= f'(0),
 \qquad {\hbox {for some}}\qquad
 (f(0), f'(0))\in D^{\alpha}.
\end{equation}
Now, define $F(\infty)=\lim_{s\rightarrow +\infty}F(s)$, i.e.,
$$
 F(\infty)
 \left(
 \begin{array}{c}
 f(0)\\[2ex]
 f'(0)
 \end{array}
 \right)
 =
 \lim_{s\rightarrow +\infty}
 F(s)
 \left(
 \begin{array}{c}
 f(0)\\[2ex]
 f'(0)
 \end{array}
 \right)
 =
 \lim_{s\rightarrow +\infty}
  \left(
 \begin{array}{c}
 e^{-i(|f(s)|^2+\alpha)\log|s|}f(s)\\[2ex]
 e^{+i ( {{s^2}\over {4}}+(2|f(s)|^2+\alpha)\log|s| )}f'(s)
 \end{array}
 \right).
$$
Theorem~\ref{yh}, part {\it {iv)}} and {\it {v)}}, asserts us the
existence of the above limits, and therefore $F(\infty)$ is a
well-defined operator. Besides, $F(\infty)$ is continuous 
(see remark~\ref{reyh0}).

We now consider the initial data (\ref{ast10}). Then, by taking
into account the continuity property of $F(\infty)$ and the
identities in (\ref{ast9}), we obtain
\begin{eqnarray*}
 & &
 F(\infty)
 \left(
 \begin{array}{c}
  f(0)\\[2ex]
  f'(0)
 \end{array}
 \right)=
 \lim_{k\rightarrow +\infty} F(\infty)
 \left(
 \begin{array}{c}
  f_{n_k}(0)\\[2ex]
  f'_{n_k}(0)
 \end{array}
 \right)
      \\[2ex]
 & &
 \qquad\qquad
 =\lim_{k\rightarrow +\infty}
 \lim_{s\rightarrow +\infty}
 \left(
 \begin{array}{c}
   e^{-i(|f_{n_k}(s)|^2+\alpha)\log|s|}
   f_{n_k}(s)\\[2ex]
   e^{i(
   {s^2\over 4}+(2|f_{n_k}(s)|^2+\alpha)\log|s|
   )}
   f'_{n_k}(s)
 \end{array}
 \right)
 =
 \left(
 \begin{array}{c}
   |f|_{+\infty} \, e^{i\theta_1}\\[2ex]
   |f'|_{+\infty} e^{i\theta_2}
 \end{array}
 \right).
\end{eqnarray*}
The limiting conditions in (\ref{ast0}) are a consequence of the
definition of $F(\infty)$ and the above identity, by taking into
account that $|f(s)|^2-|f|^2_{+\infty}=o(1),$ as $s\rightarrow
+\infty$. This concludes the proof of the existence.

We will now continue to prove the uniqueness of such this
solution. We will need the following lemma.
\begin{lemma} \label{unicidad}
 Let $f_j$, $j=1,2$, two solutions of (\ref{eqf}) such that
 \begin{equation}\label{a1}
  |f_j|_{+\infty}=|f|_{+\infty},
  \qquad\qquad |f'_j|_{+\infty}= |f'|_{+\infty}
  \qquad {\hbox {and}} \qquad
 \end{equation}
 \begin{equation}\label{a2}
 \lim_{s\rightarrow +\infty} e^{i\phi_{1,j}(s)}
 (\bar f_j f'_j)(s)= |f'|_{+\infty}|f|_{+\infty}e^{ia_{+}},
 \end{equation}
 $$
 \qquad {\hbox {with}}\qquad
 \phi_{1,j}(s)= {s^2\over 4}-\gamma_{+,j}\log|s|
 \qquad {\hbox {and}}\qquad
 \gamma_{+,j}= -3|f_j|^2_{+\infty}-\alpha,
 $$
 for some $|f|_{+\infty},\ |f'|_{+\infty}\geq 0$, and
 $a_+\in[0, 2\pi)$.
 Then,
 $$
 {{d|f_1|^2}\over {ds}}(s)= {{d|f_2|^2}\over {ds}}(s)
 \qquad {\hbox {and}} \qquad
 \Im m(\bar f_1 f'_1)(s)= \Im m (\bar f_2 f'_2)(s)
 \qquad \forall\, s\geq s_0\gg 1,
 $$
 that is
 $\bar f_1 f_1'=\bar f_2 f_2',\ $ $\forall\, s\geq s_0\gg 1$.
 A similar result can be obtained at $s=-\infty$.
\end{lemma}
\begin{proof}
Assume $f_j$, $j=1,2$, are two solutions of (\ref{eqf}) satisfying
(\ref{a1}) and (\ref{a2}), and consider $h_j=\Im m(\bar f_jf'_j)$,
$y_j=d|f_j|^2/ds$, and the associated $w_{1,j}(t)$ defined through
the change of variables (\ref{cl14}) and (\ref{cl23}).

By using (\ref{cl14}), (\ref{cl29''}) together with (\ref{cl23})
and (\ref{cl24}), after some straightforward calculations one gets
that
\begin{eqnarray*}
 & &
 \lim_{s\rightarrow +\infty} e^{i \phi_{1,j}(s)}\,
 (\bar f_j f'_{j})(s)=
 \lim_{s\rightarrow +\infty} e^{i \phi_{1,j}(s)}\,
 \left( {y_j\over 2} + i h_j  \right) (s)=
   \nonumber \\[2ex]
 & &
 \qquad\qquad
 \lim_{t\rightarrow +\infty} e^{i \phi_{1,j}(2\sqrt t)}\,
 \left( {u_j\over 2} + i v_j  \right) (t)=
 \lim_{t\rightarrow +\infty} e^{i (t- {\gamma_{+,j}\over 2}\log |t|)}\,
 e^{-i \gamma_{+,j}\log 2}
 \overline{w_{1,j}}(t).
\end{eqnarray*}
Here, we have made use that $\lambda_{+,j}=
{\sqrt{1-(\gamma_{+,j}/t)}}= 1+o(1)$, as $t\rightarrow +\infty$ in
obtaining the last identity.

Therefore, from the previous identity and the limiting condition
(\ref{a2}), we get that
\begin{equation}
 \label{a3}
 \lim_{t\rightarrow +\infty}
 e^{-i(t-{\gamma_{+,j}\over 2}\log|t|)} w_{1,j}(t)=
 b_{+}\, e^{-ia_+} e^{-i\gamma_{+,j}\log 2}.
\end{equation}

On the other hand, from (\ref{a1}) it follows that
\begin{equation}
 \label{a4}
 \tilde\gamma_{+}=\tilde\gamma_{+,j},
 \qquad
 \gamma_{+}=\gamma_{+,j}
 \qquad {\hbox {and}}\quad
 \lambda_{+}= \lambda_{+,j}={\sqrt{1-(\gamma_+/ t)}},
 \qquad j=1,2,
\end{equation}
for some $\tilde\gamma_{+}$, $\gamma_{+}$ and $\lambda_{+}$
depending only on $|f|_{+\infty}$ and $|f'|_{+\infty}$
(independent of $j$).

Next, notice that because the r.h.s. in (\ref{cl26}) is
 $O(1/\sqrt t)$, as $t\rightarrow +\infty$ and (\ref{a4}), from
 (\ref{cl26}) and (\ref{a3}), it is easy to check that
 $$
    z_{+,j}= z_{+},
    \qquad \qquad j=1,2,
 $$
for some $z_+$ depending on the fixed values $|f|_{+\infty}$,
$|f'|_{+\infty}$, $\theta_1$ and $\theta_2$.

Previous analysis shows that $w_{1,j}$, $j=1,2$, are two solutions
of the integral equation in (\ref{cl26}), for the above values of
$z_+$, $\gamma_{+}$ and $\tilde\gamma_{+}$ and, in particular,
$|w_{1,j}(t)|\leq 2\,|z_{+}|$, when $t$ is large enough.
Then, from Remark~\ref{wunicidad}, it follows that
$w_{1,1}(t)=w_{1,2}(t),\ \forall\, t\geq t_0$
and $t_0$ large enough. Therefore, after undoing the change of
variables, we get that
 $$
   y_{1}(s)=y_{2}(s)
   \qquad {\hbox {and}}\qquad
   h_1(s)= h_2(s),
   \qquad \forall s {\hbox { large enough}}.
 $$
\end{proof}

\noindent {\bf {{Uniqueness.}}} Given $|f|_{+\infty}$,
$|f'|_{+\infty}\geq 0$ and $\theta_{1}$, $\theta_{2}\in [0,2\pi)$,
assume $f_j$, $j=1,2$, are two
solutions of (\ref{eqf}) satisfying (\ref{ast0}), that is
\begin{eqnarray}
 \label{f1}
 \lim_{s\rightarrow +\infty} e^{-i \phi_{2,j}(s)}\,
 f_j(s)= |f|_{+\infty} e^{i\theta_1}
 \quad \quad {\hbox {and}}
    \\[2ex]
 \label{f2}
 \lim_{s\rightarrow +\infty} e^{-i \phi_{3,j}(s)}\,
 f'_{j}(s)= |f'|_{+\infty} e^{i\theta_2}.
\end{eqnarray}
\noindent {\underline {\it{Case $|f|_{+\infty}\neq 0$}}}:
Define  $h_j=\I m (\bar f_{j}f'_{j})$ and $y_j=d|f_j|^2/ds$. Then
$h_{j}$ and $y_{j}$ (see (\ref{cl1})) are solutions of the 
system of ODE's in (\ref{cl9}), and therefore the associated 
$w_{1,j}$, $j=1,2$, defined through the change of variables 
(\ref{cl14}) and (\ref{cl23}), is a
solution of the integral equation (\ref{cl26}), with
$z_{+,j}\in\C$, $\lambda_{+,j}=\sqrt{1-\gamma_{+,j}/t}$ and
$\tilde\gamma_{+,j}$ and $\gamma_{+,j}$ given by (\ref{cl10}) and
(\ref{cl12}).

Now, from (\ref{f1}) and (\ref{f2}), we firstly observe that
\begin{equation}
 \label{ast1}
 |f|_{+\infty}=|f_j|_{+\infty}
 \qquad {\hbox {and}}\qquad
 |f'|_{+\infty}= |f'_j|_{+\infty},
 \qquad j=1,2.
\end{equation}
Besides, from (\ref{f1}) and (\ref{f2}), we obtain 
\begin{equation}
 \label{ast3}
 \lim_{s\rightarrow +\infty}
 e^{i \phi_{1,j}(s)}
 (\bar f_j f'_j)(s)=
 |f|_{+\infty}|f'|_{+\infty} e^{i(\theta_2-\theta_1)},
\end{equation}
where recall that $\phi_{1,j}(s)=(\phi_{2,j}-\phi_{3,j})(s)=
(s^2/4)- \gamma_{+,j}\log|s|$.

Then, Lemma~\ref{unicidad} yields that
 \begin{equation}
  \label{ast6}
  {{d|f_1|^2}\over {ds}}(s)={{d|f_2|^2}\over {ds}}(s)
  \qquad {\hbox {and}} \qquad
  \I m(\bar{f_1} f'_1)(s)=  \I m(\bar{f_2} f'_2)(s),
  \qquad s {\hbox { large enough.}}
 \end{equation}
We will continue to prove that $f_1(s)=f_2(s)$. To this end, we
describe $f_j(s), \ j=1,2$ as
$$
 f_j(s)=\rho_j (s) e^{i\varphi_j(s)},\qquad \rho_j(s)\geq 0
 \qquad {\hbox {and}}\qquad
 \varphi_j(s)\in\R,
$$
so that,
$$
  |f_j|=\rho_j
  \qquad {\hbox {and}}\qquad
  h_{j}=\I m (\bar f_j f'_j)=\varphi'_j\, \rho_j,
  \qquad j=1,2.
$$
From the first identity in (\ref{ast6}) one gets that
$|f_1(s)|=|f_2(s)|\neq 0,\ \forall\, s \gg s_1$ and $s_1$ large
enough, because $|f_1|_{+\infty}=|f_2|_{+\infty}=|f|_{+\infty}\neq
0$ (see (\ref{ast1})). Then,
$$
   \rho_1(s)=\rho_2(s)\neq 0,\qquad \forall\, s\geq s_1.
$$
Now, recall that $h_1=h_2$ and $h_j=\varphi'_{j}\,\rho_j$ , so
that
$$
  \varphi_1(s)=\varphi_2(s)+\theta_{0},\qquad \forall\, s\geq s_1\gg 1,
$$
for some $\theta_0\in \R$, because $\rho_1(s)=\rho_2(s)\neq 0,
\forall\, s\geq s_1$. As a consequence,
 $$
  f_1(s)=f_2(s) \, e^{i\theta_0},
  \qquad \forall\, s\geq s_1.
 $$
Finally, by taking into account the latter identity in the
limiting condition (\ref{f1}), it is easy to check that $\theta_0=
2k\pi$, for some $k\in\Z$. Therefore, we conclude that
$$
 f_1(s)=f_2(s),\qquad \forall\, s\geq s_1
 \qquad {\hbox {c.q.d.}}
$$

\noindent {\underline {\it {Case $|f|_{+\infty}=0$}}}:
Define $g_{j}=e^{-i\alpha\log s}\, f_{j}$, for $j=1$ and $2$, 
and $\phi(s)=(s^2/4)+\alpha \log |s|$.

On the one hand, by using (\ref{f1}) and (\ref{f2}) we obtain 
$$
 |f_j|_{+\infty}=0, \qquad
 \lim_{s\rightarrow +\infty} g_{j}(s)=0
 \qquad {\hbox {and}}\qquad
 \lim_{s\rightarrow +\infty} e^{i\phi(s)} g'_{j}(s)=
 |f'|_{+\infty}\, e^{i\theta_2},
$$
so that, from (\ref{cl34''}), it is easy to see that $g_j$ are
solutions of the following integral equation:
\begin{equation}\label{ast7}
 g(s)=z_2^+ - z_1^+\int_{s}^{+\infty} e^{-i\phi(\eta)}\, d\eta
     +\int_{s}^{+\infty} e^{-i\phi(\eta)}
       \int_{\eta}^{+\infty} e^{i\phi(s')}\, ds' \, d\eta,
\end{equation}
where $z_2^+=0$, $z_1^+=|f'|_{+\infty}\, e^{i\theta_2}$, and
$$
 G(s)=-{g\over 2}\, (|g|^2(s)-|g|^2_{+\infty})
      + {{\alpha(\alpha +i)}\over {s^2}}\, g.
$$
Also, since $f_j$ solves (\ref{eqf}) and $|f_j|_{+\infty}=0$, 
$|g|^2-|g|^2_{+\infty}=O(1/|s|^2)$, as $s\rightarrow +\infty$ (see
{\it iv)} in Theorem \ref{yh}).

On the other hand, a fixed point argument proves that: there exists
$s_0$ sufficiently large such that \ref{ast7} has an unique solution
in the space 
$$
 X=\{ g\in{\mathcal {C}}([s_0,+\infty)),\ 
 \lim_{s\rightarrow +\infty} g(s)=0\ {\hbox { and }}\ 
 |g|^2-|g|^2_{+\infty}=O(1/ |s|^2)   \}.
$$
As a by-product, we obtain  $g_1(s)=g_{2}(s)$, $\forall\, s\geq
s_0$. Thus $f_1(s)=f_2(s)$, $s\geq s_0$.
\end{proof}
\begin{remark}
 \label{reselfNLS}
 If $f$ solves (\ref{eqf}) and $|f|_{+\infty}=|f|_{-\infty}=0$, then $f(s)=0$.

 Indeed, let $f$ be a solution of (\ref{eqf}) such that
$|f|_{\pm\infty}=0$. Then,  the asymptotics in
 {\it (iv)} of Theorem~\ref{selfNLS} imply that
 \begin{equation}
  \label{reself1}
  |f|^2(s)= {{2}\over {s^2}}\, {\tilde\gamma_{\pm}} +
  O\left( {1\over{|s|^2}}\right),
  \qquad {\hbox{with}}\qquad
  \tilde\gamma_{\pm}=2|f'|^2_{\pm\infty},
  \qquad s\rightarrow \pm\infty.
 \end{equation}
 Now define
 $$
   \psi(s,t)= {{1}\over {\sqrt{t}}}\, f(s/\sqrt{t})
   e^{i({s^2\over {4t}}-{\alpha\over 2}\log t)}
   \qquad t>0.
 $$
 It is easy to prove that $\psi(s,\cdot)$ is a solution of the cubic non-linear
 Schr\"odinger equation, that is, 
 $$
  i\psi_t +\psi_{ss} + {{\psi}\over 2}\,|\psi|^2=0,
  \qquad\qquad (NLS)
 $$
 and, from (\ref{reself1}),  $\psi(s,\cdot)\in L^2(\R)$.

 On the other hand, the conservation of the $L^2$-norm of the
 solutions of (NLS) yields
 \begin{eqnarray*}
  +\infty
  &>&
  \lim_{t\downarrow 0}
       ||\psi(\cdot, t)||^2_{L^2(\R)}=
  \lim_{t\downarrow 0}{1\over t}
   \int_{\R}|
   {1\over {\sqrt t}} f(s/{\sqrt{t}})
   e^{i({{s^2}\over{4t}}-{\alpha\over 2}\log {t})}
   |^2\,ds
      \\[1ex]
  &=&
  \lim_{t\downarrow 0}{1\over \sqrt{t}}\, \int_{\R}|f(\eta)|^2\, d\eta
  \geq\lim_{t\downarrow 0}
  {1\over{\sqrt{t}}}\,
  \int_{|\eta|\geq 1} |f(\eta)|^2\, d\eta ,
 \end{eqnarray*}
 with $\int_{|\eta|\geq 1}|f(\eta)|^2\, d\eta<+\infty$, because of
 (\ref{reself1}). As a consequence, we get that $|f(\eta)|=0$
 a.e. $|\eta|\geq 1$, and using (\ref{reself1}) we conclude that
 $\tilde\gamma_{\pm}=2|f'|^2_{\pm\infty}=0$.

 The result in the remark now follows from Lemma~\ref{unicidad} and
 the fact that $|f|_{\pm\infty}=|f'|_{\pm\infty}=0$
 (notice also that $f=0$ is a solution of (\ref{eqf})).
\end{remark}
 \section{Further results}
 \label{fr}
 \subsection{Symmetric solutions}\
 \label{ss}

 \smallskip
 Here, we consider two particular cases of solutions of the
 equation
 \begin{equation}
  \label{ss1}
  \Gn''={1\over 2}\,(\A+\I)\Gn\times \Gn'.
 \end{equation}
 These solutions come from the symmetry properties of the above
 equation. For these special cases, we will able to obtain more
 specific properties that the ones obtained in Section~\ref{main}.

 We study what we will refer as to {\it {``Odd Case"}} and {\it {``Mixed
 Case"}}.
 \medskip

 \noindent {\it {I. {\underline {Odd Case}}}}:
 Let $\Gn_{\delta}$ be 
 a solution of (\ref{ss1}) with the initial conditions
 \begin{equation}
  \label{ss2}
  \Gn_{\delta}(0)=(0,0,0)
  \qquad {\hbox {and}} \qquad
  \Gn_{\delta}'(0)=(0,{\sqrt{1-\delta^2}},\delta),
  \quad\delta\in[-1,1].
 \end{equation}
 Then,
 $$
   \Gn_{\delta}(s)=-\Gn_{\delta}(-s).
 $$
 Indeed, notice that if $\Gn_{\delta}$ satisfies (\ref{ss1}) and
 (\ref{ss2}), then $\tilde\Gn(s)=-\Gn_{\delta}(-s)$ is also a
 solution of (\ref{ss1}),
 $$
   \tilde\Gn(0)=-\Gn_{\delta}(0)=(0,0,0)
   \qquad {\hbox {and}}\qquad
   \tilde\Gn'(0)=\Gn_{\delta}'(0).
 $$
 Therefore, from the uniqueness assumption, we conclude that
 $\Gn_{\delta}(s)=-\Gn_{\delta}(-s)$.

 For a fixed $a\in\R$, consider the map $T^{+}$ defined by
 $$
   (\Gn_{\delta}(0)={\mathbf {0}},
   \Gn_{\delta}'(0)=(0,{\sqrt{1-\delta^2}},\delta))
   \longrightarrow A^{+}_{3}(\delta),
 $$
 with $\An^{+}(\delta)=(A^{+}_{j}(\delta))^{3}_{j=1}$ the vector
 which prescribes the asymptotic behaviour of the solution
 $\Gn_{\delta}$ at $s=+\infty$.

 We have already proved the continuity of $T^{+}$ (see
 Proposition~\ref{convergencia}). On the other hand, observe that
 $\Gn_{1}(s)=(0,0,s)$ and $\Gn_{-1}(s)=(0,0,-s)$, so that
 $A^{+}_{3}(1)=1$ and $A^{+}_{3}(-1)=-1$, respectively. As a consequence, for
 any $z^{+}_{3}\in (-1,1)$, there exists $\delta\in(-1,1)$ such
 that $A^{+}_{3}(\delta)=z^{+}_{3}$.

 In particular, if $z^{+}_{3}=0$, there exists $\delta_{0}\in(-1,1)$
  such that $A^{+}_{\delta_0}(1)=0$. Thus, we conclude the
 existence of $\Gn_{\delta_{0}}$ solution of (\ref{ss1}) and
 (\ref{ss2}) such that  is asymptotically a {\it {plane spiral}},
 that is
 $$
   \Gn_{\delta_0}(s)\approx  e^{\A\log|s|}\An^{+}(\delta_0),
   \quad {\hbox{with}}\quad
   A^{+}_{3}(\delta_0)=0,
   \quad {\hbox {as}}\quad
   s\rightarrow \pm\infty
 $$
 (recall that $\Gn_{\delta_0}(s)=-\Gn_{\delta_0}(-s)$, so that
 $\An^{+}(\delta_0)=\An^{-}(\delta_0)$, and the asymptotics in
 Proposition~\ref{asym}).
 
 Notice that, due to the invariance of (\ref{ss1}) under rotations
 with respect to the OZ-axe, from previous remark it follows that
 we can ``generate" all the solutions of (\ref{ss1}) which
 asymptotically are plane spirals.
 
 \begin{figure}[h]
  \label{Odd1}
  \begin{center}
   \scalebox{0.5}{\includegraphics{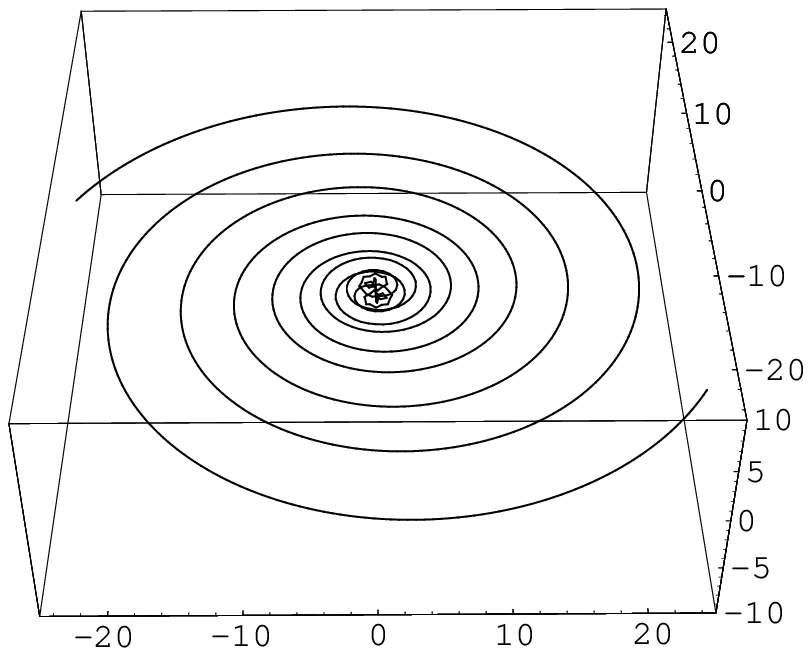}} 
   \hspace{2cm}
   \scalebox{0.5}{\includegraphics{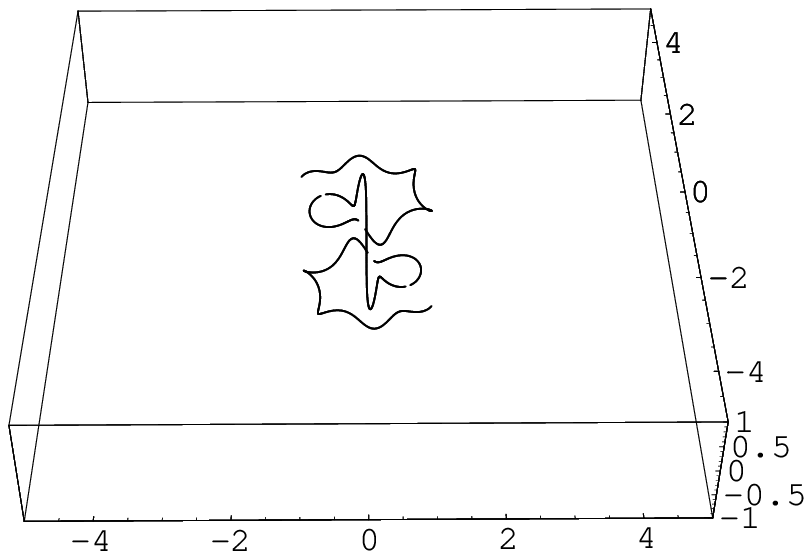}}  
  \end{center}
  \caption{$\Gn(0)={\mathbf {0}}$,
    $\Gn'(0)=(0,\sqrt{1-\delta^2},\delta)$, 
    $a=10$ and $\delta=0.956$.}
 \end{figure}
 \begin{figure}[h]
  \label{Odd3}
  \begin{center}
   \scalebox{0.55}{\includegraphics{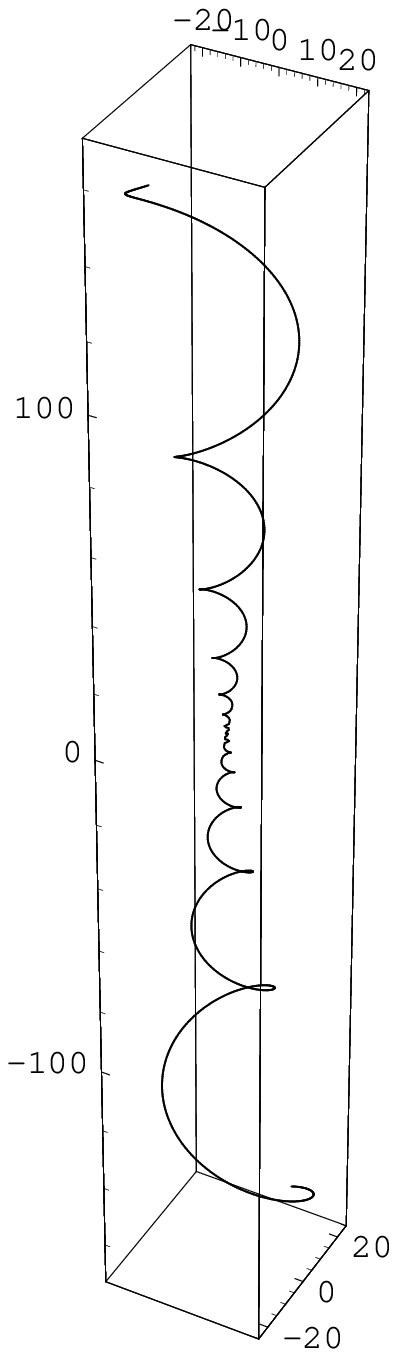}}  
   \hspace{3cm}
   \vspace{-.6cm}
   \scalebox{0.5}{\includegraphics{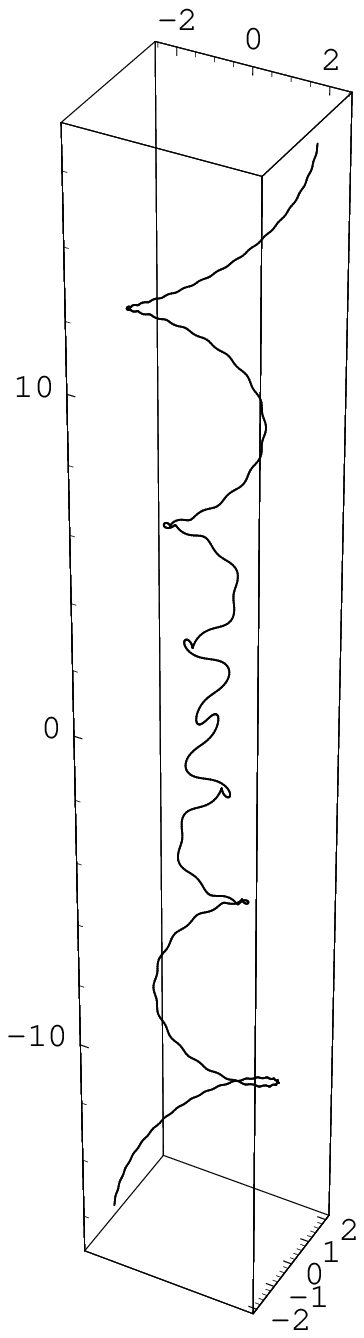}}  
  \end{center}
  \caption{$\Gn(0)={\mathbf {0}}$,
      $\Gn'(0)=(0,\sqrt{1-\delta^2},\delta)$, 
      $a=10$ and $\delta=-0.1$.}
 \end{figure}
 \begin{figure}[h]
  \label{Odd2}
  \begin{center}
   \scalebox{0.45} {\includegraphics{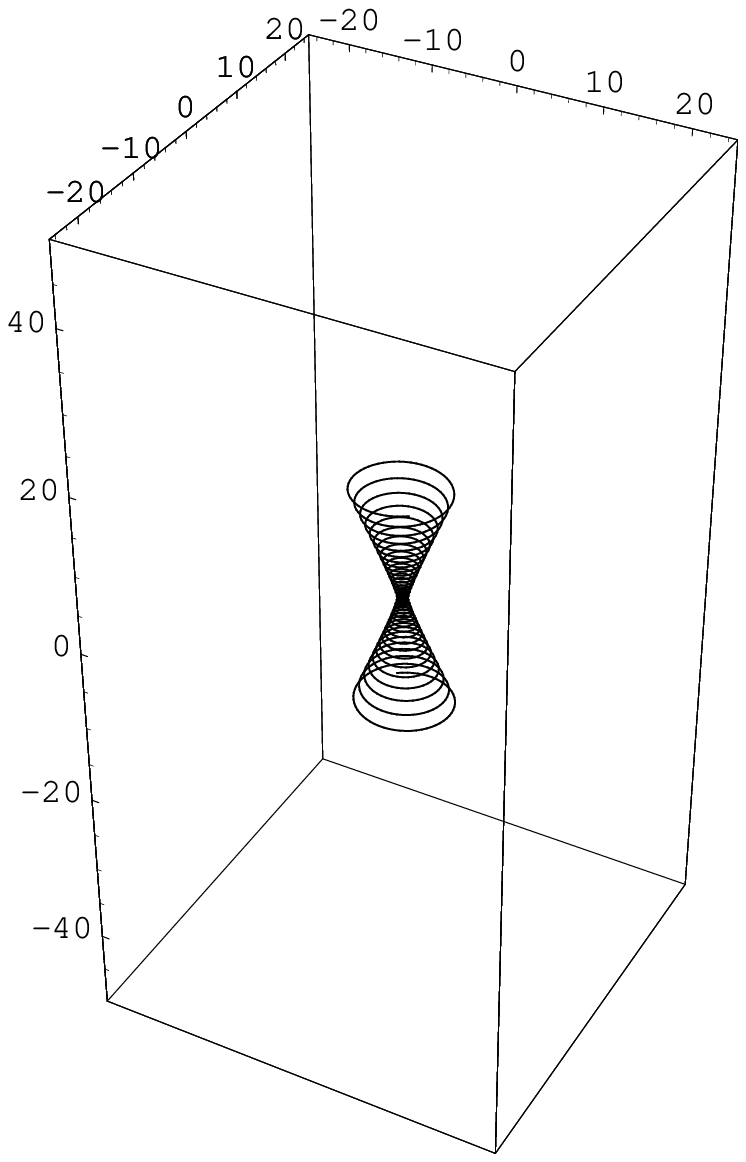}}  
    \hspace{2cm}
    \scalebox{0.6}{\includegraphics{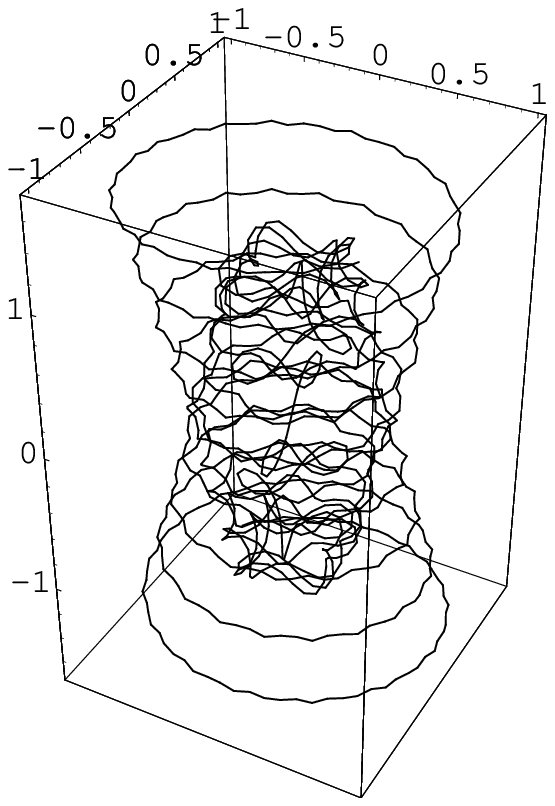}}  
  \end{center}
  \caption{$\Gn(0)={\mathbf {0}}$,
    $\Gn'(0)=(0,\sqrt{1-\delta^2},\delta)$, 
    $a=50$ and $\delta=0.9$.}
 \end{figure}
 In Figure~2, Figure~3 and Figure~4, we display the graphics of
 different solutions of (\ref{ss1}) associated to an initial data
 of the form (\ref{ss2}). The right-handside pictures represent
 the solution near the point $s=0$.
 \medskip

 \noindent {\it {II. {\underline {Mixed Case}}}}: Now assume that 
 $\Gn(s)$ is a solution of (\ref{ss1}) and define
 \begin{equation}
  \label{ss2'}
  \tilde\Gn(s)=(G_{1}(-s),G_{2}(-s), -G_{3}(-s))
 \end{equation}
 Then, a straightforward calculation gives that $\tilde\Gn(s)$ is also a
 solution of (\ref{ss1}).
 Moreover, in the particular case when $\Gn(s)$ also 
 satisfies the initial conditions
 \begin{equation}
  \label{ss3}
  \Gn(0)=(
   \displaystyle{
   {{2\,c_0}\over {\sqrt{1+a^2}}}}\, ,0,0
   )
  \qquad {\hbox{and}}\qquad
  \Gn'(0)= (0,0,1),
  \qquad (resp. \ \Gn'(0)=(0,0,-1)),
 \end{equation}
 with $c_0>0$, from the uniqueness assumption and the invariance of
 the equation (\ref{ss1}) under the transformation (\ref{ss2'}), it
 follows that $\Gn(s)$ satisfies
 \begin{equation}
  \label{ss4}
  \left\{
  \begin{array}{l}
  G_1(s)=G_1(-s)\\
  G_2(s)=G_2(-s) \\
  G_3(s)=-G_3(-s).
  \end{array}
  \right.
 \end{equation}
 In particular, from (\ref{ss4}) and the definitions of $\An^{\pm}$
 in (\ref{15''}), we get that
 $$
  A_{1}^{+}=-A_{1}^{-},\qquad
  A_{2}^{+}=-A_{2}^{-},
  \qquad {\hbox {and}}\qquad
  A_{3}^{+}=A_{3}^{-}.
 $$
 Also, from Lemma~\ref{for} and the initial conditions in
 (\ref{ss3}), it follows that
 \begin{equation}
  \label{ss6}
  c^2(s)=-aT_3(s) -\alpha
  \qquad {\hbox {with}}\qquad 
  \alpha=- a-c_0^2 
  \qquad (resp. \ \alpha=a-c_0^2).
 \end{equation}
 In \cite{GRV}, we studied the case when $a=0$. In the general
 frame when $a$ is not necessarily zero, only some of the results
 that were established there can be recovered. We continue to see
 some remarks on this respect.
 \begin{remark}
 If $a=0$ (i.e., $\A=0$), we proved in
 \cite{GRV} that given any $\theta\in(0,\pi]$, there exists $\Gn$
 solution of (\ref{ss1}) such that the angle between
 $\Bn^{+}=\An^{+}$ and $-\An^{-}=\Bn^{-}$ is precisely $\theta$.
 It would be interesting to prove that the same result holds for
 any $\A\in{\mathcal{M}}_{3\times3}$.
 \end{remark}
 \begin{remark} ({\it {Self-intersections}})

 \begin{itemize}
 \item[\it {(i)}]
 Given $a\in\R$, if $|c_0|>c_2$, with $c_2$
 large enough, then there exists $0<s\ll 1$ such that
 $$
  \Gn(s)=\Gn(-s)
  \quad (\Gn {\hbox{ has self-intersections}}).
 $$
 In fact,  arguing as in the proof of
 {\it{(ii)}} in \cite[Proposition 4]{GRV}, we prove that
 $\Gn(s)=\Gn(-s)$ if and only if $G_3(s)=0$, with $G_3(s)$ the
 solution of the following initial value problem:
 \begin{equation}
  \label{ss9}
  \left\{
  \begin{array}{l}
  {\displaystyle{
  G_3'''+\left(c^2+{s^2\over 4}  \right)G_3'-{s\over 4}\,G_3=
  -{a\over 2}\, (1-(G_3')^2)      }}
          \\[2ex]
  G_3(0)=0\quad G_3'(0)=\pm 1\quad G_3''(0)=0.

  \end{array}
  \right.
 \end{equation}
 By bearing in mind that $c^2=-aT_3(s)-\alpha$ (see
 (\ref{ss6})), similar arguments to those given in \cite{GRV}
 proves that, for values of the parameter $c_0$ large enough, the
 solution $G_3(s)$ of (\ref{ss9}) vanishes for some $0<s\ll 1$. 

 \item[\it {(ii)}]
 If $a\neq 0$, then Remark \ref{reformula} asserts
 that the solution $\G(s)$ does not self-intersect whenever $\alpha
 >0$, that is whenever $c_0^2\leq -a$ (resp. $c_0^2\leq a$)-see
 (\ref{ss6}).
 Also recall that in \cite[Proposition 4]{GRV} we proved that, if
 $a=0$ and $c_0$ is sufficiently small, then $\Gn(s)$ has no 
 self-intersections.
\end{itemize}

\end{remark}
%
 \begin{figure}[h]
  \begin{center}
   \scalebox{0.7}{\includegraphics{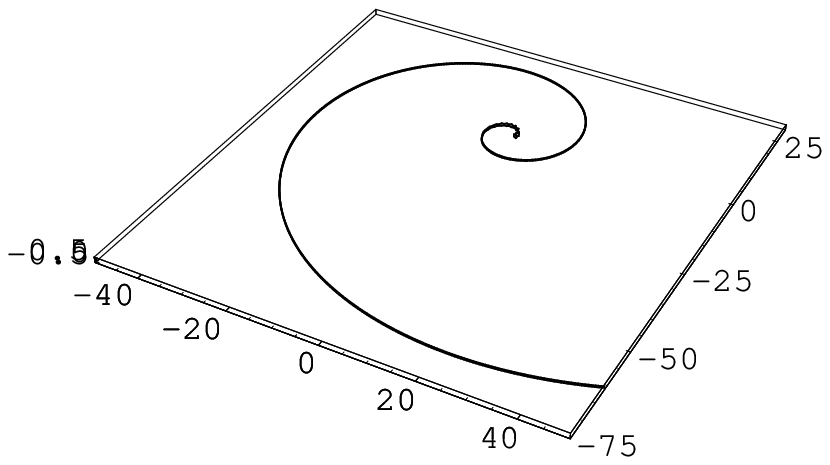}} 
   \hspace{2cm}
   \scalebox{0.65}{\includegraphics{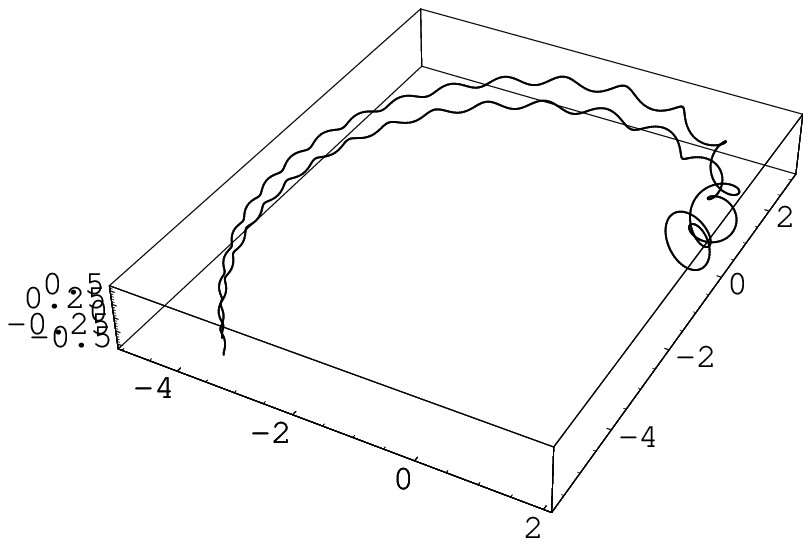}} 
  \end{center}
   \caption{$\Gn(0)=(2c_0/\sqrt{1+a^2},0,0)$, $\Gn'(0)=(0,0,1)$, 
             $a=3$, $c_0=1.8$.}
 \end{figure}
 \begin{figure}[h]
  \begin{center}
   \scalebox{0.6}{\includegraphics{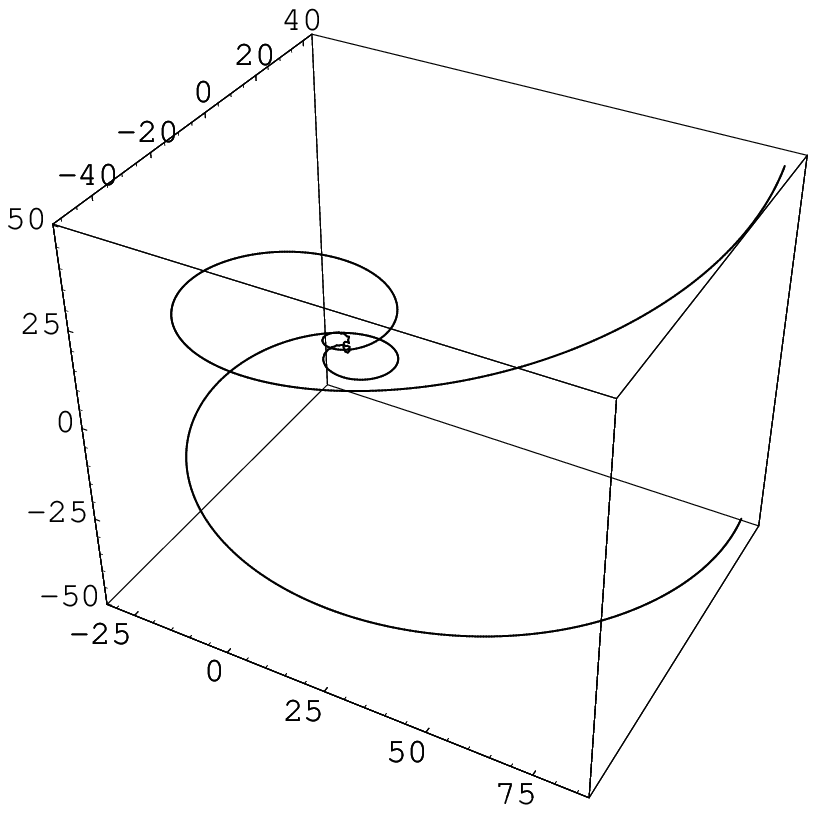}} 
   \hspace{2cm}
   \scalebox{0.6}{\includegraphics{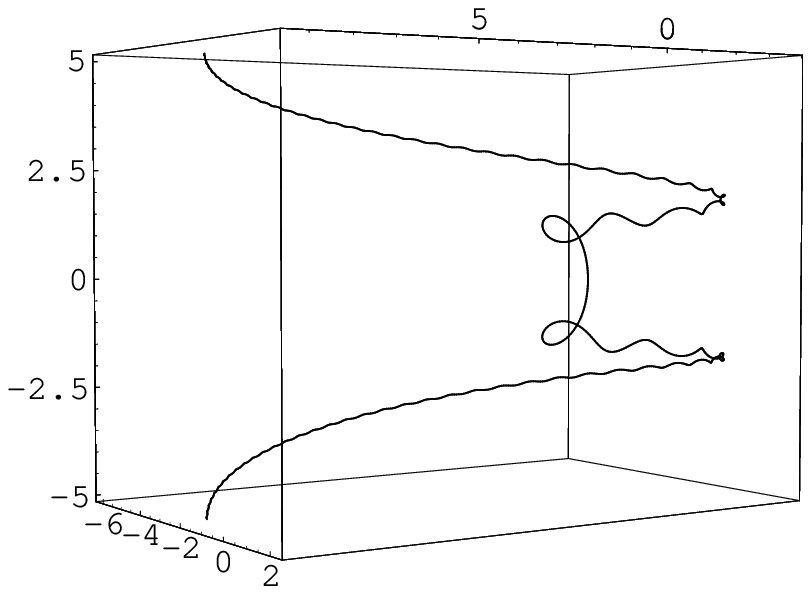}} 
  \end{center}
  \caption{$\Gn(0)=(2c_0/\sqrt{1+a^2},0,0)$, $\Gn'(0)=(0,0,1)$, 
           $a=3$, $c_0=0.4$.}
 \end{figure}
 \begin{figure}[h]
  \begin{center}
   \scalebox{0.44}{\includegraphics{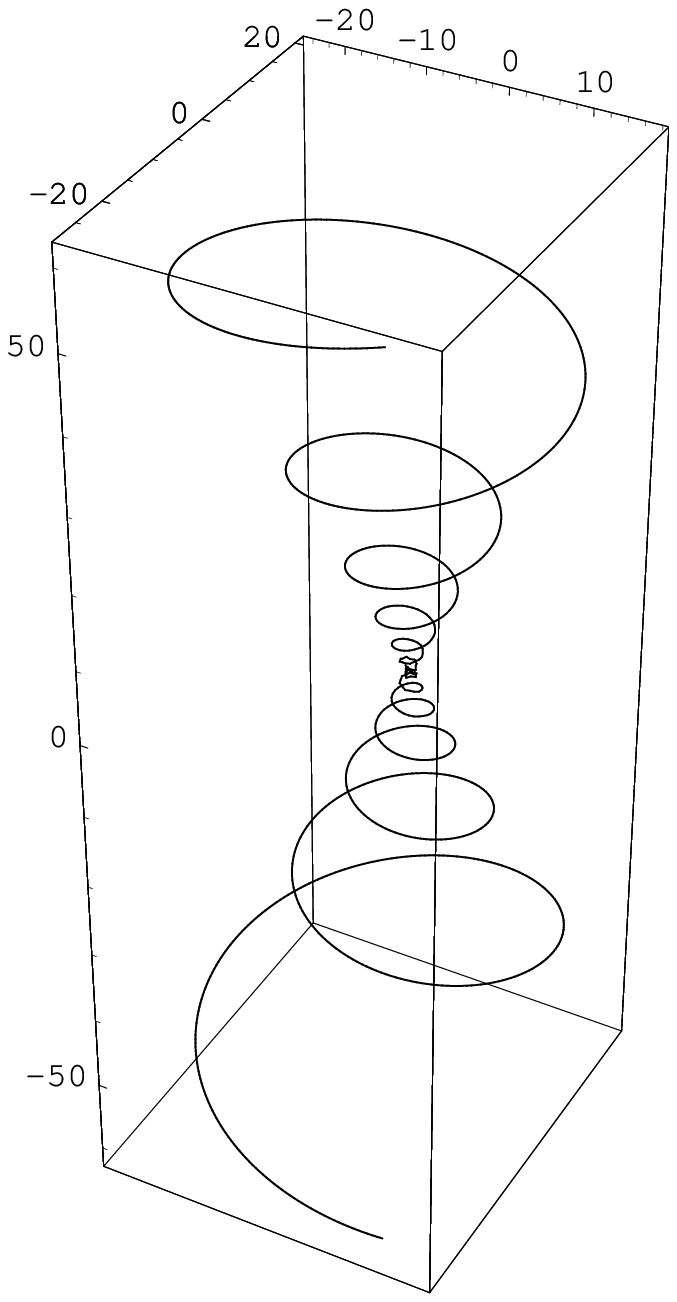}}
   \hspace{2cm}
   \scalebox{0.45}{\includegraphics{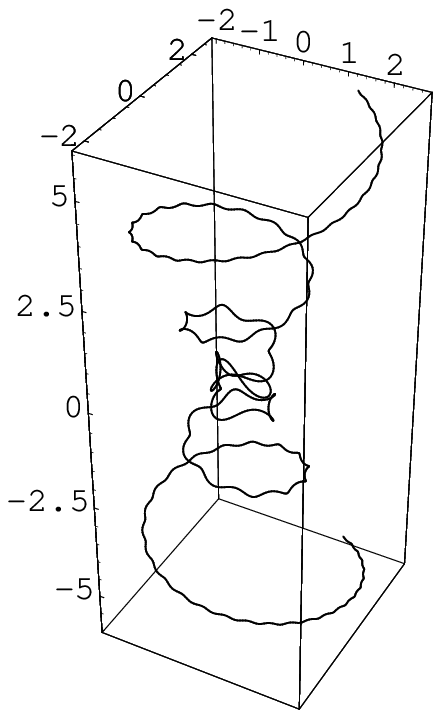}}
  \end{center}
  \caption{$\Gn(0)=(2c_0/\sqrt{1+a^2},0,0)$, $\Gn'(0)=(0,0,-1)$, 
           $a=10$, $c_0=4$.}
 \end{figure}
 Two different examples of solutions of 
 (\ref{ss1}) associated to a data of the type (\ref{ss3}) are plotted 
 in Figures~5, ~6 and ~7. Figure~5 represents 
 a plane spiral at infinity. Notice that in any of the latter figures, we can 
 clearly see the wavelike behaviour of the solutions.  As before, 
 in both figures the r.h.s. pictures are a zoom at the origin of the solution
 plotted in the l.h.s. picture.  All figures in this paper have been 
 generated using Mathematica~4.2.
 \subsection{On an ill-posedness result  
 for the initial value problem associated to 
 cubic Schr\"odinger equations}\
 \label{exun}

 \smallskip
 Here, we address the question of existence and
 uniqueness for the initial value problem (IVP) related
 to cubic Schr\"odinger equations with the principal
 value distribution as initial data:
 \begin{eqnarray}
  \label{exun1}
  \left\{
  \begin{array}{l}
   i\psi_t+\psi_{ss} +
   \displaystyle{
   {\psi\over 2}\,(|\psi|^2+{\alpha\over t})}
   =0,\qquad t>0, \quad s\in \R,\qquad
   \alpha\in\R
        \\[2ex]
   \psi(s,0)= c_1 \,p.v. (1/s),
   \qquad {\hbox {with}}\qquad c_1\in\C \setminus\{0\}.
  \end{array}
  \right.
 \end{eqnarray}
More precisely, the following ill-posedness result is true:
 \begin{proposition}
 \label{exunT2}
  Let $\alpha \geq 0$ and consider the IVP defined in
 (\ref{exun1}). Then, 
 either there is no weak solution $\psi$ for the IVP (\ref{exun1})
 in the class
$$
 t^{1\over 2} \psi, \quad 
 t^{3\over 2} |\psi|^2 \psi\in L^{\infty}([0,+\infty):{\mathcal{S}}'(\R))
 \qquad {\hbox {with}}\qquad
$$
$$
 \lim_{t\downarrow 0}\psi(s,0)= c_1\, p.v. (1/s),
 \quad {\hbox {in}}\quad {\mathcal {S}}'(\R),
 \quad c_1\in \C\setminus \{0\},
$$
 or there is more than one.
\end{proposition}
\begin{proof}[Proof of Proposition~\ref{exunT2}]
 The proof follows the argument in the proof of Theorem~1.5 in
 \cite{KPV}. Briefly, assuming the existence of a unique solution of
 the problem (\ref{exun1}), then it has to be a self-similar solution
 of the form
 $$
   \psi(s,t)={1\over {\sqrt{t}}}\, e^{i{{s^2}\over {4t}}} f(s/\sqrt{t}),
 $$
 with $f$ an odd solution of 
 $$
  f''+ i{s\over 2}\,f' +{f\over 2}\,(|f|^2+\alpha)=0,
  \qquad \alpha\in\R,
 $$
 On the other hand, since $\alpha\geq 0$, we get that
 $2|f|_{+\infty}+\alpha\neq 0$, and therefore, from the asymptotics in
 {\it (ii)} in this section, it follows that either there is no limit
 for the solution $\psi(s,t)$, as $t\rightarrow 0$, or the limit is
 identically zero, which contradicts the fact that $\psi(s,t)$ solves 
 (\ref{exun1}).
\end{proof}
Of special interest in Proposition~\ref{exunT2} is the case 
$\alpha=0$. If this is the case, the equation in (\ref{exun1}) reduces
to the commonly referred in the literature as to be the cubic
non-linear Schr\"odinger equation:
$$
 i\psi_t +\psi_{ss} + {{\psi}\over 2}\,|\psi|^2=0.
 \qquad\qquad\qquad  (NLS)
$$
In this setting, it is important to be mentioned that in
\cite{KPV} Kenig, Ponce and Vega proved an analogous ill-posedness 
result to the one in Proposition~\ref{exunT2} ($\alpha=0$) for the delta
distribution as initial datum. They left opened the same question
for the principal value.

 \subsection{On the ill-posedness of the Localized
  Induction Approximation}\
 \label{u}
 \smallskip

 We consider the following IVP:
 \begin{eqnarray}
 \label{u1}
 \left\{
 \begin{array}{l}
 \Xn_t= \Xn_s\times \Xn_{ss},\qquad s\in\R,\quad t>0,
    \\[1ex]
 \Xn(s,0)=\Xn_{0}(s).
 \end{array}
 \right.
 \end{eqnarray}
 Take $\widetilde{\Xn}(s,t)=\sqrt{t} \widetilde{\Gn}(s/{\sqrt{t}})$ with
 $\widetilde\Gn$ solving
 \begin{equation}
  \label{u2}
  {1\over 2}\, \widetilde\Gn -{s\over 2}\, \widetilde\Gn'
  =\tilde c\,\widetilde\bn
 \end{equation}
 and the initial conditions
 \begin{equation}
  \label{u3}
  \widetilde\Tn(0)={\mathbf{e_1}},\qquad
  \widetilde\nn(0)={\mathbf{e_2}},\qquad
  \widetilde\bn(0)={\mathbf{e_3}},
 \end{equation}
 where $\{ {\mathbf{e_1}},{\mathbf{e_2}},{\mathbf{e_3}} \}$ is the canonical
 basis in $\R^3$.

 Observe that evaluating (\ref{u2}) at $s=0$, we obtain 
 \begin{equation}
  \label{u3'}
  \widetilde\Gn(0)= 2\tilde c_0\,\widetilde\bn(0)=(0,0,2\tilde c_0),
 \end{equation}
 with $\tilde c_0=\tilde c(0)$ and $\tilde c(s)$ 
 the curvature associated to the curve $\widetilde\Gn(s)$.

 In \cite{GRV}, it was proved that $\widetilde\Xn$ solves (\ref{u1})
 with initial data
 \begin{equation}
  \label{u4}
  {\widetilde\Xn}_{0}(s)= s\left(\An^{+}\chi_{[0,+\infty)}(s)
  +\An^{-}\chi_{(-\infty,0]}(s) \right),
 \end{equation}
 for some $\An^{+}$ and $\An^{-}$ unitary vectors in $\R^3$,
 depending on the parameter $\tilde c_0$.

 Now define $\Gn(s)$ as
 \begin{equation}
 \label{u5}
 \Gn(s)=\widetilde\Gn(s)\chi_{[0,+\infty)}(s)
       +\rho\,\widetilde\Gn(s)\chi_{(-\infty,0]}(s)
 \qquad {\hbox {with}}\qquad
 \rho=
 \left(
 \begin{array}{ccc}
 -1&0&0  \\
 0&-1&0  \\
 0&0& 1  \\
 \end{array}
 \right).
 \end{equation}
 We claim that:
 \begin{equation}
  \label{u6}
  \Xn(s,t)={\sqrt{t}}\Gn(s/{\sqrt{t}})
 \end{equation}
 is also a solution of (\ref{u1}) with initial data
 \begin{equation}
 \label{u6''}
 \Xn_{0}(s)=\widetilde{\Xn}_{0}(s)\chi_{[0,+\infty)}(s)+
 \rho\,\widetilde{\Xn}_{0}(s)\chi_{(-\infty,0]}(s).
 \end{equation}
 We prove the claim by direct computation.

 Firstly, notice that (\ref{u6}) is a solution of LIA if and
 only if $\Gn$ satisfies
 \begin{equation}
  \label{u6'}
  {1\over 2}\, \Gn(s/{\sqrt{t}})- {s\over
  {2\sqrt{t}}}\,\Gn'(s/{\sqrt{t}})=
  (\Gn'\times\Gn'')(s/{\sqrt{t}}),
  \qquad t>0.
 \end{equation}
 Now, from (\ref{u5}), we obtain 
 $$
  \Gn'(s)=\widetilde\Gn'(s)\chi_{(0,+\infty)}(s)+
          \rho\,\widetilde\Gn'(s)\chi_{(-\infty,0)}(s)+
          [\widetilde\Gn(0)-\rho\,\widetilde\Gn(0)]\,\delta_{0},
 $$
 with $\widetilde\Gn(0)=\rho\,\widetilde\Gn(0)$, because of (\ref{u3'}).
 Therefore,
 \begin{equation}
  \label{u7}
  \Gn'(s)=\widetilde\Gn'(s)\chi_{(0,+\infty)}(s)+
          \rho\,\widetilde\Gn'(s)\chi_{(-\infty,0)}(s).
 \end{equation}
 Next, differentiating the above identity we obtain 
 $$
  \Gn''(s)=\widetilde\Gn''(s)\chi_{(0,+\infty)}(s)+
           \rho\,\widetilde\Gn''(s)\chi_{(-\infty,0)}(s)+
           [\widetilde\Gn'(0)-\rho\,\widetilde\Gn'(0)]\,\delta_{0},
 $$
 and from (\ref{u3})
 $$
  \widetilde\Gn'(0)={\mathbf {e_1}}= -\rho\,\widetilde\Gn'(0),
 $$
 so that
 \begin{equation}
  \label{u8}
  \Gn''(s)=\widetilde\Gn''(s)\chi_{(0,+\infty)}(s)+
           \rho\,\widetilde\Gn''(s)\chi_{(-\infty,0)}(s)+
           2{\mathbf {e_{1}}}\delta_{0}.
 \end{equation}
 We firstly compute $(\Gn'\times\Gn'')(s)$. From (\ref{u7}) and
 (\ref{u8}), we obtain 
 $$
  (\Gn'\times\Gn'')(s)=
   (\widetilde\Gn'\times\widetilde\Gn'')(s)\chi_{(0,+\infty)}(s)+
   \rho \,(\widetilde\Gn'\times\widetilde\Gn'')(s)\chi_{(-\infty,0)}(s)+
   2\Gn'(s)\times {\mathbf{e_1}}\delta_{0.}
 $$
 Notice that $\Gn'(s)\times {\mathbf {e_1}}$ can be extended to a
 continuous function up to $s=0$ because, from (\ref{u7}) and
 (\ref{u3}),
 $$
  \lim_{s\rightarrow 0^{\pm}}\Gn'(s)\times {\mathbf {e_1}}=0.
 $$
 As a by-product
 \begin{eqnarray*}
 (\G'\times\G'')(s)
  &=&
  (\widetilde\Gn'\times \widetilde\Gn'')\chi_{(0,+\infty)}(s)+
  \rho\, (\widetilde\Gn'\times \widetilde\Gn'')\chi_{(-\infty,0)}(s)
      \\
  &=&
   (\tilde c\,\widetilde\bn)(s)\chi_{(0+\infty)}(s)+
   \rho \,(\tilde c\,\widetilde\bn)(s)\chi_{(-\infty,0)}(s).
 \end{eqnarray*}
 Moreover, from the initial conditions in (\ref{u3}), it is easy
 to check that $\Gn'\times\Gn''$ can be extended continuously at
 $s=0$ and
 $$
  \lim_{s\rightarrow 0^{\pm}} (\Gn'\times\Gn'')(s)=
  \tilde c_0\,\widetilde\bn(0).
 $$
 In order to compute $(\Gn'- s \Gn)/2$ we use (\ref{u5}),
 (\ref{u7})and (\ref{u2}). Then, we obtain the following
 \begin {eqnarray*}
 {1\over 2}\, \Gn(s)-{s\over 2}\Gn'(s)
 &=&
  \left({1\over 2}\, \widetilde\Gn(s)
  -{s\over 2}\widetilde\Gn'(s)\right)\chi_{(0,+\infty)}(s)+
  \rho \,\left({1\over 2}\, \widetilde\Gn(s)
  -{s\over 2}\widetilde\Gn'(s)\right)\chi_{(0,+\infty)}(s)
     \\
 &=&
 (\tilde c\,\widetilde\bn)(s)\chi_{(0,+\infty)}(s) +
 \rho\,(\tilde c\,\widetilde\bn)(s)\chi_{(-\infty,0)}(s),
 \end{eqnarray*}
 and
 $$
  \lim_{s\rightarrow 0^{\pm}}
  \left( {1\over2}\,\Gn(s)-{s\over 2}\,\Gn'(s)
  \right)=\tilde c_0 \,\widetilde\bn(0).
 $$
 From previous identities we conclude (\ref{u6'}), and therefore
 $\X(s,t)=\sqrt{t}\Gn(s/\sqrt{t})$ satisfies
 \begin{equation}
  \label{u9}
  \Xn_{t}=\Xn_s\times\Xn_{ss},\qquad s\in\R\quad t>0.
 \end{equation}
 Notice that $\X(s,t)=\sqrt{t}\Gn(s/\sqrt{t})$ is continuous for
 all $s$ and $t>0$ (see (\ref{u5}) and (\ref{u3'})), and
 $\Xn_s(s,t)=\Gn'(s/{\sqrt{t}})$ is a real analytic function
 except at $s=0$, where it has a jump singularity.
 In fact, using (\ref{u7}) and (\ref{u3}) it is easy to check that
 $$
  \lim_{s\rightarrow 0^{+}} \Gn'(s/{\sqrt{t}})={\mathbf{e_1}},
  \qquad {\hbox {and}}\qquad
  \lim_{s\rightarrow 0^{-}} \Gn'(s/{\sqrt{t}})=-{\mathbf{e_1}}
 $$
 Then, here (\ref{u9}) is understood for
 $$
  \Xn(s,t)\in {\mathcal{C}}
  \left( (0,+\infty): {\hbox {Lip}}^{1}(\R) \right)
 $$
 and
 $$
   (\Xn_t)(s,\cdot)\in{\mathcal{C}}(\R),
   \qquad {\hbox {for all}}\qquad
   t>0.
 $$
 Finally, using the same arguments as in the proof of
 in \cite[Proposition~1]{GRV}, we observe that
 $$
  |\Xn(s,t)-\Xn_{0}(s)|\leq 2\tilde c_0{\sqrt t},
  \qquad \forall\, s\in\R.
 $$
 This concludes the proof of the claim.

 The claim asserts the existence of a (singular) solution of LIA
 associated to a given initial data $\Xn_{0}(s)$ as in
 (\ref{u6''}). Such this solution is different from the (regular)
 solution obtained in \cite{GRV} for the same data.
 The following proposition gathers the above obtained ill-posedness result:
\begin{proposition}
 \label{uT1}
 Given $(\An^{+},\An^{-})\in {\mathbb{S}}^{1}\times
 {\mathbb{S}}^{1}$ with $\An^{+}-\An^{-}\neq {\mathbf {0}}$, consider the IVP:
 \begin{eqnarray}
 \label{u10}
 \left\{
 \begin{array}{l}
 \Xn_t= \Xn_s\times \Xn_{ss},\qquad s\in\R,\quad t>0,
    \\[1ex]
 \Xn(s,0)=s\left(\An^{+}\chi_{[0,+\infty)}(s)
  +\An^{-}\chi_{(-\infty,0]}(s) \right).
 \end{array}
 \right.
 \end{eqnarray}
 Then, there exist more that one solution for the IVP (\ref{u10})
 in the class
 $$
   \Xn(s,t)\in {\mathcal{C}}
   \left( (0,+\infty): {\hbox {Lip}}^{1}(\R) \right),
 $$
 $$
  (\Xn)_{t}(s,\cdot)\in{\mathcal{C}}(\R)
  \quad {\hbox{for all}}\quad
  t>0.
 $$
\end{proposition}
%
 Figure~8 illustrates the vortex line evolution for
 same value of $c_0$ (i.e, $c_0=0.8$). Precisely, the
 subsequent plots correspond to vortex position at time
 $t=10^{-4}$, $0.1$, $1$, $1.5$, $2$ and $2.5$ 
 (c.f. \cite[Figure~1]{GRV}).
 \begin{figure}[h]
  \label{ill-LIA2}
  \begin{center}
   \scalebox{0.7}{\includegraphics{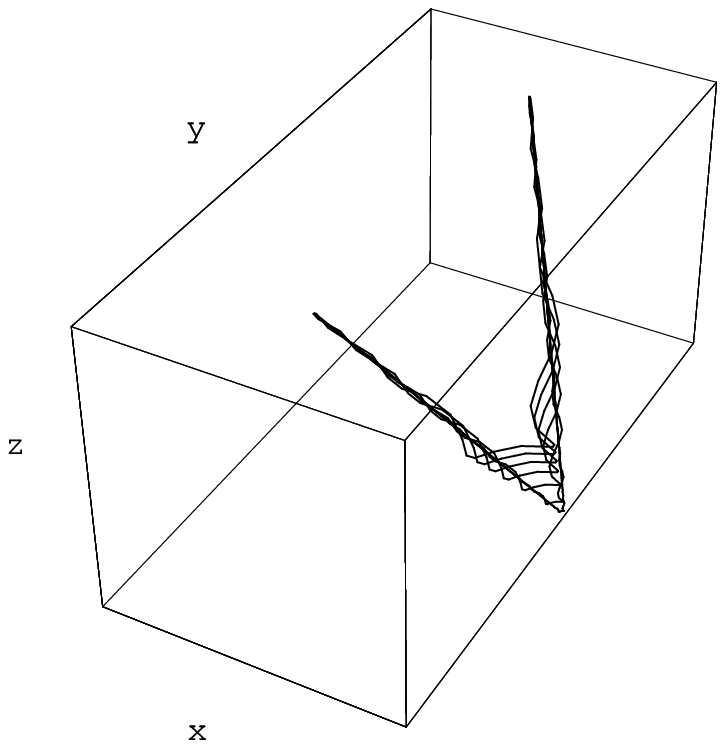}}
   \hspace{1cm}
   \scalebox{0.65}{\includegraphics{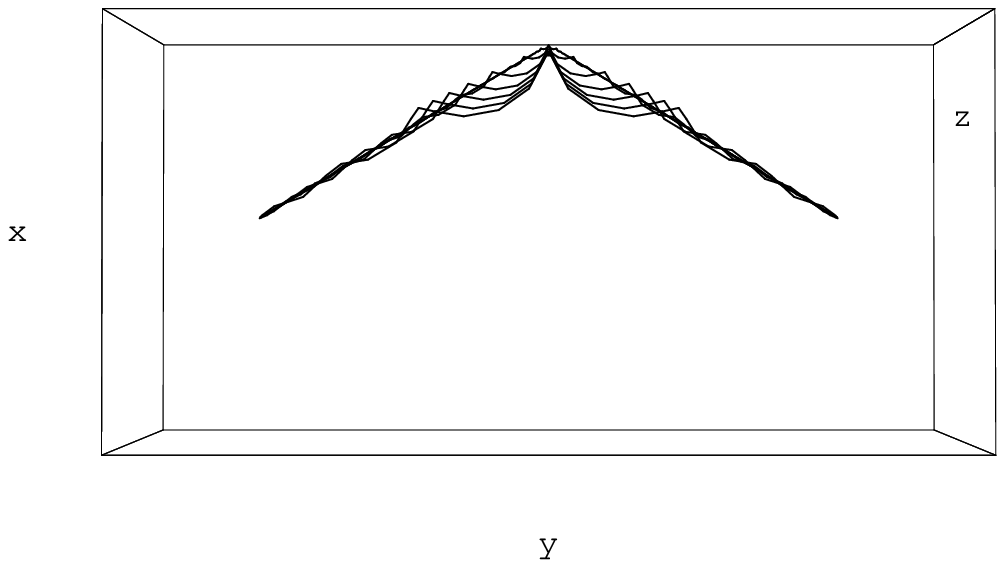}}
  \end{center}
  \caption{Singular vortex line evolution for $c_0=0.8$. Different views.}
 \end{figure}

\end{document}